\documentclass[11pt,twoside]{article}

\usepackage{mathrsfs}
\usepackage{amssymb,amsmath,amsfonts,amsthm,color,mathrsfs}
\usepackage[Symbol]{upgreek}
\usepackage{txfonts}
\usepackage[nottoc,notlot,notlof]{tocbibind}
\usepackage[active]{srcltx}
\usepackage[top=1in, bottom=1in, left=1.25in, right=1.25in]{geometry}

\theoremstyle{plain}

\newtheorem{thm}{Theorem}[section]
\newtheorem{lem}[thm]{Lemma}
\newtheorem{prop}[thm]{Proposition}
\newtheorem{cor}[thm]{Corollary}

\newtheorem{defn}[thm]{Definition}
\newtheorem{rem}[thm]{Remark}

\numberwithin{equation}{section}
\allowdisplaybreaks

\begin{document}
\title{Large Time Behavior and Convergence for the Camassa-Holm Equations with Fractional Laplacian Viscosity}
\author{ Zaihui Gan\thanks{ Corresponding author:
ganzaihui2008cn@tju.edu.cn(Zaihui Gan)}, Yong He\thanks{Yong He: heyong80@tju.edu.cn}, Linghui Meng\thanks{Linghui Meng: mlh1722604275@tju.edu.cn}
\\  {\small   Center for Applied Mathematics, Tianjin University,  Tianjin 300072, China}
 }
\date{}
\maketitle {
 {
\begin{quote}
{\bf Abstract.}
In this paper, we consider the $n$-dimensional ($n=2,3$) Camassa-Holm equations with fractional Laplacian viscosity in the whole space. In stark contrast to the Camassa-Holm equations without any nonlocal effect, to our best knowledge, little has been known on the large time behavior and convergence for the nonlocal equations under study. We first study the large time behavior of solutions. We then discuss the relation between the equations under consideration and the imcompressible Navier-Stokes equations with fractional Laplacian viscosity (INSF). The main difficulty to achieve them lies in the fractional Laplacian viscosity. Fortunately, by employing some properties of fractional Laplacian, in particular, the fractional Leibniz chain rule and the fractional Gagliardo-Nirenberg-Sobolev type estimates, the high and low frequency splitting method and the Fourier splitting method, we first establish the large time behavior concerning non-uniform decay and algebraic decay of solutions to the nonlocal equations under study.  In particular, under the critical case $s=\dfrac{n}{4}$, the nonlocal version of Ladyzhenskaya's inequality is skillfully used, and the smallness of initial data in several Sobolev spaces is required to gain the non-uniform decay and algebraic decay. On the other hand, by means of the fractional heat kernel estimates, we figure out the relation between the nonlocal equations under consideration and the equations (INSF). Specifically, we prove that the solution to the Camassa-Holm equations with nonlocal viscosity converges strongly as the filter parameter $\alpha\rightarrow~0$ to a solution of the equations (INSF).
\\[0.15cm]
{\bf Key Words.} Camassa-Holm equations; Fractional Laplacian viscosity; Large time behavior; Convergence
\\[0.15cm]
{\bf MSC(2000).} 35G25; 35K55\\
\end{quote}}}
%\begin{abstract}

%\end{abstract}

\section{Introduction}
\allowdisplaybreaks
 In this article, we investigate the following Camassa-Holm equations with fractional Laplacian viscosity in $\mathbb{R}^{n}~(n=2,3)$:
\begin{equation}\label{1.1}
\left\{
\begin{array}{ll}
\displaystyle \mathbf{v}_{t}+\mathbf{u}\cdot \nabla \mathbf{v}+\mathbf{v}\cdot \nabla \mathbf{u}^{T}+\nabla p=-\nu (-\Delta)^{\beta}\mathbf{v},~(t,x)\in \mathbb{R}^{+}\times \mathbb{R}^{n},
\\[0.2cm]
\displaystyle \mathbf{u}-\alpha^{2} \Delta \mathbf{u} =\mathbf{v},
\\[0.2cm]
\displaystyle\mbox{div}~ \mathbf{v}=\mbox{div}~ \mathbf{u}=0,
\end{array}
\right.
\end{equation}
with the initial condition
\begin{equation}\label{1.2}
\mathbf{v}(0,x)=\mathbf{v}_{0}(x),~\mathbf{u}(0,x)=\mathbf{u}_{0}(x),~x\in\mathbb{R}^{n}.
 \end{equation}
 Here, $\mathbf{v},~\mathbf{u}$ denotes the fluid velocity field and the filtered fluid velocity, respectively, and $p$ the scalar pressure. $\alpha$ is a length scale parameter representing the width of the filter, and $\nu>0$ is the viscosity coefficient which is fixed in our discussions. In particular, the divergence free condition $\mbox{div}~ \mathbf{v}=0$ indicates the imcompressibility of the fluid,  $(-\Delta)^{\beta}$ denotes the fractional power of the Laplacian in $\mathbb{R}^{n}$,  $\displaystyle\frac{n}{4}\leq \beta< 1$ and  $ n=2,3$. Recall that the Camassa-Holm equations with Laplacian viscosity (equations \eqref{1.1} with $\beta=1$) read
 \begin{equation}\label{1.3}
\left\{
\begin{array}{ll}
\mathbf{v}_{t}+\mathbf{u}\cdot \nabla \mathbf{v}+\mathbf{v}\cdot \nabla \mathbf{u}^{T}+\nabla p= \nu  \Delta \mathbf{v},~(t,x)\in \mathbb{R}^{+}\times \mathbb{R}^{n},
\\[0.2cm]
\mathbf{u}-\alpha^{2} \Delta \mathbf{u} =\mathbf{v},
\\[0.2cm]
\mbox{div}~ \mathbf{v}=0.
 \end{array}
\right.
\end{equation}
 As it is well-known that the system \eqref{1.3} rose from work on shallow water equations \cite{Camassa-1}. Specifically, it was introduced in \cite{Holm-Marsden} as a natural mathematical generalization of the integrable inviscid one-dimensional Camassa-Holm equation discovered in \cite{Camassa-1} through a variational formulation and with a lagrangian averaging. It could be used as a closure model for the mean effects of subgrid excitations, and be also viewed as a filtered Navier-Stokes equations with the parameter $\alpha$ in the filter, which obeys a modified Kelvin circulation theorem along filtered velocities \cite{Holm-Marsden}. Numerical examples that seem to justify this intuition were reported in \cite{Chen-holm}. Formally, the system \eqref{1.3} reduces to the imcompressible Navier-Stokes equations as $\alpha\rightarrow 0$:
\begin{equation}\label{1.4}
\left\{
\begin{array}{ll}
\mathbf{v}_{t}+\mathbf{v}\cdot \nabla \mathbf{v} +\nabla p=\nu \Delta \mathbf{v},
\\[0.2cm]
\mbox{div}~ \mathbf{v}=0.
\end{array}
\right.
 \end{equation} \\
 \indent For the fractional Laplacian in the whole space, there are several different ways to define it \cite{Caffarelli-1,Landkof,Nezza}. For example, for a function $f\in \mathscr{C}$, the integral fractional Laplacian $(-\Delta)^{\beta}$ at the point $x$ can be defined as
\begin{equation}\label{1.5}
\left.
\begin{array}{ll}
\displaystyle I_{\beta}f(x)&\triangleq(-\Delta)^{\beta}f(x):= C_{n,\beta}~\mbox{P.V.}\int_{\mathbb{R}^{n} }\frac{f(x)-f(\xi)}{|x-\xi|^{n+2\beta}}d\xi,
\\[0.20cm]
 &\displaystyle := C_{n,\beta}~\mbox{P.V.}\int_{\mathbb{R}^{n} }\frac{f(x)-f(\xi)}{|x-\xi|^{n+2\beta}}d\xi
\\[0.35cm]
&\displaystyle := C_{n,\beta}\lim\limits_{\varepsilon\rightarrow 0^+}\int_{|\xi|>\varepsilon }\frac{f(x+\xi)-f(x)}{|\xi|^{n+2\beta}}d\xi,
 \end{array}
\right.
\end{equation}
or
 \begin{equation}\label{1.6}
 \displaystyle I_{\beta}f(x)\triangleq (-\Delta)^{\beta}f(x):=\frac{C_{n,\beta}}{2}\int_{\mathbb{R}^{n} }\frac{2f(x)-f(x+y)-f(x-y)}{|y|^{n+2\beta}}dy,\end{equation}
where
 the parameter $\beta$ is a real number with $0<\beta<1$, $\mbox{P.V.}$ is a commonly used abbreviation for "in the principle value sense" (as defined by the latter equation), and $C_{n,\beta}$ is some normalization constant depending only on $n$ and $\beta$, precisely given by
  \begin{equation}\label{1.7}
 \displaystyle C_{n,\beta}=\left(\int_{\mathbb{R}^{n} }\frac{1-\cos(\zeta_{1})}{|\zeta|^{n+2\beta}}d \zeta\right)^{-1}.\end{equation}
\indent Before going further, we collect some facts on the fractional Sobolev spaces $W^{\beta,p}(\mathbb{R}^{n})$ and $H^{\beta}(\mathbb{R}^{n})$, as well as the definition of the fractional fractional Laplacian \cite{Nezza}.\\
\begin{defn}\label{d1.1}
 In the whole space, for $\beta\in (0,1)$, if $f\in  \mathscr{S}(\mathbb{R}^{n}) $, let
  $\Lambda^{\gamma}=(-\Delta)^{\beta}$ with $\gamma=2\beta$, and
\begin{equation*}
\widehat{\Lambda^{2\beta}f} (\xi)=\widehat{(-\Delta)^{\beta}f} (\xi)=|\xi|^{2\beta}\widehat{f}(\xi),\end{equation*}
the domain of definition of the fractional Laplacian, $\mathcal{D}\left(\Lambda^{\beta}\right) $ is endowed with a natural norm $\|\cdot\|_{\mathcal{D}\left(\Lambda^{\beta}\right) }$ and is a Hilbert space. The norm of $u$ in $\mathcal{D}\left(\Lambda^{\beta}\right) $ is defined by
\begin{equation}\label{1.8}
 \|u\|_{\mathcal{D}\left(\Lambda^{\beta}\right)}:=
\left\|\Lambda^{\beta}u\right\|_{L^2( \mathbb{R}^n) }. \end{equation}
It should be pointed out that in the whole space, if any function $\psi\in  \mathscr{S}(\mathbb{R}^{n}) $, $\mathcal{D}\left(\Lambda^{\beta}\right) $ is equivalent to the fractional Sobolev space $\dot{H}^{\beta}(\mathbb{R}^n) $, defined as the completion of $C^{\infty}_{0}(\mathbb{R}^n)$ with the norm
\begin{equation}\label{1.9}
\|\psi\|_{\dot{H}^{\beta}(\mathbb{R}^n)}=\left(\int_{\mathbb{R}^n}|\xi|^{2\beta}\left|\widehat{\psi}
\right|^{2}d\xi\right)^{\frac{1}{2}}=\left\|(-\Delta)^{\frac{\beta}{2}}\psi\right\|_{L^{2}(\mathbb{R}^n)}.
\end{equation}
On the other hand, the norm $\|u\|_{H^{\beta}(\mathbb{R}^n)}$ in the fractioal Laplacian Sobolev space $H^{\beta}(\mathbb{R}^n)$ is represented as
\begin{equation}\label{1.10}
\|\mathbf{u}\|^{2}_{H^{\beta}(\mathbb{R}^n)}:=2C(n,\beta)^{-1}\left\|\Lambda^{\beta} u\right\|_{L^2( \mathbb{R}^n) }^2+\|u\|_{L^2(\mathbb{R}^n)}^2.
\end{equation}
In particular, the norm of $\mathcal{D}\left(\Lambda^{2}\right)=\mathcal{D}(-\Delta)$ is equivalent to the $H^{2}(\mathbb{R}^n)$ norm.\hfill$\Box$
\end{defn}

\begin{defn}\label{d1.2}
Let $\beta\in (0,1)$. For any $p\in [1,\infty)$, we define $W^{\beta,p}(\mathbb{R}^{n})$ as follows
\begin{equation}\label{1.11}
W^{\beta,p}(\mathbb{R}^{n}):=\left\{u\in L^{p}(\mathbb{R}^{n}):~~\frac{|u(x)-u(y)|}{|x-y|^{\frac{n}{p}+\beta}}\in L^{p}(\mathbb{R}^{n}\times\mathbb{R}^{n})\right\},
\end{equation}
i.e., an intermediary Banach space between $L^{p}(\mathbb{R}^{n})$ and $W^{1,p}(\mathbb{R}^{n})$, endowed with the natural norm
\begin{equation}\label{1.12}
\|u\|_{W^{\beta,p}(\mathbb{R}^{n})}:=\left(\int_{\mathbb{R}^{n}}|u|^{p}dx
+\int_{\mathbb{R}^{n}}\int_{\mathbb{R}^{n}}\frac{|u(x)-u(y)|^{p}}{|x-y|^{n+\beta p}}
dxdy\right)^{\frac{1}{p}},
\end{equation}
where the term
\begin{equation}\label{1.13}
[u]_{W^{\beta,p}(\mathbb{R}^{n})}:=\left(\int_{\mathbb{R}^{n}}\int_{\mathbb{R}^{n}}
\frac{|u(x)-u(y)|^{p}}{|x-y|^{n+\beta p}}
dxdy\right)^{\frac{1}{p}}
\end{equation}
is the so-called Gagliardo (semi) norm of $u$. \\
\indent However, there is another case for {\bf $\beta\in (1,\infty)$} and $\beta$ is not an integer. In this case, we write $\beta=m+m'$, where $m$ is an integer and $m'\in (0,1)$. The space $W^{\beta,p}(\mathbb{R}^{n})$ consists of those equivalence classes of functions $u\in W^{m,p}(\mathbb{R}^{n})$ whose distributional derivatives $D^{\alpha}u$, with $|\alpha|=m$, belong to $W^{m',p}(\mathbb{R}^{n})$, namely
\begin{equation*}
W^{\beta,p}(\mathbb{R}^{n}):=\left\{u\in W^{m,p}(\mathbb{R}^{n}):~~D^{\alpha}u\in W^{m',p}(\mathbb{R}^{n})~~\mbox{for~~any}~~\alpha~s.t.~~|\alpha|=m\right\},
\end{equation*}
 and this is a Banach space with respect to the norm
 \begin{equation}\label{1.14}
\|u\|_{W^{\beta,p}(\mathbb{R}^{n})}:=\left( \left\|u\right\|^{p}_{W^{m,p}(\mathbb{R}^{n})}+\sum\limits_{|\alpha|=m}
\left\|D^{\alpha}u\right\|^{p}_{W^{m',p}(\mathbb{R}^{n})}\right)^{\frac{1}{p}}.
\end{equation}
 Clearly, if $\beta=m$ is an integer, the space $W^{\beta,p}(\mathbb{R}^{n})$ coincides with the Sobolev space $W^{m,p}(\mathbb{R}^{n})$.\\
 \indent Note that for any $\beta>0$, the space $C^{\infty}_{0}(\mathbb{R}^{n})$ of smooth functions with compact support is dense in $W^{\beta,p}(\mathbb{R}^{n})$, and $W^{\beta,p}_{0}(\mathbb{R}^{n})=W^{\beta,p}(\mathbb{R}^{n})$, where $W^{\beta,p}_{0}(\mathbb{R}^{n})$ denotes the closure of $C^{\infty}_{0}(\mathbb{R}^{n})$ in the space $W^{\beta,p}(\mathbb{R}^{n})$.\\
  \indent In particular, for $\beta\in (0,1)$ and $p=2$, the fractional Sobolev spaces $W^{\beta,2}(\mathbb{R}^{n})$ and $W^{\beta,2}_{0}(\mathbb{R}^{n})$ turn out to be Hilbert spaces, which are usually labeled by $W^{\beta,2}(\mathbb{R}^{n})=H^{\beta}(\mathbb{R}^{n})$ and $W^{\beta,2}_{0}(\mathbb{R}^{n})=H^{\beta}_{0}(\mathbb{R}^{n})$. That is,
\begin{equation}\label{1.15}
H^{\beta}(\mathbb{R}^{n}):=\left\{u\in L^{2}(\mathbb{R}^{n}):~~\frac{|u(x)-u(y)|}{|x-y|^{\frac{n}{2}+\beta}}\in L^{2}(\mathbb{R}^{n}\times\mathbb{R}^{n})\right\},
\end{equation}
i.e., an intermediary Hilbert space between $L^{2}(\mathbb{R}^{n})$ and $H^{1}(\mathbb{R}^{n})$, endowed with the natural norm
\begin{equation}\label{1.16}
\|u\|_{H^{\beta}(\mathbb{R}^{n})}:=\left(\int_{\mathbb{R}^{n}}|u|^{2}dx
+\int_{\mathbb{R}^{n}}\int_{\mathbb{R}^{n}}\frac{|u(x)-u(y)|^{2}}{|x-y|^{n+2\beta}}
dxdy\right)^{\frac{1}{2}},
\end{equation}
where the term
\begin{equation}\label{1.17}
[u]_{H^{\beta}(\mathbb{R}^{n})}:=\left(\int_{\mathbb{R}^{n}}\int_{\mathbb{R}^{n}}
\frac{|u(x)-u(y)|^{2}}{|x-y|^{n+2\beta}}
dxdy\right)^{\frac{1}{2}}
\end{equation}
is the so-called seminorm of $u$. \\
\indent There is an alternative definition of the space $H^{\beta}(\mathbb{R}^{n})$ via the Fourier transform. For any real $\beta\geq 0$, we may define
\begin{equation}\label{1.18}
\widehat{H}^{\beta}(\mathbb{R}^{n}):=\left\{u\in L^{2}(\mathbb{R}^{n}):~~
\int_{\mathbb{R}^{n}}\left(1+|\xi|^{2\beta}\right)|\mathcal{F}u(\xi)|^{2}d\xi<\infty
\right\}.
\end{equation}
In the same manner, for $\beta<0$ there is an analogous definition for $H^{\beta}(\mathbb{R}^{n})$:
 \begin{equation}\label{1.19}
H^{\beta}(\mathbb{R}^{n}):=\left\{u\in \mathscr{S}'(\mathbb{R}^{n}):~~
\int_{\mathbb{R}^{n}}\left(1+|\xi|^{2 }\right)^{\beta}|\mathcal{F}u(\xi)|^{2}d\xi<\infty
\right\}.
\end{equation}
 On the other hand, let $\beta\in (0,1)$ and let $(-\Delta)^{\beta}:~~\mathscr{S}\rightarrow L^{2}(\mathbb{R}^{n})$ be the fractional Laplacian operator defined by \eqref{1.6}. Then
 \\[0.2cm]
 \indent (1)  For any $u\in \mathscr{S}$,
  \begin{equation}\label{1.20}
  (-\Delta)^{\beta}u=\mathcal{F}^{-1}\left[|\xi|^{2\beta}\left(\mathcal{F}u\right)\right],~~\forall \xi \in \mathbb{R}^{n}.
\end{equation}
 \indent (2)  The fractional Sobolev space $H^{\beta}(\mathbb{R}^{n})$ defined in \eqref{1.15} coincides with $\widehat{H}^{\beta}(\mathbb{R}^{n})$ defined in \eqref{1.18}. In particular, for any $u\in H^{\beta}(\mathbb{R}^{n})$
 \begin{equation}\label{1.21}
 [u]^{2}_{H^{\beta}(\mathbb{R}^{n})}=2C(n,\beta)^{-1}
 \int_{\mathbb{R}^{n}}|\xi|^{2\beta}\left|\mathcal{F}u(\xi)\right|^{2}d\xi,
 \end{equation}
 where $C(n,\beta)$ is defined by \eqref{1.7}.
 \\[0.2cm]
 \indent (3)  For $u\in H^{\beta}(\mathbb{R}^{n})$,
 \begin{equation}\label{1.22}
 [u]^{2}_{H^{\beta}(\mathbb{R}^{n})}=2C(n,\beta)^{-1}
 \left\|(-\Delta)^{\frac{\beta}{2}}u\right\|_{L^{2}(\mathbb{R}^{n})}^{2},
 \end{equation}
 where  $C(n,\beta)$ is defined by \eqref{1.7}.\hfill$\Box$
\end{defn}

\indent Recently, a great attention has been devoted to the study of nonlocal problems driven by fractional Laplacian type operators in the literature, not only for a pure academic interest, but also for the various applications in different fields. It is well-known that fractional Laplacian $(-\Delta)^{\beta}$ is a spatial integro-differential operator, and that it can be used to describe the spatial nonlocality and power law behaviors in various science and engineering problems. In the recent two decades, fractional Laplacian has been utilized to model energy dissipation of acoustic propagation in human tissue \cite{Chen-Holm}, turbulence diffusion \cite{Chen-W}, contaminant transport in ground water \cite{Pang-Chen}, non-local heat conduction \cite{Bobaru,Elia,Mongiovi}, and electromagnetic fields on fractals \cite{Tarasov}.\\
\indent Before going further, we first speak of some results for the system \eqref{1.3}. The non-uniform decay and algebraic decay were considered in \cite{Clayton-08}. Concerning the convergence from  \eqref{1.3} to the imcompressible Navier-Stokes equations \eqref{1.4}, the authors in
\cite{Gallay,Giga} proved how the solution of \eqref{1.3} approaches a solution to  \eqref{1.4} weakly when the filter parameter $\alpha$ tends to zero. Bjorland and Schonbek in \cite{Clayton-08} showed how solutions to \eqref{1.3} approach solutions to the Navier-Stokes equation \eqref{1.4} strongly as $\alpha$ tends to zero when the solutions to \eqref{1.4} is sufficiently regular. In \cite{Bjorland}, Bjorland investigated the relationship between solutions of the Navier-Stokes equations \eqref{1.4} and the Camassa-Holm equations \eqref{1.3} by describing the way a solution of \eqref{1.3} approaches the fixed point zero, i.e., computing the first and second order decay asymptotics for solutions with small initial data. However, there are some similar results for the Navier-Stokes equations \eqref{1.4}. Decay results for \eqref{1.4} has been studied in the literature \cite{Carpio,Gallay,Giga,Schonbek-1,Schonbek-2}. The asymptotic behavior of the 2-D vorticity equation for \eqref{1.4} has been investigated in \cite{Carpio,Gallay,Giga}. In \cite{Carpio}, Carpio studied the asymptotic behavior for the vorticity equation for \eqref{1.4} in two and three space dimensions.
 Gallay and Wayne in \cite{Gallay} calculated the asymptotics by applying invariant manifold technique to the semiflow governing the vorticity problem for \eqref{1.4}. The large time behavior of the vorticity of two-dimensional viscous flow for \eqref{1.4} was established by Giga and Kambe in \cite{Giga}. \\
\indent In stark contrast to those works on the study of the Camassa-Holm equations \eqref{1.3} in recent decades, little has been known on the space-fractional derivative viscosity in the literature despite that non-standard diffusions are very natural also for these problems. In particular, the study of the Camassa-Holm equations with fractional Laplacian viscosity \eqref{1.1} is more challenging due to the vector integral expression and nonlocal property.\\
 \indent The aim of this paper is twofold. We first intend to establish the large time behavior of solutions to the nonlocal Camassa-Holm equations \eqref{1.1} , which  concerns the non-uniform decay and algebraic decay. Our second goal is to discuss the relation betweem the equations under study and the incompressible Navier-Stokes equations with fractional Laplacian viscosity:
\begin{equation}\label{1.23}
\left\{
\begin{array}{ll}
\mathbf{v}_{t}+\mathbf{v}\cdot \nabla \mathbf{v} +\nabla p=-\nu (-\Delta)^{\beta} \mathbf{v},
\\[0.2cm]
\mbox{div}~ \mathbf{v}=0.
\end{array}
\right.
 \end{equation}
 To achieve these results, the main difficulty lies in the fractional Laplacian viscosity. Fortunately, with the help of some properties of fractional Laplacian introduced in \cite{Caffarelli-1,P-I}, in particular, the fractional Leibniz chain rule and the fractional Gagliardo-Nirenberg-Sobolev type estimates, we first establish the large time behavior concerning non-uniform decay and algebraic decay of solutions to the nonlocal equations under study by applying the high and low frequency splitting method first used in \cite{Masuda} and the Fourier splitting method introduced in \cite{Kozono-1,Kozono-2}. In particular, under the critical case $s=\dfrac{n}{4}$, the nonlocal version of Ladyzhenskaya's inequality is skillfully used, and the smallness of initial data in several Sobolev spaces is required to gain the non-uniform decay and algebraic decay. On the other hand, by means of the fractional heat kernel estimates \cite{Dong-Li-1} and Leray projector, we figure out the relation between the nonlocal equations \eqref{1.1} and the imcompressible viscous nonlocal Navier-Stokes equations \eqref{1.23}. Specifically, we prove that the solution to \eqref{1.1} converges strongly as the filter parameter $\alpha\rightarrow~0$ to a solution of the imcompressible nonlocal equations \eqref{1.23}.
 \\[0.3cm]
 \indent We now give some remarks on the nonlocal Camassa-Holm equations \eqref{1.1}.
 \begin{rem}\label{r1.3} \rm
  Fractional power of the Laplacian arises in a numerous variety of equations in mathematical physics and related fields \cite{Abdelouhab-Bona,Caffarelli-1,Elgart,Frohlich,Kenig-Martel,Lieb-Yau,Majda-Mc,Weinstein}. In stark contrast to the problem on the large time behavior for the Camassa-Holm equations without any nonlocal term \eqref{1.3}, it seems fair to say that extremely little is known about the large time behavior for the solutions to the nonlocal equations \eqref{1.1}-\eqref{1.2}in two and three space dimensions. Indeed, to our best knowledge, the only example in \cite{Gui-Liu} for which some results of the Camassa-Holm equations with fractional dissipation in one space dimension have been shown is the following:
 \\[0.30cm]
\indent $\bullet$ Global well-posedness and blow-up of solutions to the Camassa-Holm equations with fractional dissipation under the supercritical case: $\gamma\in \left[\frac{1}{2},1\right)$.
\\[0.3cm]
\indent $\bullet$ The zero filter limit of the Camassa-Holm equation with fractional dissipation, as well as the possible blow-up of solutions under the subcritical case: $0\leq \gamma <\frac{1}{2}$.\hfill $\Box$\\
 \end{rem}
 \begin{rem}\label{r1.4} \rm
 It should be pointed out that the initimate relation between the Navier-Stokes equations and the viscous Camassa-Holm equations gives hope that a similar program for Navier-Stokes equations may be realized for the viscous Camassa-Holm equations. As a matter of fact, such attempts are met with resistance from the filter in the viscous Camassa-Holm equations. On one hand, in a functional setting the filter eases problems by smoothing the solution. On the other hand, in a dynamical setting the filter adds complication to the problems. In particular, the filters does not scale well with the other parts of the equations, and the resulting nonlinear term has dependence on the scaled time variable which has not been presented in the Navier-Stokes equations.\hfill$\Box$
 \end{rem}
 \indent We end this section by describing the notation we shall use in this paper.
 \\[0.3cm]
 {\bf Notations}
 \\[0.3cm]
$\mathscr{S}(\mathbb{R}^{n})$ denotes the Schwartz calss. The $i^{th}$ component of $\mathbf{v}\cdot \nabla \mathbf{u}^{T}$ is denoted by $\left(\mathbf{v}\cdot \nabla \mathbf{u}^{T}\right)_{i}=\sum\limits_{j=1}^{n}v_{j}\partial_{i}u_{j}$. Let
     $\displaystyle\left \langle \mathbf{u},\mathbf{v}\right\rangle=\int_{ \mathbb{R}^{n} }\mathbf{u}\cdot \mathbf{v}dx$ and $\displaystyle\Sigma =\left\{\phi \in C^{\infty}_{0}( \mathbb{R}^{n} )|\nabla\cdot \phi=0\right\}.$ $L^{p}_{0}(\mathbb{R}^{n})$ denotes the closure of  $C^{\infty}_{0}(\mathbb{R}^{n})$ in the space  $L^{p}(\mathbb{R}^{n})$ and $H^{m}_{0}(\mathbb{R}^{n})$  the completion of $C^{\infty}_{0}( \mathbb{R}^{n} )$ in the norm $\|\cdot\|_{H^{m}(\mathbb{R}^{n})}$. We denote by $L^{p}( \mathbb{R}^{n} )$ the standard Lebesgue space, and
    $L^{p}_{\sigma}( \mathbb{R}^{n} )$ the completion of $\Sigma$ in the norm $\|\cdot\|_{L^{p}( \mathbb{R}^{n} )}$. The completion of $\Sigma$ under the $\mathcal{D}\left(\Lambda^{\beta}\right)( \mathbb{R}^{n} )$-norm is denoted by $\mathcal{D}_{\sigma}\left(\Lambda^{\beta}\right)( \mathbb{R}^{n} )$ and $\left(\mathcal{D}_{\sigma}\left(\Lambda^{\beta}\right)\right)'( \mathbb{R}^{n} )$ is its dual space. The completion of $\Sigma$ under the $H^{m}( \mathbb{R}^{n} )$-norm will be denoted by  $H^{m}_{\sigma}( \mathbb{R}^{n} )$ and
      $(H^{m}_{\sigma}( \mathbb{R}^{n} ))'$  be the corresponding dual space. $\mathcal{F}(\phi)$ or $\hat{\phi}$ denotes the Fourier transform of a function $\phi$, with $\mathcal{F}^{-1}(\phi)$ or $ \breve{\phi} $  the inverse Fourier transform. For $a\lesssim b$, we mean that there is a uniform constant $C$, which may be different on different lines, such that $a\leq C b$. For $\beta=1$, $\mathcal{D} (-\Delta)=H^{2}(\mathbb{R}^{n})\cap H^{1}_{0}(\mathbb{R}^{n}) $ and $\mathcal{D} (\Lambda) (\mathbb{R}^{n})=H^{1}_{0}(\mathbb{R}^{n}) $. Generally, the letter $C$ will denote a generic constant. \hfill$\Box$
      \\[0.3cm]
           \indent The rest of the paper is organized as follows: in Section 2 we collect some preliminaries. In Section 3 the non-uniform decay is established. Subsequently, in Section 4 we show the algebraic decay. In the last section (Section 5), we prove the convergence from the solution of \eqref{1.1}-\eqref{1.2} to the imcompressible Navier-Stokes equations with nonlocal viscosity \eqref{1.23}.
\section{Preliminaries}
\allowdisplaybreaks
In this section, we collect several preliminary results.
\begin{lem} \label{l2.1}
 Let $\mathbf{u}$ and $\mathbf{v}$ be smooth divergence free functions with compact support. Then one has
$$
\left.
\begin{array}{ll}
&\displaystyle\mathbf{u}\cdot\nabla \mathbf{v}+ \sum \limits^{n}_{j=1} v_{j}\nabla u_{j}=-\mathbf{u}\times(\nabla\times \mathbf{v})+\nabla(\mathbf{v}\cdot \mathbf{u}),
\\[0.3cm]
&\displaystyle\left\langle \mathbf{u}\cdot\nabla \mathbf{v},\mathbf{u}\right\rangle+\left\langle \mathbf{v}\cdot\nabla \mathbf{u}^{T},\mathbf{u}\right\rangle=0,
\\[0.3cm]
&\displaystyle\left\langle \mathbf{u}\times(\nabla\times \mathbf{v}),\mathbf{u}\right\rangle=0.
\end{array}
\right.
$$
\end{lem}
{\bf Proof.} ~By direct calculation, it is easy to achieve these expected identities.\hfill$\Box$
\begin{lem} \label{l2.2}
For $n=2,3$ and $0<\beta < 1$,  let $\mathbf{u}$ and $\mathbf{v}$ be smooth divergence free functions with compact supports. Then if $(\mathbf{v},\mathbf{u})$ solves \eqref{1.1}-\eqref{1.2} , there holds
\begin{equation}\label{2.1}
\frac{1}{2}\frac{d}{dt}\left(\langle \mathbf{u},\mathbf{u}\rangle+\alpha^{2}\langle \nabla \mathbf{u},\nabla \mathbf{u}\rangle\right) +\nu\left(\left\langle\Lambda^{\beta}  \mathbf{u},\Lambda^{\beta} \mathbf{u}\right\rangle+\alpha^{2}\left\langle\nabla\Lambda^{\beta}   \mathbf{u},\nabla\Lambda^{\beta} \mathbf{u}\right\rangle\right)=0,
 \end{equation}
and
\begin{equation}\label{2.2}
\left.
\begin{array}{ll}
 &\displaystyle\|\mathbf{u}(\cdot,t)\|_{L^2(\mathbb{R}^{n})}^{2}+\alpha^{2}\|\nabla \mathbf{u}(\cdot,t)\|_{L^2(\mathbb{R}^{n})}^{2}
 \\[0.30cm]
&\displaystyle\qquad\qquad +2\nu\int_0^t \left\|\Lambda^{\beta} \mathbf{u}(\cdot,t)\right\|_{L^2(\mathbb{R}^{n})}^{2}dt
+2\nu\alpha^{2}\int_0^t\left\|\nabla\Lambda^{\beta} \mathbf{u}(\cdot,t)\right\|_{L^2(\mathbb{R}^{n})}^{2}dt
\\[0.40cm]
 &\displaystyle\qquad  \quad\leq\|\mathbf{u}_{0}\|_{L^2(\mathbb{R}^{n})}^{2}+\alpha^{2}\|\nabla \mathbf{u}_{0}\|_{L^2(\mathbb{R}^{n})}^{2}.
 \end{array}
\right.
\end{equation}
\end{lem}
 {\bf Proof.}\quad Thanks to $\mbox{div}~ \mathbf{v}=\mbox{div}~ \mathbf{u}=0, $ making inner product with $\mathbf{u}$ on the both sides in the first equation in \eqref{1.1} gives
$$\left\langle \mathbf{v}_{t},\mathbf{u}\right\rangle+\left\langle \mathbf{u}\cdot \nabla \mathbf{v}+\mathbf{v}\cdot \nabla \mathbf{u}^{T},\mathbf{u}\right\rangle+ \left\langle\nabla p,\mathbf{u}\right\rangle+\nu\left\langle(-\Delta)^{\beta} \mathbf{v},\mathbf{u}\right\rangle=0.$$
Note that Lemma \ref{l2.1}, one deduces by integrating by parts
$$\langle \mathbf{v}_{t},\mathbf{u}\rangle+\nu\left\langle(-\Delta)^{\beta}\mathbf{v},  \mathbf{u}\right\rangle=\langle \mathbf{v}_{t},\mathbf{u}\rangle+\nu\left\langle \Lambda^{\beta} \mathbf{v},  \Lambda^{\beta} \mathbf{u}\right\rangle=0.$$
This together with the second equation in \eqref{1.1} concludes that
\begin{equation*}
\left.
\begin{array}{ll}
&\displaystyle\left\langle \mathbf{u}_{t}-\alpha^{2} \Delta \mathbf{u}_{t},\mathbf{u}\right\rangle+\nu\left\langle(-\Delta)^{\beta}(\mathbf{u}-\alpha^{2}\Delta \mathbf{u}), \mathbf{u}\right\rangle
\\[0.3cm]
&\displaystyle\qquad=\frac{1}{2}\frac{d}{dt}\left(\langle \mathbf{u},\mathbf{u}\rangle+\alpha^{2}\langle\nabla \mathbf{u},\nabla \mathbf{u}\rangle\right) +\nu\left(\left\langle\Lambda^{\beta} \mathbf{u},\Lambda^{\beta}  \mathbf{u}\right\rangle+\alpha^{2}\left\langle\nabla\Lambda^{\beta}  \mathbf{u},\nabla\Lambda^{\beta}  \mathbf{u}\right\rangle\right)
\\[0.3cm]
&\displaystyle\qquad=0.
\end{array}
\right.
\end{equation*}
This is the equality \eqref{2.1}. \eqref{2.2} follows by integrating both sides of \eqref{2.1} with respect to $t$.\hfill$\Box$
\\[0.3cm]
\indent Before going further, we introduce the following notion of weak solutions to the Camassa-Holm equations with fractional Laplacian viscosity \eqref{1.1}-\eqref{1.2} in $\mathbb{R}^{n}~(n=2,3)$.

\begin{defn}\label{d2.3}
Let $\frac{n}{4}\leq \beta< 1$ with $n=2,3$. A weak solution to \eqref{1.1}-\eqref{1.2} is a pair of functions $(\mathbf{v},\mathbf{u})$ such that
\begin{eqnarray*}
&&
\mathbf{v}\in \emph{L}^{\infty}\left([0,T];\emph{L}_{\sigma}^{2}(\mathbb{R}^{n})\right)\cap \emph{L}^{2}\left([0,T];\mathcal{D}_{\sigma}\left(\Lambda^{\beta}\right)(\mathbb{R}^{n})\right),
\\[0.2cm]
&&\partial_{t}\mathbf{v}\in\emph{L}^{2}\left([0,T];\mathcal{B}' \right),
\\[0.2cm]
&&\mathbf{u}\in\emph{L}^{\infty}\left([0,T]; \emph{H}_\sigma^{2} (\mathbb{R}^{n})\right)\cap \emph{L}^{2}\left([0,T];\mathcal{D}_\sigma\left(\Lambda^{2+\beta}\right)(\mathbb{R}^{n})\right),
\\[0.2cm]
 &&\mathbf{v}(0,x)=\mathbf{v}_{0}(x).
\end{eqnarray*}
Here, $\displaystyle \mathcal{B}=\left\{\begin{array}{ll}
\mathcal{D}_{\sigma}\left(\Lambda^{\beta}\right)(\mathbb{R}^{n})~~&\hbox{for}~~\frac{n}{4}< \beta < 1,
\\[0.3cm]
H_{\sigma}^{\frac{n}{4}}(\mathbb{R}^{n})~~&\hbox{for}~~\beta=\frac{n}{4}.\end{array}\right.$
\\[0.20cm]
In addition, for every ${\phi}\in\emph{L}^{2}\left([0,T];\mathcal{E}\right)$ with ${\phi(T)}=0$, there holds
\begin{eqnarray*}
\vspace{5pt}
&&-\int_{0}^{T}\left\langle\mathbf{ v},\partial_{t}\phi\right\rangle ds+\int_{0}^{T}\left\langle \mathbf{u}\cdot\nabla \mathbf{v},\phi \right\rangle ds +\int_{0}^{T}\left\langle\phi\cdot\nabla \mathbf{u},\mathbf{v}\right\rangle ds+\nu\int_{0}^{T}\left\langle \Lambda ^{\beta}\mathbf{v},\Lambda ^{\beta} \phi\right\rangle ds\\
\\
&&\qquad\qquad\qquad\quad={\left\langle \mathbf{v}_{0},\phi(0)\right\rangle}.
\end{eqnarray*}
In particular, for $t\in [0,T]$ there holds
$$\left\langle\mathbf{u},\phi\right\rangle+\alpha^{2}\left\langle\nabla \mathbf{u}, \nabla\phi\right\rangle=\langle \mathbf{v},\phi\rangle. $$
Here, $\displaystyle \mathcal{E}=\left\{\begin{array}{ll}
\mathcal{D}_{\sigma}\left(\Lambda\right)(\mathbb{R}^{n})~~&\hbox{for}~~\frac{n}{4}< \beta < 1,
\\[0.3cm]
H_{\sigma}^{1}(\mathbb{R}^{n})~~&\hbox{for}~~\beta=\frac{n}{4}.\end{array}\right.$
$$\eqno \Box$$
\end{defn}

\indent  Let $\mathscr{S}(\mathbb{R}^{n}) $ be the Schwartz class. The nonlocal operator $(-\Delta)^{\beta}$ is defined for any $g\in \mathscr{S}(\mathbb{R}^{n}) $ through the Fourier transform: if $(-\Delta)^{\beta}g=h$, then
$$
 \widehat{h}(\xi)=|\xi|^{2\beta}\widehat{g}(\xi).\eqno(D-1)$$
 It should be pointed out that if $\psi$ and $\phi$ belong to the Schwartz class $\mathscr{S}(\mathbb{R}^{n})$, (D-1) together with Plancherel's theorem yields
$$\int_{\mathbb{R}^{n} }(-\Delta)^{\beta}\psi\phi dx=\int_{\mathbb{R}^{n} }|\xi|^{2\beta}\widehat{\psi}\widehat{\phi} d\xi=\int_{\mathbb{R}^{n} }(-\Delta)^{\frac{\beta}{2}}\psi(-\Delta)^{\frac{\beta}{2}}\phi dx.$$

Thanks to Theorem 3.1 and Theorem 4.1 in \cite{Gan-Lin-Tong} (see also \cite{Gan-HMW}), using the energy method and a bootstrap argument, we obtain the following proposition concerning the existence, uniqueness and regularity of a weak solution to \eqref{1.1}-\eqref{1.2}:
 \begin{prop}\label{p2.4}
 Let $\frac{n}{4}\leq \beta < 1$ with $n=2,3$. Assume that
 \\[0.2cm]
 \indent (1) for $\displaystyle\frac{n}{4}<\beta < 1$,  $\mathbf{v}_{0}\in \mathcal{D}_{\sigma}\left(\Lambda^{M} \right)(\mathbb{R}^{n} )$, $M\geq 0$,\\
  and\\
   \indent (2) for $\displaystyle \beta=\frac{n}{4}$, $\mathbf{v}_{0}\in H_{\sigma}^{M}(\mathbb{R}^{n})$, $M\geq 0$, and in addition, there exists an $\varepsilon^{*}=\varepsilon^{*}(\alpha,\nu,n)$ sufficiently small such that $\left\| \mathbf{v}_{0}\right\|_{H_{0}^{M}(\mathbb{R}^{n}) }\leq \varepsilon^{*}$.
   \\[0.2cm]
     Then there exists a unique weak solution to \eqref{1.1}-\eqref{1.2} in the sense of Definition \ref{d2.3}. In addition, this solution satisfies the energy estimate \eqref{2.2}, and for all $m+2k\beta\leq M$ there holds
\begin{equation}\label{2.3}
\left\|\partial^{k}_{t}\nabla^{m}\mathbf{v}\right\|^{2}_{L^{2}( \mathbb{R}^{n} )}+\nu\int_{0}^{T}
\left\|\partial^{k}_{t}\nabla^{m}\Lambda^{\beta} \mathbf{v}\right\|^{2}_{L^{2}( \mathbb{R}^{n} )}dt\leq C\left(n,\beta,\alpha,\nu,\| \mathbf{v}_{0}\|_{\mathcal{A}}\right),\end{equation}
 where  $\displaystyle \mathcal{A}=\left\{\begin{array}{ll}
\mathcal{D}\left(\Lambda^{M}\right)(\mathbb{R}^{n})~~&\hbox{for}~~\frac{n}{4}< \beta < 1,
\\[0.3cm]
H_{0}^{M}(\mathbb{R}^{n})~~&\hbox{for}~~\beta=\frac{n}{4},\end{array}\right.$  $m$ and $k$ are both non-negative integers. \hfill$\Box$
$$\qquad \eqno \Box$$
 \end{prop}
 By applying the Gagliardon-Nirenberg-Sobolev inequality to the bound \eqref{2.3} in Proposition \ref{p2.4}, we achieve a corollary which describes the action of the filter.
  \begin{cor}\label{c2.5} \rm
   For $\frac{n}{4}\leq \beta< 1$ with $n=2,3$, let $(\mathbf{v},\mathbf{u})$ be the solution to the Cauchy problem \eqref{1.1}-\eqref{1.2} constructed in Proposition \ref{p2.4}. Then the following estimates hold for all $m+2k\beta\leq M$:
\begin{equation}\label{2.4}
\left.
\begin{array}{ll}
 \displaystyle\left\|\partial^{k}_{t}\nabla^{m} \mathbf{u}\right\|^{2}_{L^{2}(\mathbb{R}^{n})}
+2\alpha^2\left\|\partial^{k}_{t}\nabla^{m+1}  \mathbf{u}\right\|^{2}_{L^{2}(\mathbb{R}^{n})}
+\alpha^4\left\|\partial^{k}_{t}\nabla^{m+2}  \mathbf{u}\right\|^{2}_{L^{2}(\mathbb{R}^{n})}
=\left\|\partial^{k}_{t}\nabla^{m} \mathbf{v}\right\|^{2}_{L^{2}(\mathbb{R}^{n})},
\end{array}
\right.
\end{equation}
\begin{equation}\label{2.5}
\left.
\begin{array}{ll}
 \displaystyle\left\|\partial^{k}_{t}\nabla^{m}\mathbf{u}\right\|^{2}_{L^{n}(\mathbb{R}^{n})}
+\left\|\partial^{k}_{t}\nabla^{m}\Lambda^{\beta} \mathbf{u}\right\|^{2}_{L^{n}(\mathbb{R}^{n})}
+\left\|\partial^{k}_{t}\nabla^{m+1}\mathbf{u}\right\|^{2}_{L^{n}(\mathbb{R}^{n})}
  \lesssim\left\|\partial^{k}_{t}\nabla^{m}  \mathbf{v}\right\|^{2}_{L^{2}(\mathbb{R}^{n})},
\end{array}
\right.
\end{equation}
\begin{equation}\label{2.6}
\left.
\begin{array}{ll}
 \displaystyle\left\|\partial^{k}_{t}\nabla^{m}\mathbf{u}\right\|^{2}_{L^{n}(\mathbb{R}^{n})}
+\nu\int^{t}_{0}\left\|\partial^{k}_{t}\nabla^{m}\Lambda^{\beta} \mathbf{u}\right\|^{2}_{L^{n}(\mathbb{R}^{n})}ds
\leq C\left(n,\beta,\alpha,\nu,m,k,\|\mathbf{v}_{0}\|_{\mathcal{D}\left(\Lambda^{M} \right)(\mathbb{R}^{n})}\right).
\end{array}
\right.
\end{equation}
Here, $m$ and $k$ are both non-negative integers.
$$\qquad \eqno \Box$$
  \end{cor}
 We then claim a lemma concerning the Helmholtz equation $\mathbf{u}-\alpha^{2}\Delta \mathbf{u}=\mathbf{v}$.
\begin{lem}\label{l2.6}\rm
 Let $n=2,3$ and $\frac{n}{4}\leq \beta < 1$. Given $\mathbf{v}\in w^{\beta,p}(\mathbb{R}^{n})$ with $1<p<\infty$, there exists a weak solution  $\mathbf{u}\in W^{2,p}(\mathbb{R}^{n})$ to the Helmholtz equation $\mathbf{u}-\alpha^{2}\Delta \mathbf{u}=\mathbf{v}$ such that the following estimates hold:
 \begin{equation*}
\left.
\begin{array}{ll}
\|\mathbf{u}\|_{L^p(\mathbb{R}^{n})}\leq \|\mathbf{v}\|_{L^p(\mathbb{R}^{n})},
\\[0.3cm]
\|\mathbf{u}\|_{L^q(\mathbb{R}^{n})}\leq \dfrac{C(n,p,q)}{\alpha^{ 2\gamma_{1}}}\|\mathbf{v}\|_{L^p( \mathbb{R}^{n})},
\\[0.3cm]
\|\nabla \mathbf{u}\|_{L^q( \mathbb{R}^{n})}\leq \dfrac{C(n,p,q)}{\alpha^{1+2\gamma_{2}}}\|\mathbf{v}\|_{L^p( \mathbb{R}^{n})},
\\[0.3cm]
\|\Delta \mathbf{u}\|_{L^q( \mathbb{R}^{n})}\leq \dfrac{C(n,p,q)}{\alpha^{2+2\gamma_{3}-\beta}}\left\|\Lambda^{\beta}\mathbf{v} \right\|_{L^p( \mathbb{R}^{n})},
\end{array}
\right.
\end{equation*}
where $w^{\beta,p}(\mathbb{R}^{n})$ is defined by Definition \ref{d1.2},$\gamma_{1}=\frac{n}{2}\left(\frac{1}{p}-\frac{1}{q}\right)<1,~
\gamma_{2}=\frac{n}{2}\left(\frac{1}{p}-\frac{1}{q}\right)<\frac{1}{2},
~\gamma_{3}=\frac{n}{2}\left(\frac{1}{p}-\frac{1}{q}\right)<\frac{\beta}{2}$. In particular, there holds that for $0<\alpha<1$:
\begin{equation*}
\left.
\begin{array}{ll}
&\displaystyle\|\mathbf{u}\|_{L^q(\mathbb{R}^{n})}\leq \dfrac{C(n,p,q)}{\alpha^{1+\gamma_{1}}}\|\mathbf{v}\|_{L^p( \mathbb{R}^{n})},
\\[0.3cm]
&\displaystyle\|\nabla \mathbf{u}\|_{L^q( \mathbb{R}^{n})}\leq \dfrac{C(n,p,q)}{\alpha^{\frac{3}{2}+\gamma_{2}}}\|\mathbf{v}\|_{L^p( \mathbb{R}^{n})},
\\[0.3cm]
&\displaystyle\|\Delta \mathbf{u}\|_{L^q( \mathbb{R}^{n})}\leq \dfrac{C(n,p,q)}{\alpha^{2-\frac{\beta}{2}+\gamma_{3}}}\left\|\Lambda^{\beta}\mathbf{v}\right\|_{L^p( \mathbb{R}^{n})}.
\end{array}
\right.
\end{equation*}
In addition, if $n\left(\frac{2}{p}-1\right)<\beta$, then the solution is unique.
\end{lem}
 {\bf Proof.}   Note that $1-\alpha^{2}\Delta$ is a strictly positive, compact and self-adjoint operator, using standard elliptic theory and making suitable scaling on spatial variables, Sobolev embedding theorem and interpolation inequalities deduce the expected estimates.\hfill$\Box$\\
 \begin{lem} \label{l2.7}\rm
For $\frac{n}{4}\leq \beta < 1$, $n=2,3$, let
  $$E(t)\in C^{1}\left([0,\infty)\right) ,~   \psi\in C^{1}\left([0,\infty), C^{1}\cap L^{2}( \mathbb{R}^{n})\right) ,~\widetilde{\psi} \in C^{1}\left([0,\infty),  L^{\infty}( \mathbb{R}^{n})\right).$$
  Solutions of \eqref{1.1}-\eqref{1.2} constructed in Proposition \ref{p2.4} admit the following genearlized energy inequalities:
\begin{equation}\label{2.7}
\left.
\begin{array}{ll}
& \displaystyle E(t)\left\|\psi(t)\ast\mathbf{ v}(t)\right\|^{2}_{L^{2}(\mathbb{R}^{n})}
\\[0.3cm]
&\displaystyle \qquad=E(s)\left\|\psi(s)\ast \mathbf{v}(s)\right\|^{2}_{L^{2}(\mathbb{R}^{n})}+\displaystyle\int^{t}_{s}E'(\tau)\left\|\psi(\tau)\ast \mathbf{v}(\tau)\right\|^{2}_{L^{2}(\mathbb{R}^{n})}d\tau
\\[0.3cm]
&\displaystyle\qquad\quad +2\displaystyle\int^{t}_{s}E(\tau)\left\langle\psi'(\tau)\ast \mathbf{v}(\tau),\psi(\tau)\ast \mathbf{v}(\tau)\right\rangle
d\tau
\\[0.3cm]
&\displaystyle\qquad\quad-2\nu\displaystyle\int^{t}_{s}E(\tau)\left\|\Lambda^{\beta}\psi(\tau)\ast \mathbf{v}(\tau)\right\|^{2}_{L^{2}(\mathbb{R}^{n})}d\tau
\\[0.3cm]
&\displaystyle\qquad\quad-2 \displaystyle\int^{t}_{s}E(\tau)
\left\langle \mathbf{u}\cdot \nabla \mathbf{v}, \psi(\tau)\ast\psi(\tau)\ast \mathbf{v}(\tau)\right\rangle d\tau
\\[0.3cm]
 &\displaystyle\qquad\quad-2 \displaystyle\int^{t}_{s}E(\tau)
\left\langle \mathbf{v}\cdot \nabla \mathbf{u}^{T}, \psi(\tau)\ast\psi(\tau)\ast \mathbf{v}(\tau)\right\rangle d\tau,
\end{array}
\right.
\end{equation}
and
\begin{equation}\label{2.8}
\left.
\begin{array}{ll}
 & \displaystyle E(t)\left\|\widetilde{\psi}(t)  \hat{\mathbf{v}}(t)\right\|^{2}_{L^{2}(\mathbb{R}^{n})}
 \\[0.3cm]
 &\displaystyle\qquad=E(s)\left\|\widetilde{\psi}(s) \hat{\mathbf{v}}(s)\right\|^{2}_{L^{2}(\mathbb{R}^{n})}+\displaystyle\int^{t}_{s}
E'(\tau)\left\|\widetilde{\psi}(\tau)\hat{\mathbf{v}}(\tau)\right\|^{2}_{L^{2}(\mathbb{R}^{n})}d\tau
\\[0.3cm]
 &\displaystyle\qquad\quad +2\displaystyle\int^{t}_{s}E(\tau)\left\langle\widetilde{\psi}' (\tau) \hat{\mathbf{v}}(\tau),\widetilde{\psi} (\tau) \hat{\mathbf{v}}(\tau)  \right\rangle d\tau
 \\[0.3cm]
 &\displaystyle\qquad\quad -2\nu\displaystyle\int^{t}_{s}E(\tau)\left\|\xi^{\beta}\widetilde{\psi} (\tau) \hat{\mathbf{v}}(\tau) \right\|^{2}_{L^{2}(\mathbb{R}^{n})}d\tau
 \\[0.3cm]
&\displaystyle\qquad\quad -2\displaystyle\int^{t}_{s}E(\tau)\left\langle \mathcal{F}(\mathbf{u}\cdot \nabla \mathbf{v}) \widetilde{\psi}^{2}(\tau) \hat{\mathbf{v}}(\tau)\right \rangle
d\tau
\\[0.3cm]
&\displaystyle\qquad\quad -2\displaystyle\int^{t}_{s}E(\tau)\left\langle \mathcal{F}(\mathbf{v}\cdot \nabla \mathbf{u}^{T}) \widetilde{\psi}^{2}(\tau) \hat{\mathbf{v}}(\tau) \right\rangle
d\tau.\end{array}
\right.
\end{equation}
\end{lem}
{\bf Proof.}\quad Multiplying the first equation in \eqref{1.1} by $E(t)\psi\ast\psi\ast \mathbf{v}$ and  integrating in space variable $x$ yields, after some integration by parts,
\begin{equation}\label{2.9}
\left.
\begin{array}{ll}
&\displaystyle\int _{\mathbb{R}^{n}}\partial_{t}\mathbf{v}\cdot E(t)\psi\ast\psi\ast \mathbf{v}dx+\displaystyle\int _{\mathbb{R}^{n}}\mathbf{u}\cdot \nabla \mathbf{v} \cdot E(t)\psi\ast\psi\ast \mathbf{v}dx
\\[0.3cm]
&\displaystyle\qquad+\displaystyle\int _{\mathbb{R}^{n}}\mathbf{v}\cdot \nabla \mathbf{u}^{T} \cdot E(t)\psi\ast\psi\ast \mathbf{v}dx+\displaystyle\int _{\mathbb{R}^{n}}\nabla p \cdot E(t)\psi\ast\psi\ast \mathbf{v}dx
\\[0.3cm]
&\displaystyle\qquad+\nu\displaystyle\int _{\mathbb{R}^{n}}(-\Delta)^{\beta}\mathbf{v}\cdot E(t)\psi\ast\psi\ast \mathbf{v}dx=0.
\end{array}
\right.
\end{equation}
Rearranging \eqref{2.9} gives rise to
\begin{equation}\label{2.10}
\left.
\begin{array}{ll}
&\displaystyle\dfrac{d}{dt}\left(E(t)\left\|\psi(t)\ast \mathbf{v}(t) \right\|^{2}_{L^{2}(\mathbb{R}^{n})}\right)
\\[0.3cm]
&\displaystyle\qquad=E'(t)\left\|\psi(t)\ast \mathbf{v}(t) \right\|^{2}_{L^{2}(\mathbb{R}^{n})}+2E(t)\displaystyle\int _{\mathbb{R}^{n}}\left(\psi'(t)\ast \mathbf{v}(t))(\psi(t)\ast \mathbf{v}(t)\right)dx
\\[0.3cm]
&\displaystyle\qquad\quad-2\nu E(t)\displaystyle\int _{\mathbb{R}^{n}}\left|\Lambda^{\beta}\psi(t)\ast \mathbf{v}(t)\right|^{2}dx-2E(t)\displaystyle\int _{\mathbb{R}^{n}}\mathbf{u}\cdot \nabla \mathbf{v}\psi(t)\ast\psi(t)\ast \mathbf{v}(t) dx
\\[0.3cm]
&\displaystyle\qquad\quad-2E(t)\displaystyle\int _{\mathbb{R}^{n}}\mathbf{v}\cdot \nabla \mathbf{u}^{T}\psi(t)\ast\psi(t)\ast \mathbf{v}(t) dx.
\end{array}
\right.
\end{equation}
Integrating \eqref{2.10} over $(s,t)$ concludes \eqref{2.7}. To attain \eqref{2.8}, note that $\mbox{div}~\mathbf{ v}=0$, making the Fourier transform on the both sides of the first equation in \eqref{1.1} with respect to the space variable $x$, then multiplying the resulting equation by $E(t)\widetilde{\psi}^{2}(t)\hat{\mathbf{v}}(t)$, one deduces \eqref{2.8}.\hfill$\Box$\\
\begin{lem}[\cite{Grafakos-1,Grafakos-2,Kato-Ponce}]\label{l2.8}\rm
 Let $\Lambda^{\beta}=(-\Delta)^{\frac{\beta}{2}}$ be the standard Riesz potential of order $\beta\in \mathbb{R}$, $\beta>0$, $1<p,~p_{1},~ p_{2},~p_{3},~p_{4}< \infty$, and $\frac{1}{p}=\frac{1}{p_{1}}+\frac{1}{p_{2}}=\frac{1}{p_{3}}+\frac{1}{p_{4}}$. Then the following bilinear estimate holds for all $f, g\in \mathscr{S}(\mathbb{R}^{n})$:
$$\left\|\Lambda^{\beta}(fg)\right\|_{L^{p}(\mathbb{R}^{n})}\leq C\left\|\Lambda^{\beta}f\right\|_{L^{p_{1}}(\mathbb{R}^{n})}\|g\|_{L^{p_{2}}(\mathbb{R}^{n})}
+C\|f\|_{L^{p_{3}}(\mathbb{R}^{n})}\left\|\Lambda^{\beta}g\right\|_{L^{p_{4}}(\mathbb{R}^{n})}.\eqno \Box $$
\end{lem}
\begin{lem}[\cite{FGO}]\label{l2.9}\rm
  Let $\Lambda^{\beta}=(-\Delta)^{\frac{\beta}{2}}$ be the standard Riesz potential of order $\beta\in \mathbb{R}$, $\beta_{1}, \beta_{2}\in [0,1]$, $\beta=\beta_{1}+ \beta_{2}$, and $p, ~p_{1},~p_{2}\in (1,\infty) $  such
  that $\frac{1}{p}=\frac{1}{p_{1}}+\frac{1}{p_{2}}$. Then the following bilinear estimate holds for all $f,g\in \mathscr{S}(\mathbb{R}^{n})$ with $n\geq 1$:
  $$\left\|\Lambda^{\beta}(fg)-f\Lambda^{\beta}g-g\Lambda^{\beta}f\right\|_{L^{p}(\mathbb{R}^{n})}\leq C\left\|\Lambda^{\beta_{1}}f\right\|_{L^{p_{1}}(\mathbb{R}^{n})}
  \left\|\Lambda^{\beta_{2}}g\right\|_{L^{p_{2}}(\mathbb{R}^{n})}. \eqno \Box$$
  \end{lem}
\indent We now give a nonlocal Sobolev type imbedding result.
\begin{lem}\label{l2.10}\rm
For $0<\beta< 1$ and $n=2,3$,   (1)  the inclusion $\mathcal{D}_{\sigma}\left(\Lambda^{\beta} \right)(\mathbb{R}^{n})\hookrightarrow L^{2}_{\sigma}(\mathbb{R}^{n})$ is compact.\quad (2)  the imbedding $\mathcal{D}_{\sigma}\left(\Lambda^{\gamma} \right)(\mathbb{R}^{n})\hookrightarrow H^{\gamma}_{\sigma}(\mathbb{R}^{n})$ is continuous for all $\gamma\geq 0$.
\end{lem}
{\bf Proof.}\quad It is easy to check it by using standard functional analysis method (see also \cite{Nezza}).\hfill$\Box$\\
\begin{lem} \label{l2.11}\rm
 For $\frac{n}{4}< \beta<1$ with $n=2,3$, let \vspace*{2 ex} $A=\frac{n}{2}+1-2\beta$. Direct calculation gives\\
 \vspace{5pt}
  $(\mbox{I})$\quad $ \beta\geq 1-\beta,~\dfrac{n-2}{2n} <\dfrac{2\beta-1}{n}<\dfrac{1}{n},  ~\dfrac{2\beta-1}{n}=\dfrac{1}{2}-\dfrac{A}{n},$\\
  \vspace{5pt}
  \indent\quad $\dfrac{n}{2}-1<A< 1<\dfrac{n}{2}<2~\mbox{for}~n=3$,\\
  \vspace{5pt}
  \indent\quad $\dfrac{n}{2}-1<A< 1=\dfrac{n}{2}<2~\mbox{for}~n=2$. \\
 \vspace{5pt}
 $\mbox{(II)}$\quad Due to $\frac{n}{2\beta-1}=\frac{2n}{n-2A}$, for $\dfrac{n}{4}\leq \beta<1$ with $n=3$, and $\dfrac{n}{4}< \beta<1$ with $n=2$, we have
 the following fractional Sobolev-type continuous imbedding between $\mathcal{D}\left(\Lambda^{A} \right)(\mathbb{R}^{n})$ and $L^{ \frac{n}{2\beta-1} }(\mathbb{R}^{n})$:
   $$\mathcal{D}\left(\Lambda^{A} \right)(\mathbb{R}^{n})\hookrightarrow L^{ \frac{n}{2\beta-1} }(\mathbb{R}^{n}).\eqno \Box$$
   \end{lem}

The following Lemma concerns the nonlocal version of the known estimates given in Ladyzhenskaya-Shkoller-Seregin \cite{Ladyzhenskaya-1,Ladyzhenskaya-2,Ladyzhenskaya-3,Ladyzhenskaya-4}.
\begin{lem} \label{l2.12}
For $n=2,3$ and $u(x)\in H^{1}_{0}(\mathbb{R}^{n})$, $\forall$ $\varepsilon>0$, the following estimates hold:
$$\|u\|^{2}_{L^{4}(\mathbb{R}^{n})}\leq \varepsilon \|\nabla u\|^{2}_{L^{2}(\mathbb{R}^{n})}+\varepsilon^{-1}\| u\|^{2}_{L^{2}(\mathbb{R}^{n})}\quad \mbox{for}~n=2,\eqno(E-1)$$
 $$\|u\|^{2}_{L^{4}(\mathbb{R}^{n})}\leq 3^{-\frac{1}{4}}\sqrt{2\varepsilon }\|\nabla u\|^{2}_{L^{2}(\mathbb{R}^{n})}+\sqrt{2}(3^{ \frac{5}{2}}\varepsilon)^{ -\frac{1}{6}}\| u\|^{2}_{L^{2}(\mathbb{R}^{n})} \quad \mbox{for}~n=3.\eqno(E-2)$$
 The above inequalities (E-1) and (E-2) can be generalized to the following nonlocal version (fractional power Sobolev-type) estimates.\\
\indent $\heartsuit$ For $\displaystyle\frac{n}{4}<\beta< 1$ and $u\in \mathcal{D}(\Lambda^{\beta} )(\mathbb{R}^{n})$, the following estimates hold:
$$\|u\|^{2}_{L^{4}(\mathbb{R}^{n})}\leq  \varepsilon \left\|\Lambda^{\beta} u\right\|^{2}_{L^{2}(\mathbb{R}^{n})}+\varepsilon^{-1}\| u\|^{2}_{L^{2}(\mathbb{R}^{n})}\quad \mbox{for}~n=2,\eqno(E-3)$$
$$\|u\|^{2}_{L^{4}(\mathbb{R}^{n})}\leq  C(\beta)\varepsilon \left\|\Lambda^{\beta}  u\right\|^{2}_{L^{2}(\mathbb{R}^{n})}+C(\varepsilon) \| u\|^{2}_{L^{2}(\mathbb{R}^{n})}\quad \mbox{for}~n=3.\eqno(E-4)$$
Here, $\varepsilon$, $C(\beta)$ and $C(\varepsilon)$ are constants; $C(s)$ depends only on spatial dimensions and $\beta$, and  $C(\varepsilon)=O(\varepsilon^{-\frac{1}{3}})$.\\
\indent $\heartsuit$ For the critical case $\displaystyle s=\frac{n}{4} $ and $\displaystyle u\in \mathcal{D}\left(\Lambda^{\frac{n}{4}} \right)(\mathbb{R}^{n})$, the following estimates hold:
$$\|u\|^{2}_{L^{4}(\mathbb{R}^{n} )}\leq C_{n}\left(\left\|\Lambda^{\frac{n}{4}}  u\right\|^{2}_{L^{2}(\mathbb{R}^{n})}+\| u\|^{2}_{L^{2}(\mathbb{R}^{n})}\right) \quad \mbox{for}~n=2,3.\eqno(E-5)$$
Here, $C_{n}$ is a constant depending only on space dimensions $n$.\hfill$\Box$
\end{lem}
\section{Non-Uniform Decay }
\allowdisplaybreaks
In this section we consider the non-uniform decay of the Cauchy problem for the Camassa-Holm equations \eqref{1.1}-\eqref{1.2} in $\mathbb{R}^{n}$ if the initial data is assumed only in $L^{2}(\mathbb{R}^{n})$. In particular, if $\mathbf{v}_{0}\in D_{\sigma}\left(\Lambda \right)(\mathbb{R}^{n})$ for $\dfrac{n}{4}< \beta < 1$, and $\mathbf{v}_{0}\in H_{\sigma}^{1}(\mathbb{R}^{n})$ for $\beta=\dfrac{n}{4}$,  then one deduces that the $L^{2}$-norm of the solution to \eqref{1.1}-\eqref{1.2} decays to zero as time $t$ tends to infinity. Unfortunately, we can't determine the decay rate without more information on the initial data. We now formulate the non-uniform decay result as follows.
\begin{thm}\label{t3.1}\rm
For $\dfrac{n}{4}\leq \beta < 1$, $n=2,3$, let $\mathbf{v}$ be the solution to the Cauchy problem \eqref{1.1}-\eqref{1.2} constructed in Proposition \ref{p2.4}. Then
\\[0.3cm]
\indent (I)\quad If $\mathbf{v}_{0}\in D_{\sigma}\left(\Lambda \right)(\mathbb{R}^{n})$ for $\dfrac{n}{4}< \beta < 1$, and if $\mathbf{v}_{0}\in H_{\sigma}^{1}(\mathbb{R}^{n})$ for $\beta=\dfrac{n}{4}$,  then
  $$\lim\limits_{t\rightarrow\infty}\|\mathbf{v}(t)\|_{L^{2} (\mathbb{R}^{n})}=0;$$
\indent (II)\quad If $\mathbf{v}_{0}\in D_{\sigma}\left(\Lambda \right)(\mathbb{R}^{n})$ for $\dfrac{n}{4}< \beta < 1$, and if $\mathbf{v}_{0}\in H_{\sigma}^{1}(\mathbb{R}^{n})$ for $\beta=\dfrac{n}{4}$,  then
$$\lim\limits_{t\rightarrow\infty}\dfrac{1}{t}\displaystyle\int^{t}_{0}\|\mathbf{v}(\tau)\|_{L^{2} (\mathbb{R}^{n})}d\tau=0;$$
\indent (III)\quad If $\mathbf{v}_{0}\in D_{\sigma}\left(\Lambda \right)(\mathbb{R}^{n})$ for $\dfrac{n}{4}< \beta < 1$, and if $\mathbf{v}_{0}\in H_{\sigma}^{1}(\mathbb{R}^{n})$ for $\beta=\dfrac{n}{4}$ with  $\|\mathbf{v}_{0}\|_{L^{2} (\mathbb{R}^{n})}\leq \varepsilon^{*}$  for an $\varepsilon^{*}=\varepsilon^{*}(\alpha,\nu,n)$ sufficiently small,  then there exists no function $G(t,s): ~\mathbb{R}^{+}\times \mathbb{R}^{+}\rightarrow \mathbb{R}^{+}$ admitting the following two properties simultaneously:
 $$(1)~~\|\mathbf{v}(t)\|_{L^{2} (\mathbb{R}^{n})}\leq G\left(t, \|\mathbf{v}_{0}\|_{L^{2} (\mathbb{R}^{n})}\right), ~\mbox{and}~~(2) ~~\mbox{for~ all}~s,~\lim\limits_{t\rightarrow\infty} G(t,s)=0. $$
\end{thm}
 {\bf Proof.} We shall follow the idea introduced in \cite{Masuda,Ogawa}. The idea is to split the energy into low and high frequency parts firsty used in \cite{Masuda}, to use a cut-off function and the generalized energy inequalities, and then to show that both the high and low frequency terms approach zero.
 \\[0.3cm]
\indent We first show (I).
\\[0.3cm]
\indent Due to $\|\mathbf{v}(t)\|_{L^{2} (\mathbb{R}^{n})}=\|\hat{\mathbf{v}}(t)\|_{L^{2} (\mathbb{R}^{n})}$, it suffices to show that $\lim\limits_{t\rightarrow\infty}\|\hat{\mathbf{v}}(t)\|_{L^{2} (\mathbb{R}^{n})}=0$. Splitting the energy into low and high frequency parts gives rise to
\begin{equation}\label{3.1}
\|\hat{\mathbf{v}}\|_{L^{2}(\mathbb{R}^{n})}\leq \|\phi\hat{\mathbf{v}}\|_{L^{2}(\mathbb{R}^{n})}
+\|(1-\phi)\hat{\mathbf{v}}\|_{L^{2}(\mathbb{R}^{n})},\end{equation}
where $\phi=e^{-\nu|\xi|^{2\beta}}$. In the following, we shall divide the proof into two steps.\\
\\
 {\bf Step 1. }\quad Estimate the low frequency part of the energy $\|\hat{\mathbf{v}}\|_{L^{2}(\mathbb{R}^{n})}$.
 \\[0.30cm]
\indent Fix $t$ temporarily, then make the choice of $E=1$ and $\psi(\tau)=\mathcal{F}^{-1}\left[e^{-\nu|\xi|^{2\beta}(t+1-\tau)}\right]$ in \eqref{2.7}. Note that $\psi$ and $\mathcal{F}(\psi)$ are rapidly decreasing functions for $\tau<t+1$, the relation $\hat{\psi}'(\tau)=\nu|\xi|^{2\beta}\hat{\psi}$ assures that the third and fourth terms on the right hand side of \eqref{2.7} add to zero. By Plancherel's theorem, it follows from \eqref{2.7} and $\phi=e^{-\nu|\xi|^{2\beta}}=\hat{\psi}(t)$ that
\begin{equation}\label{3.2}
\left.
\begin{array}{ll}
&\displaystyle\left\|\phi\hat{\mathbf{v}}(t) \right\|^{2}_{L^{2}(\mathbb{R}^{n})}\leq  \left\|e^{-\nu|\xi|^{2\beta}(t-s)}\phi\hat{\mathbf{v}}(s)\right\|^{2}_{L^{2}(\mathbb{R}^{n})}
\\[0.3cm]
&\qquad\qquad\qquad+2\displaystyle\int^{t} _{s}\left|\left\langle\breve{\phi}^{2}\ast \left(\mathbf{u}\cdot \nabla \mathbf{v}+\mathbf{v}\cdot \nabla \mathbf{u}^{T}\right),e^{-2\nu(-\Delta)^{\beta}(t-\tau)}\mathbf{v}(\tau)\right\rangle\right|d\tau.
\end{array}
\right.
\end{equation}
Due to Lemma 2.7, by the aid of  H\"{o}lder's inequality, Young's inequality and Gagliardo-Nirenberg-Sobolev inequality, we have for $\dfrac{n}{4}\leq \beta < 1$
\begin{equation}\label{3.3}
\left.
\begin{array}{ll}
&\left|\left\langle\breve{\phi}^{2}\ast\mathbf{ u}\cdot \nabla \mathbf{v} ,e^{-2\nu(-\Delta)^{\beta}(t-\tau)}\mathbf{v}(\tau)\right\rangle\right|
\\[0.3cm]
&\qquad\leq\left\|\breve{\phi}^{2}\ast \mathbf{u} \cdot \nabla \mathbf{v}\right\| _{L^{2}(\mathbb{R}^{n})}\left\| e^{-2\nu (-\Delta)^{\beta}(t-\tau)}\mathbf{v}(\tau)\right  \| _{L^{2}(\mathbb{R}^{n})}
\\[0.3cm]
&\qquad\leq\left\|\breve{\phi}^{2}\ast \mathbf{u} \right\| _{L^{\infty}(\mathbb{R}^{n})}\left\| \nabla \mathbf{v}\right  \| _{L^{2}(\mathbb{R}^{n})}\left\|   \mathbf{v}\right  \| _{L^{2}(\mathbb{R}^{n})}
\\[0.3cm]
&\qquad\leq\left\|\breve{\phi}^{2} \right\| _{L^{\frac{2n}{n+2\beta}}(\mathbb{R}^{n})}
\left\|\mathbf{u}\right\| _{L^{\frac{2n}{n-2\beta}}(\mathbb{R}^{n})}
\left\| \nabla \mathbf{v}\right  \| _{L^{2}(\mathbb{R}^{n})}\left\|   \mathbf{v}\right  \| _{L^{2}(\mathbb{R}^{n})}
\\[0.3cm]
&\qquad\leq C(\phi)\left\|   \mathbf{v}\right  \| _{L^{2}(\mathbb{R}^{n})}
\left\|\Lambda^{\beta}\mathbf{u}\right\| _{L^{2}(\mathbb{R}^{n})}
\left\| \nabla \mathbf{v}\right  \| _{L^{2}(\mathbb{R}^{n})}.
  \end{array}
\right.
\end{equation}
In the same manner, one deduces
\begin{equation}\label{3.4}
\left.
\begin{array}{ll}
&\left|\left\langle\breve{\phi}^{2}\ast \mathbf{v}\cdot \nabla \mathbf{u}^{T} ,e^{-2\nu(-\Delta)^{\beta}(t-\tau)}\mathbf{v}(\tau)\right\rangle\right|
\\[0.25cm]
&\qquad\leq C(\phi)
 \| \mathbf{v}  \| _{L^{2}(\mathbb{R}^{n})}\left\| \Lambda^{\beta} \mathbf{u} \right \| _{L^{2}(\mathbb{R}^{n})}\| \nabla \mathbf{v} \| _{L^{2}(\mathbb{R}^{n})}.
 \end{array}
\right.
\end{equation}
Thanks to the triangle inequality, H\"{o}lder's inequality, Proposition \ref{p2.4},  \eqref{2.2} and \eqref{3.2}, one achieves
\begin{equation}\label{3.5}
\left.
\begin{array}{ll}
&\left\| \phi \hat{\mathbf{v}}(t)\right \|^{2} _{L^{2}(\mathbb{R}^{n})}
 \lesssim \left\| e^{-\nu|\xi|^{2\beta}(t-s)} \phi \hat{\mathbf{v}}(s)\right\|^{2}
  _{L^{2}(\mathbb{R}^{n})}
  \\[0.25cm]
 & \qquad\qquad\qquad\quad+2C(\phi)C\left(\|\mathbf{v}_{0} \| _{L^{2}(\mathbb{R}^{n})}\right)\left(\displaystyle\int^{t} _{s}\| \Lambda^{\beta} \mathbf{u}  \| _{L^{2}(\mathbb{R}^{n})}^{2}d\tau\right)^{\frac{1}{2}}
  \left(\displaystyle\int^{t} _{s}\| \nabla \mathbf{v}  \| _{L^{2}(\mathbb{R}^{n})}^{2}d\tau\right)^{\frac{1}{2}}.
 \end{array}
\right.
\end{equation}
Since
$$\lim\limits_{t\rightarrow\infty}\left\| e^{-\nu|\xi|^{2\beta}(t-s)} \phi \hat{\mathbf{v}}(s)\right\|^{2}
  _{L^{2}(\mathbb{R}^{n})}=0,$$
letting $t\rightarrow\infty$ gives rise to
\begin{equation*}
\left.
\begin{array}{ll}
\vspace{5pt}
&\limsup\limits_{t\rightarrow\infty}\left\| \phi \hat{\mathbf{v}}(t)\right \|^{2} _{L^{2}(\mathbb{R}^{n})}\\
&\qquad\leq C(\phi)\left(\|\mathbf{v}_{0} \| _{L^{2}(\mathbb{R}^{n})}\right)\left(\displaystyle\int^{\infty} _{s}\left\|  \Lambda^{\beta}\mathbf{u }\right \| _{L^{2}(\mathbb{R}^{n})}^{2}d\tau\right)^{\frac{1}{2}}
  \left(\displaystyle\int^{\infty} _{s}\|\nabla  \mathbf{v } \| _{L^{2}(\mathbb{R}^{n})}^{2}d\tau\right)^{\frac{1}{2}}.
 \end{array}
\right.
\end{equation*}
Recall Proposition \ref{p2.4}, for $\mathbf{v}_{0}\in D_{\sigma}\left(\Lambda \right)(\mathbb{R}^{n})$ with $\dfrac{n}{4}< \beta < 1$, and for $\mathbf{v}_{0}\in H_{\sigma}^{1}(\mathbb{R}^{n})$ with $\beta=\dfrac{n}{4}$, there holds
\begin{equation}\label{3.6}
\left\{
\begin{array}{ll}
\displaystyle\left\| \mathbf{v}\right\|^{2}_{L^{2}(\mathbb{R}^{n})}+\nu\int_{0}^{T}
\| \Lambda^{\beta} \mathbf{v}\|^{2}_{L^{2}(\mathbb{R}^{n})}dt\leq C\left(n,\beta,\alpha,\nu,\| \mathbf{v}_{0}\|_{\mathcal{A}_{1}}\right),
\\[0.3cm]
\displaystyle\left\|\nabla\mathbf{v}\right\|^{2}_{L^{2}(\mathbb{R}^{n})}+\nu\int_{0}^{T}
\| \nabla  \Lambda^{\beta}  \mathbf{v}\|^{2}_{L^{2}(\mathbb{R}^{n})}dt\leq C\left(n,\beta,\alpha,\nu,\| \mathbf{v}_{0}\|_{\mathcal{A}_{2}}\right).
\end{array}
\right.
\end{equation}
Here,  $\displaystyle \mathcal{A}_{1}=L^{2}_{0}(\mathbb{R}^{n}) $, and $\displaystyle \mathcal{A}_{2}=\left\{\begin{array}{ll}
\mathcal{D}\left(\Lambda \right)(\mathbb{R}^{n})~~&\hbox{for}~~\dfrac{n}{4}< \beta < 1,
\\[0.3cm]
H_{0}^{1}(\mathbb{R}^{n})~~&\hbox{for}~~\beta=\dfrac{n}{4}.\end{array}\right.$\\
By interpolation inequality and H\"{o}lder's inequality, \eqref{3.6} yields that
$$\int_{0}^{T}
\| \nabla \mathbf{v}\|^{2}_{L^{2}(\mathbb{R}^{n})}dt\leq\left(\int_{0}^{T}
\| \Lambda^{\beta} \mathbf{v}\|^{2}_{L^{2}(\mathbb{R}^{n})}dt \right)^{\frac{1}{2}}
\left(\int_{0}^{T}
\| \nabla\Lambda^{\beta} \mathbf{v}\|^{2}_{L^{2}(\mathbb{R}^{n})}dt \right)^{\frac{1}{2}}.\eqno(3.6-a)$$
This together with \eqref{2.2} and \eqref{2.3}
implies that $\left\|  \Lambda^{\beta}\mathbf{u} \right \| _{L^{2}(\mathbb{R}^{n})}^{2}$ and $\|  \nabla \mathbf{v}  \| _{L^{2}(\mathbb{R}^{n})}^{2}$ are both integrable on the positive real line. Letting $s\rightarrow\infty$ then gives
\begin{equation}\label{3.7}
 \limsup\limits_{t\rightarrow\infty}\left\| \phi \hat{\mathbf{v}}(t) \right\|^{2} _{L^{2}(\mathbb{R}^{n})}=0.\end{equation}
{\bf Step 2}\quad We now estimate the high-frequency part of  the energy $\|\mathbf{v}(t)\|_{L^{2} (\mathbb{R}^{n})}$.
\\[0.30cm]
  \indent Put $\widetilde{\psi}=1-e^{-\nu|\xi|^{2\beta}}=1-\phi$ in \eqref{2.8}. Let $B_{G}(t)=\left\{\xi:~|\xi|\leq G(t)\right\}$, where $G(t)$ will be determined later. Note that $\langle \mathbf{u}\cdot \nabla \mathbf{v},\mathbf{v}\rangle=0$, replacing $\widetilde{\psi}^{2}$ by $1-\widetilde{\psi}^{2}$ in the fourth term on the right hand side of \eqref{2.8} yields
\begin{equation}\label{3.8}
\left.
\begin{array}{ll}
\vspace{5pt}
& E(t)\left\|\widetilde{\psi}(t)  \hat{\mathbf{v}}(t)\right\|^{2}_{L^{2}(\mathbb{R}^{n})}\\
\vspace{5pt}
&\qquad\leq E(s)\left\|\widetilde{\psi}(s) \hat{\mathbf{v}}(s)\right\|^{2}_{L^{2}(\mathbb{R}^{n})}+\displaystyle\int^{t}_{s}
E'(\tau)\int _{B_{G}(\tau)}\left|\widetilde{\psi}(\tau)\hat{\mathbf{v}}(\tau)\right |^{2}d\xi d\tau\\
\vspace{5pt}
&\qquad\quad +\displaystyle\int^{t}_{s}\left(E'(\tau)-2\nu E(\tau)G^{2\beta}(\tau)\right)\int _{B^{c}_{G}(\tau)}\left|\widetilde{\psi}(\tau)\hat{\mathbf{v}}(\tau)\right |^{2}d\xi d\tau\\
\vspace{5pt}
&\qquad\quad+2 \displaystyle\int^{t}_{s}E(\tau)
\left|\left\langle \mathcal{F}\left( \mathbf{u}\cdot \nabla \mathbf{v}+\mathbf{v}\cdot \nabla \mathbf{u}^{T}\right), \left(1-\widetilde{\psi}^{2}(\tau)\right) \hat{ \mathbf{v}}(\tau)\right\rangle\right| d\tau\\
\vspace{5pt}
&\qquad\quad+2 \displaystyle\int^{t}_{s}E(\tau)
\left|\left\langle \mathcal{F}\left(\mathbf{v}\cdot \nabla \mathbf{u}^{T}\right),  \hat{\mathbf{v}}(\tau)\right\rangle \right| d\tau.
\end{array}
\right.
\end{equation}
Since
\begin{equation*}
\left.
\begin{array}{ll}
\vspace{5pt}
\psi(\tau)=\mathcal{F}^{-1}\left[e^{-\nu|\xi|^{2\beta}(t+1-\tau)}\right],
~~\phi=e^{-\nu|\xi|^{2\beta}}=\hat{\psi}(t),\\
\vspace{5pt}
\widetilde{\psi}=1-\phi,~~
\widetilde{\psi}(\tau)=\mathcal{F}^{-1}\left[1-e^{-\nu|\xi|^{2\beta}(t+1-\tau)}\right],
 \end{array}
\right.
\end{equation*}
  $\mathcal{F}(\phi)=1-\widetilde{\psi}^{2}$ and $\phi=1-\widetilde{\psi}$ are rapidly decreasing functions,
applying H\"{o}lder's inequality, the Plancherel's Theorem, Young's inequality and Gagliardo-Nirenberg-Sobolev inequality, we obtain
\begin{equation}\label{3.9}
\left.
\begin{array}{ll}
& \left|\left\langle \mathcal{F}\left(\mathbf{u}\cdot \nabla \mathbf{v}+\mathbf{v}\cdot \nabla \mathbf{u}^{T}\right),\left(1-\widetilde{\psi}^{2}(\tau)\right)\hat{\mathbf{v}}(\tau)\right\rangle\right|
\\[0.3cm]
&\qquad=\left|\left\langle\left(1-\widetilde{\psi}^{2}(\tau)\right) \mathcal{F}\left(\mathbf{u}\cdot \nabla \mathbf{v}+\mathbf{v}\cdot \nabla \mathbf{u}^{T}\right),\hat{\mathbf{v}}(\tau)\right\rangle\right|
\\[0.3cm]
& \qquad\leq \left\|1-\widetilde{\psi}^{2}(\tau)\right\|_{L^{2}(\mathbb{R}^{n})}
\left\|\hat{\mathbf{v}}(\tau)\right\|_{L^{2}(\mathbb{R}^{n})}
\\[0.3cm]
&\qquad\qquad
\cdot\left(\left\|\mathcal{F}\left(\mathbf{u}\cdot \nabla \mathbf{v}\right)\right\|_{L^{\infty}(\mathbb{R}^{n})}+\left\|\mathcal{F}\left(\mathbf{v}\cdot \nabla \mathbf{u}^{T}\right)\right\|_{L^{\infty}(\mathbb{R}^{n})}\right)
\\[0.3cm]
&\qquad\leq C(\phi)
\left(\left\|\hat{ \mathbf{u}}\ast \xi  \hat{\mathbf{v}}\right\|_{L^{\infty}(\mathbb{R}^{n})}+\left\| \hat{\mathbf{v}}\ast \xi \hat{\mathbf{u}}\right\|_{L^{\infty}(\mathbb{R}^{n})}\right)
\|\mathbf{v}\|_{L^{2}(\mathbb{R}^{n})}
\\[0.3cm]
&\qquad\leq C(\phi)\|\mathbf{v}\|_{L^{2}(\mathbb{R}^{n})}
\left(\left\| \hat{ \mathbf{u}}\ast \xi  \hat{\mathbf{v}}\right\|_{L^{\infty}(\mathbb{R}^{n})}+\left\|  \hat{\mathbf{v}}\ast \xi  \hat{\mathbf{u}}\right\|_{L^{\infty}(\mathbb{R}^{n})}\right)
\\[0.3cm]
&\qquad\leq C(\phi)\|\mathbf{v}\|_{L^{2}(\mathbb{R}^{n})}\left(\| \mathbf{ u} \|_{L^{2}(\mathbb{R}^{n})}
 \| \nabla \mathbf{v} \|_{L^{2}(\mathbb{R}^{n})}+\|\mathbf{v}\|_{L^{2}(\mathbb{R}^{n})}\| \nabla \mathbf{u} \|_{L^{2}(\mathbb{R}^{n})} \right).
 \end{array}
\right.
\end{equation}
In the same manner, one concludes
\begin{equation}\label{3.10}
\left.
\begin{array}{ll}
\vspace{5pt}
& \left| \left\langle \mathcal{F}\left(\mathbf{v}\cdot \nabla \mathbf{u}^{T}\right),\hat{\mathbf{v}}(\tau)\right\rangle\right|\leq\left\| \mathbf{v}\cdot \nabla \mathbf{u}^{T} \right\|_{L^{2}(\mathbb{R}^{n})}\|\mathbf{v}\|_{L^{2}(\mathbb{R}^{n})}\\
\vspace{5pt}
&\qquad\qquad\qquad\qquad \leq C\| \mathbf{v}\|_{L^{\frac{2n}{n-2\beta}}(\mathbb{R}^{n})}\| \nabla \mathbf{u}\|_{L^{\frac{n}{\beta}}(\mathbb{R}^{n})}\|\mathbf{v}\|_{L^{2}(\mathbb{R}^{n})}\\
\vspace{5pt}
&\qquad\qquad\qquad\qquad \leq C\left\| \Lambda ^{\beta}\mathbf{v}\right\|^{2}_{L^{2}(\mathbb{R}^{n})}\| \mathbf{v}\|_{L^{2}(\mathbb{R}^{n})}.
\end{array}
\right.
\end{equation}
Choosing $E(t)=(1+t)^{k}$, and $G^{2\beta}(t)=\dfrac{k}{2\nu (1+t)}$ in \eqref{3.8} such that $E'-2\nu EG^{2\beta}=0$, then taking $k>0$ sufficiently large yields
\begin{eqnarray*}
 && \left\|(1-\phi)\hat{\mathbf{v}}(\tau)\right\|^{2}_{L^{2}(\mathbb{R}^{n})}\\
  &&\qquad\leq \frac{(1+s)^{k}}{(1+t)^{k}} \left\|(1-\phi)\hat{\mathbf{v}}(s)\right\|^{2}_{L^{2}(\mathbb{R}^{n})}\\
 &&\qquad \quad+\displaystyle\int^{t}_{s}\frac{k(1+\tau)^{k-1}}{(1+t)^{k}}
 \displaystyle\int_{B_{G}(\tau)}\left|(1-\phi)\hat{\mathbf{v}}(\tau)\right|^{2}d\xi d\tau\\
 &&\qquad\quad+C\|\mathbf{v}_{0}\|_{L^{2}(\mathbb{R}^{n})}\displaystyle\int^{t}_{s}
 \frac{(1+\tau)^{k}}{(1+t)^{k}}\cdot \left( \|\mathbf{\mathbf{u}}\|_{L^{2}(\mathbb{R}^{n})}\|\nabla \mathbf{v}\|_{L^{2}(\mathbb{R}^{n})}\right.
 \\
&&\qquad\qquad\qquad\qquad\qquad\qquad \left.+ \|\mathbf{\mathbf{v}}\|_{L^{2}(\mathbb{R}^{n})}\|\nabla \mathbf{u}\|_{L^{2}(\mathbb{R}^{n})}+ \left\|\Lambda^{\beta}\mathbf{v}\right\|^{2}_{L^{2}(\mathbb{R}^{n})}\right)d\tau.\qquad\qquad
 \end{eqnarray*}
For $\xi\in B_{G}(t)$ and $t$ sufficiently large, there holds that $\widetilde{\psi}=|1-\phi|\leq \nu |\xi|^{2\beta} $. In particular, $\displaystyle|1-\phi|^{2}\leq \frac{k^{2}}{4(1+t)^{2}}$. Thus the second term on the right hand side of the above inequality is bounded as follows:
\begin{eqnarray*}
  &&  \displaystyle\int^{t}_{s}\frac{k^{3}(1+\tau)^{-3}}{4}
 \displaystyle\int_{B_{G}(\tau)}| \hat{\mathbf{v}}(\tau)|^{2}d\xi d\tau
 \\[0.1cm]
  && \qquad\qquad\leq\displaystyle\int^{t}_{s}\frac{k^{3}(1+\tau)^{-3}}{4}
  \| \mathbf{v}(\tau)\|^{2}_{L^{2}(\mathbb{R}^{n})}d\tau
  \\[0.1cm]
 &&\qquad\qquad \leq\frac{k^{3} }{4} \| \mathbf{v}_{0}\|^{2}_{L^{2}(\mathbb{R}^{n})}\displaystyle\int^{t}_{s} (1+\tau)^{-3} d\tau
 \\[0.1cm]
 &&\qquad\qquad \leq\frac{k^{3} }{8} \| \mathbf{v}_{0}\|^{2}_{L^{2}(\mathbb{R}^{n})}  (1+s)^{-2} .
  \end{eqnarray*}
Letting $t\rightarrow \infty$ gives rise to
\begin{equation}\label{3.11}
\left.
\begin{array}{ll}
\vspace{5pt}
& \limsup\limits_{t\rightarrow \infty} \left\|(1-\phi) \hat{\mathbf{v}}(\tau)\right\|^{2}_{L^{2}(\mathbb{R}^{n})}\\
\vspace{5pt}
&\qquad\leq\dfrac{k^{3} }{8} \| \mathbf{v}_{0}\|^{2}_{L^{2}(\mathbb{R}^{n})}  (1+s)^{-2}\\
 \vspace{5pt}
&\qquad\quad
 +C\| \mathbf{v}_{0}\|_{L^{2}(\mathbb{R}^{n})} \left(\displaystyle\int^{\infty}_{s} \|\nabla \mathbf{v}\|^{2}_{L^{2}(\mathbb{R}^{n})} d\tau+\displaystyle\int^{\infty}_{s} \left\| \Lambda^{\beta} \mathbf{v}\right\|^{2}_{L^{2}(\mathbb{R}^{n})} d\tau
 +\displaystyle\int^{\infty}_{s} \| \nabla \mathbf{u}\|^{2}_{L^{2}(\mathbb{R}^{n})} d\tau
 \right).
 \end{array}
\right.
\end{equation}
Recall Proposition \ref{p2.4},  for $\mathbf{v}_{0}\in D_{\sigma}\left(\Lambda \right)(\mathbb{R}^{n})$ with $\dfrac{n}{4}< \beta < 1$, and for $\mathbf{v}_{0}\in H_{\sigma}^{1}(\mathbb{R}^{n})$ with $\beta=\dfrac{n}{4}$, there holds
\begin{equation}\label{3.12}
\| \mathbf{v}\|^{2}_{L^{2}(\mathbb{R}^{n})}+\nu\int_{0}^{T}
\left\| \Lambda^{\beta} \mathbf{v}\right\|^{2}_{L^{2}(\mathbb{R}^{n})}dt\leq C\left(n,\beta,\alpha,\nu,\| \mathbf{v}_{0}\|_{H^{1}_{0}(\mathbb{R}^{n})}\right).\end{equation}
Thanks to \eqref{2.2}, \eqref{3.6} and \eqref{3.12}, by interpolation inequality, one deduces
 that $\left\|\Lambda^{\beta} \mathbf{v}\right\|^{2}_{L^{2}(\mathbb{R}^{n})}$,~~$\|\nabla \mathbf{v}\|^{2}_{L^{2}(\mathbb{R}^{n})}$ and $\|\nabla \mathbf{u}\|^{2}_{L^{2}(\mathbb{R}^{n})}$
are all integrable on the real line. Letting $s\rightarrow\infty$ gives
$$\limsup\limits_{t\rightarrow \infty} \left\|(1-\phi) \hat{\mathbf{v}}(\tau)\right\|^{2}_{L^{2}(\mathbb{R}^{n})}=0.$$
Combining this with \eqref{3.7} and the Plancherel's theorem finishes the proof of (I).
\\[0.3cm]
{\bf We next show (II)}.
\\[0.3cm]
\indent According to (I), given an $\epsilon>0$ we can choose $s$ large enough such that $\|\mathbf{v}\| _{L^{2}(\mathbb{R}^{n})}\leq \varepsilon$ for $\tau>s$. Thus there holds
\begin{equation}\label{3.13}
\left.
\begin{array}{ll}
\vspace{5pt}
&\displaystyle \frac{1}{t}\int^{t}_{0}\|\mathbf{v}(\tau)\| _{L^{2}(\mathbb{R}^{n})}d\tau\\
 \vspace{5pt}
 &\qquad =\displaystyle\frac{1}{t}\int^{s}_{0}\|\mathbf{v}(\tau)\| _{L^{2}(\mathbb{R}^{n})}d\tau+\frac{1}{t}\displaystyle\int^{t}_{s}\|\mathbf{v}(\tau)\| _{L^{2}(\mathbb{R}^{n})}d\tau\\
 \vspace{5pt}
 &\qquad \leq\displaystyle\frac{1}{t} \int^{s}_{0}\|\mathbf{v}(\tau)\| _{L^{2}(\mathbb{R}^{n})}d\tau+\varepsilon\frac{t-s}{t}.
 \end{array}
\right.
\end{equation}
Since $\varepsilon$ can be chosen arbitrarily, letting $t\rightarrow \infty$
   finishes the proof of (II).
   \\[0.3cm]
{\bf We are now in the position to show (III)}.
\\[0.3cm]
   \indent Let $\mathbf{u}_{0}(x)$ be any smooth function with compact support, and $\displaystyle\mathbf{u}^{\varepsilon}_{0}(x)=\varepsilon^{\frac{n}{2}}\mathbf{u}_{0}(\varepsilon x)$. In addition, let
   $\mathbf{v}^{\varepsilon}_{0}=\mathbf{u}^{\varepsilon}_{0}-\alpha^{2}\Delta \mathbf{u}^{\varepsilon}_{0}$
  and $\mathbf{v}^{\varepsilon}$ be the solution of \eqref{1.1}-\eqref{1.2} given by Proposition \ref{p2.4} corresponding to the initial data $\mathbf{v}_{0}$. For any $\varepsilon>0$, a straightforward computation shows that
  \begin{equation}\label{3.14}
  \left\|\mathbf{u}^{\varepsilon}_{0}\right\| _{L^{2}(\mathbb{R}^{n})}=\|\mathbf{u}_{0}\| _{L^{2}(\mathbb{R}^{n})},~ ~ \left\|\nabla ^{m}\mathbf{u}^{\varepsilon}_{0}\right\| _{L^{2}(\mathbb{R}^{n})}=\varepsilon^{m} \left\|\nabla \mathbf{u}_{0}\right\| _{L^{2}(\mathbb{R}^{n})},\end{equation}
 \begin{equation}\label{3.15}
\left.
\begin{array}{ll}
 \displaystyle\left\|\mathbf{v}^{\varepsilon}_{0}\right\| ^{2}_{L^{2}(\mathbb{R}^{n})} &\displaystyle=\left\|\mathbf{u}^{\varepsilon}_{0}\right\| ^{2}_{L^{2}(\mathbb{R}^{n})}+\alpha^{2}\left\|\nabla \mathbf{u}^{\varepsilon}_{0}\right\| ^{2}_{L^{2}(\mathbb{R}^{n})}+\alpha^{4} \left\|\Delta \mathbf{u}^{\varepsilon}_{0}\right\| ^{2}_{L^{2}(\mathbb{R}^{n})}
 \\[0.3cm]
&  =\displaystyle\|\mathbf{u} _{0}\| ^{2}_{L^{2}(\mathbb{R}^{n})}+\alpha^{2}\varepsilon^{2}\left\|\nabla \mathbf{u} _{0}\right\| ^{2}_{L^{2}(\mathbb{R}^{n})}+\alpha^{4}\varepsilon^{4}\left\|\Delta \mathbf{u} _{0}\right\| ^{2}_{L^{2}(\mathbb{R}^{n})},
 \end{array}
\right.
\end{equation}
  and
 \begin{equation}\label{3.16}
\left.
\begin{array}{ll}
\left\|\nabla \mathbf{v}^{\varepsilon}_{0}\right\| ^{2}_{L^{2}(\mathbb{R}^{n})}&=\left\|\nabla \mathbf{u}^{\varepsilon}_{0}\right\| ^{2}_{L^{2}(\mathbb{R}^{n})}
  +2\alpha^{2}\left\|\nabla^{2} \mathbf{u}^{\varepsilon}_{0}\right\| ^{2}_{L^{2}(\mathbb{R}^{n})}+\alpha^{4} \left\|\nabla\Delta \mathbf{u}^{\varepsilon}_{0}\right\| ^{2}_{L^{2}(\mathbb{R}^{n})}
   \\[0.3cm]
 &\displaystyle=\varepsilon^{2}\left\|\nabla^{2} \mathbf{u}_{0}\right\| ^{2}_{L^{2}(\mathbb{R}^{n})}
  +2\alpha^{2}\varepsilon^{4}\left\|\Delta \mathbf{u}_{0}\right\| ^{2}_{L^{2}(\mathbb{R}^{n})}+\alpha^{4} \varepsilon^{6}\left\|\nabla\Delta \mathbf{u}_{0}\right\| ^{2}_{L^{2}(\mathbb{R}^{n})}.\end{array}
\right.
\end{equation}
It follows from \eqref{3.15}, \eqref{3.16} and Corollary \ref{c2.5} that there exists a constant $C=C\left(\| \mathbf{u} _{0}\| ^{2}_{H^{3}_{\sigma}(\mathbb{R}^{n})}\right)$ such that for all $\varepsilon>0$,
  \begin{equation}\label{3.17}
 \left \|  \mathbf{v}^{\varepsilon} _{0} \right\| ^{2}_{L^{2}(\mathbb{R}^{n})}\leq C,~~\left\| \nabla  \mathbf{v}^{\varepsilon} _{0} \right\| ^{2}_{L^{2}(\mathbb{R}^{n})}\leq C\varepsilon^{2}.
  \end{equation}
 We then claim
  \begin{equation}\label{3.18}
   \frac{d}{dt}\left(\left\| \mathbf{u}^{\varepsilon} \right\| ^{2}_{L^{2}(\mathbb{R}^{n})}+\alpha^{2}\left\| \nabla \mathbf{u}^{\varepsilon} \right\| ^{2}_{L^{2}(\mathbb{R}^{n})}\right)\geq -C\varepsilon^{2},
   \end{equation}
  which is equivalent to
  \begin{equation}\label{3.19}
\left.
\begin{array}{ll}
&\left\| \mathbf{u}^{\varepsilon} \right\| ^{2}_{L^{2}(\mathbb{R}^{n})}+\alpha^{2}\left\| \nabla \mathbf{u}^{\varepsilon} \right\| ^{2}_{L^{2}(\mathbb{R}^{n})}
\\[0.3cm]
  &\qquad \geq\left\|\mathbf{u}^{\varepsilon}_{0} \right\| ^{2}_{L^{2}(\mathbb{R}^{n})}+\alpha^{2}\left\| \nabla \mathbf{u}^{\varepsilon}_{0} \right\| ^{2}_{L^{2}(\mathbb{R}^{n})}-C\varepsilon^{2}t
  \\[0.3cm]
 &\qquad=\left\| \mathbf{u}_{0} \right\| ^{2}_{L^{2}(\mathbb{R}^{n})}+\alpha^{2}\varepsilon^{2}\left\| \nabla \mathbf{u}_{0} \right\| ^{2}_{L^{2}(\mathbb{R}^{n})}-C\varepsilon^{2}t
 \\[0.3cm]
&\qquad\geq \left\| \mathbf{u}_{0} \right\| ^{2}_{L^{2}(\mathbb{R}^{n})} -C\varepsilon^{2}t.
\end{array}
\right.
\end{equation}
Thanks to \eqref{3.18} and \eqref{3.19}, we conclude that there is not a function $G(t,s)$ continuous and approaching zero in $t$ for each fixed $s$, such that
\begin{equation}\label{3.20}
\| \mathbf{v}\|_{L^{2}(\mathbb{R}^{n})}\leq G\left(t,\| \mathbf{u}_{0} \|_{L^{2}(\mathbb{R}^{n})}\right).
\end{equation}
Otherwise, if there were such a function, then at some $t_{0}$ it would admit the bound
 \begin{equation}\label{3.21}
  G\left(t_{0},\| \mathbf{u}_{0} \|_{L^{2}(\mathbb{R}^{n})}\right)\leq \frac{1}{2}\| \mathbf{u}_{0} \|_{L^{2}(\mathbb{R}^{n})}.
  \end{equation}
 Choosing $\varepsilon$ sufficiently small in \eqref{3.19}, in particular, $\displaystyle\varepsilon^{2\beta}\leq \frac{\| \mathbf{u}_{0} \|^{2}_{L^{2}(\mathbb{R}^{n})}}{4Ct_{0}} $, one deduces that
  $$  G\left(t_{0},\| \mathbf{u}_{0} \|^{2}_{L^{2}(\mathbb{R}^{n})}\right)\geq \| \mathbf{v} \|^{2}_{L^{2}(\mathbb{R}^{n})} \geq \frac{3}{4}\| \mathbf{u}_{0} \|^{2}_{L^{2}(\mathbb{R}^{n})}.$$
  This is contradictory to \eqref{3.21}. \\
  \indent Once we have shown \eqref{3.18} or \eqref{3.19}, the proof of (III) will be finished.\\
  \indent We are now in the position to show \eqref{3.18}.\\
  \indent Note that $\mathbf{v}^{\varepsilon}$ is a solution of \eqref{1.1}-\eqref{1.2}, multiplying the first equation for $\mathbf{v}^{\varepsilon}$ in (1.1) by $\Delta \mathbf{v}^{\varepsilon}$, then integrating by parts yields
  \begin{equation}\label{3.22}
  \frac{1}{2}\frac{d}{dt}\left\| \nabla  \mathbf{v}^{\varepsilon} \right\| ^{2}_{L^{2}(\mathbb{R}^{n})}+\nu\left\| \Lambda^{\beta}\nabla  \mathbf{v}^{\varepsilon} \right\| ^{2}_{L^{2}(\mathbb{R}^{n})}=\left\langle \mathbf{u}^{\varepsilon}\cdot\nabla \mathbf{v}^{\varepsilon}, \Delta \mathbf{v}^{\varepsilon}\right\rangle+\left\langle \mathbf{v}^{\varepsilon}\cdot\nabla \mathbf{u}^{\varepsilon T}, \Delta \mathbf{v}^{\varepsilon}\right\rangle.
  \end{equation}
  \\[0.2cm]
  We then deal with the two terms on the right hand side of \eqref{3.22} through two cases:\\
  \\
   \indent  {\bf Case (I)}\quad $\dfrac{n}{4}< \beta< 1$ for $n=2,3$;
   \\[0.2cm]
   \indent {\bf Case (II)}\quad $ \beta=\dfrac{n}{4}$ for $n=2,3$.
    \\[0.3cm]
We first consider Case {\bf (I)}\quad $\dfrac{n}{4}< \beta< 1$ for $n=2,3$.
 \\[0.3cm]
  \indent In this case, notice that $\left\langle \mathbf{u}^{\varepsilon}\cdot\nabla \mathbf{v}^{\varepsilon}, \mathbf{v}^{\varepsilon}\right\rangle=0$, ~~ $\displaystyle\frac{\beta}{n}=\frac{1}{2}-\frac{n/2-\beta}{n}$, ~$\displaystyle\frac{n}{2}-1<\frac{n}{2}-\beta<\frac{n}{4}$,~ $\displaystyle\frac{n}{2}-1<\frac{n}{2}-2\beta+1<1$, H\"{o}lder's inequality, Sobolev inequality and  Cauchy-Schwartz inequality yield that
 \begin{equation}\label{3.23}
\left.
\begin{array}{ll}
&\quad\displaystyle\left|\left\langle \mathbf{u}^{\varepsilon}\cdot\nabla \mathbf{v}^{\varepsilon}, \Delta \mathbf{v}^{\varepsilon}\right\rangle\right|=\left| \left\langle \nabla \mathbf{u}^{\varepsilon}\cdot\nabla \mathbf{v}^{\varepsilon}, \nabla \mathbf{v}^{\varepsilon}\right\rangle\right|
\\[0.3cm]
 &\displaystyle\qquad\qquad\qquad\quad\leq C\left\| \nabla   \mathbf{ v}^{\varepsilon} \right\| _{L^{2}(\mathbb{R}^{n})}
 \left \| \nabla \mathbf{u}^{\varepsilon}\nabla v^{\varepsilon} \right\| _{L^{2}(\mathbb{R}^{n})}
  \\[0.3cm]
 &\displaystyle\qquad\qquad\qquad\quad\leq C\left\| \nabla    \mathbf{v}^{\varepsilon} \right\| _{L^{2}(\mathbb{R}^{n})}
  \left\| \nabla \mathbf{v}^{\varepsilon} \right\| _{L^{\frac{2n}{n-2\beta}}(\mathbb{R}^{n})} \left\| \nabla \mathbf{u}^{\varepsilon} \right\| _{L^{\frac{n}{\beta}}(\mathbb{R}^{n})}
  \\[0.3cm]
 &\displaystyle\qquad\qquad\qquad\quad\leq C\left\| \nabla    \mathbf{v}^{\varepsilon}\right \| _{L^{2}(\mathbb{R}^{n})}
  \left\| \Lambda^{\beta}\nabla \mathbf{v}^{\varepsilon} \right\| _{L^{2}(\mathbb{R}^{n})}\left \| \mathbf{v}^{\varepsilon} \right\| _{L^{2}(\mathbb{R}^{n})}
  \\[0.3cm]
 &\displaystyle\qquad\qquad\qquad\quad\displaystyle\leq \frac{\nu}{4}\left\|\Lambda^{\beta} \nabla \mathbf{v}^{\varepsilon} \right\| ^{2}_{L^{2}(\mathbb{R}^{n})}+C\left\| \mathbf{v}^{\varepsilon} \right\| ^{2}_{L^{2}(\mathbb{R}^{n})}\left\| \nabla \mathbf{v}^{\varepsilon}\right \| ^{2}_{L^{2}(\mathbb{R}^{n})}.
\end{array}
\right.
\end{equation}
On the other hand, Lemma \ref{l2.8} and  Lemma \ref{l2.11} ensure that
  \begin{equation}\label{3.24}
\left.
\begin{array}{ll}
&\left|\left\langle \mathbf{v}^{\varepsilon}\cdot\nabla \mathbf{u}^{\varepsilon T}, \Delta \mathbf{v}^{\varepsilon}\right\rangle\right|
\\[0.3cm]
  &\qquad \leq\left|\left\langle\Lambda^{1-\beta} \left(\mathbf{v}^{\varepsilon}\cdot\nabla \mathbf{u}^{\varepsilon T}\right),  \Lambda^{\beta}\nabla \mathbf{v}^{\varepsilon}\right\rangle\right|
  \\[0.3cm]
   &\qquad \leq\left\|\Lambda^{1-\beta} \left(\mathbf{v}^{\varepsilon}\cdot\nabla \mathbf{u}^{\varepsilon T}\right) \right\|_{L^{2}(\mathbb{R}^{n})}\left\|\Lambda^{\beta}\nabla \mathbf{v}^{\varepsilon} \right\|_{L^{2}(\mathbb{R}^{n})}
   \\[0.3cm]
    &\qquad\displaystyle \leq C\left(\left\|\Lambda^{1-\beta} \mathbf{v}^{\varepsilon}\right \|_{L^{ \frac{2n}{n-4\beta+2} }(\mathbb{R}^{n})}\left\|\nabla \mathbf{u}^{\varepsilon} \right\|_{L^{\frac{2n}{4\beta-2} }(\mathbb{R}^{n})}\right.
    \\[0.3cm]
   &\qquad \qquad\quad\displaystyle \left.+\left\| \mathbf{v}^{\varepsilon} \right\|_{L^{ \frac{2n}{n-2\beta} }(\mathbb{R}^{n})}\left\|\nabla \mathbf{u}^{\varepsilon} \right\|_{L^{ \frac{n}{\beta} }(\mathbb{R}^{n})}\right)\cdot\left \|\Lambda^{\beta}\nabla \mathbf{v}^{\varepsilon} \right\|_{L^{2}(\mathbb{R}^{n})}
   \\[0.3cm]
    &\qquad \leq \displaystyle\frac{\nu}{4}\left\|\Lambda^{\beta} \nabla \mathbf{v}^{\varepsilon} \right \| ^{2}_{L^{2}(\mathbb{R}^{n})}+C\left\| \Lambda^{\beta}\mathbf{v}^{\varepsilon} \right\| ^{2}_{L^{2}(\mathbb{R}^{n})}\left\| \Lambda^{ \frac{n}{2}-2\beta+1}\nabla \mathbf{u}^{\varepsilon} \right\| ^{2}_{L^{2}(\mathbb{R}^{n})}
    \\[0.3cm]
  &\qquad \leq \displaystyle\frac{\nu}{4}\left\|\Lambda^{\beta} \nabla \mathbf{v}^{\varepsilon} \right\| ^{2}_{L^{2}(\mathbb{R}^{n})}+C\left\| \Lambda^{\beta}\mathbf{v}^{\varepsilon}\right \| ^{2}_{L^{2}(\mathbb{R}^{n})}\| \mathbf{v}^{\varepsilon} \| ^{2}_{L^{2}(\mathbb{R}^{n})}.
  \end{array}
\right.
\end{equation}
  Combining \eqref{3.22} with \eqref{3.23} and \eqref{3.24} gives rise to
   \begin{equation}\label{3.25}
   \frac{d}{dt}\left\| \nabla  \mathbf{v}^{\varepsilon} \right\| ^{2}_{L^{2}(\mathbb{R}^{n})}+\frac{\nu}{2}\left\| \Lambda^{\beta}\nabla  \mathbf{v}^{\varepsilon} \right\| ^{2}_{L^{2}(\mathbb{R}^{n})}\leq C\left\| \nabla \mathbf{v}^{\varepsilon} \right\| ^{2}_{L^{2}(\mathbb{R}^{n})}\left\|  \mathbf{v}^{\varepsilon} \right\| ^{2}_{L^{2}(\mathbb{R}^{n})}.\end{equation}
   We next consider Case {\bf  (II)}\quad $ \beta=\dfrac{n}{4}$ for $n=2,3$.
   \\[0.3cm]
    \indent In this case, thanks to $\left\langle \mathbf{u}^{\varepsilon}\cdot\nabla \mathbf{v}^{\varepsilon}, \mathbf{v}^{\varepsilon}\right\rangle=0$ and $\frac{1}{4}=\frac{1}{2}-\frac{\beta}{n}$, note that Lemma \ref{l2.12}, applying H\"{o}lder's inequality and Gagliardo-Nirenberg-Sobolev inequality imply
  \begin{equation}\label{3.26}
\left.
\begin{array}{ll}
&\left|\left\langle \mathbf{u}^{\varepsilon}\cdot\nabla \mathbf{v}^{\varepsilon}, \Delta \mathbf{v}^{\varepsilon}\right\rangle\right|\leq\left| \left\langle \nabla \mathbf{u}^{\varepsilon}\cdot\nabla \mathbf{v}^{\varepsilon}, \nabla \mathbf{v}^{\varepsilon}\right\rangle\right|
\\[0.3cm]
 &\qquad\qquad\qquad\quad\leq\left\| \nabla   \mathbf{u}^{\varepsilon} \right\| _{L^{2}(\mathbb{R}^{n})}
  \left\| \nabla \mathbf{v}^{\varepsilon} \right\| ^{2}_{L^{4}(\mathbb{R}^{n})}
  \\[0.3cm]
 &\qquad\qquad\qquad\quad\leq C\left\| \nabla    \mathbf{u}^{\varepsilon} \right\| _{L^{2}(\mathbb{R}^{n})}
  \left(\left\| \nabla \mathbf{v}^{\varepsilon} \right\| ^{2}_{L^{2}(\mathbb{R}^{n})}+\left\| \Lambda^{\frac{n}{4}}\nabla \mathbf{v}^{\varepsilon} \right\| ^{2}_{L^{2}(\mathbb{R}^{n})}\right),
     \end{array}
\right.
\end{equation}
and
  \begin{equation}\label{3.27}
\left.
\begin{array}{ll}
&\left|\left\langle \mathbf{v}^{\varepsilon}\cdot\nabla \mathbf{u}^{\varepsilon T}, \Delta \mathbf{v}^{\varepsilon}\right\rangle\right|
\\[0.3cm]
 &\qquad \leq\left\| \nabla \mathbf{v}^{\varepsilon} \right\| ^{2}_{L^{4}(\mathbb{R}^{n})}\left\| \nabla \mathbf{u}^{\varepsilon} \right\|_{L^{2}(\mathbb{R}^{n})}+\left\| \Delta \mathbf{u}^{\varepsilon} \right\|_{L^{2}(\mathbb{R}^{n})}\left\| \mathbf{v}^{\varepsilon}\nabla v^{\varepsilon} \right\| _{L^{2}(\mathbb{R}^{n})}
 \\[0.3cm]
 &\qquad \leq\left\| \nabla \mathbf{v}^{\varepsilon} \right\| ^{2}_{L^{4}(\mathbb{R}^{n})}\left\| \nabla \mathbf{u}^{\varepsilon} \right\|_{L^{2}(\mathbb{R}^{n})}+\left\| \Delta \mathbf{u}^{\varepsilon} \right\|_{L^{2}(\mathbb{R}^{n})}\left\| \mathbf{v}^{\varepsilon} \right\| _{L^{4}(\mathbb{R}^{n})}\left\|  \nabla v^{\varepsilon} \right\| _{L^{4}(\mathbb{R}^{n})}
  \\[0.3cm]
 &\qquad \leq C\left\| \nabla \mathbf{u}^{\varepsilon} \right\|_{L^{2}(\mathbb{R}^{n})}\left(\left\| \nabla \mathbf{v}^{\varepsilon} \right\| ^{2}_{L^{2}(\mathbb{R}^{n})}+\left\| \Lambda^{\frac{n}{4}}\nabla \mathbf{v}^{\varepsilon} \right\| ^{2}_{L^{2}(\mathbb{R}^{n})}\right)
 \\[0.3cm]
  &\qquad\qquad+\left\| \Delta \mathbf{u}^{\varepsilon}\right \|_{L^{2}(\mathbb{R}^{n})}\left(\left\|  \mathbf{v}^{\varepsilon}  \right\|^{2} _{L^{2}(\mathbb{R}^{n})}+\left\|\Lambda^{\frac{n}{4}} \mathbf{v}^{\varepsilon}  \right\|^{2} _{L^{2}(\mathbb{R}^{n})}\right)^{\frac{1}{2}}
  \left(\left\|  \nabla \mathbf{v}^{\varepsilon} \right\|^{2} _{L^{2}(\mathbb{R}^{n})}+\left\|\Lambda^{\frac{n}{4}} \nabla \mathbf{v}^{\varepsilon} \right\|^{2} _{L^{2}(\mathbb{R}^{n})}\right)^{\frac{1}{2}}
  \\[0.3cm]
 &\qquad \leq C\left\| \nabla \mathbf{u}^{\varepsilon} \right\|_{L^{2}(\mathbb{R}^{n})}\left(\left\| \nabla \mathbf{v}^{\varepsilon} \right\| ^{2}_{L^{2}(\mathbb{R}^{n})}+\left\| \Lambda^{\frac{n}{4}}\nabla \mathbf{v}^{\varepsilon} \right\| ^{2}_{L^{2}(\mathbb{R}^{n})}\right)
 \\[0.3cm]
  &\qquad\qquad+\left\| \Delta \mathbf{u}^{\varepsilon}\right \|_{L^{2}(\mathbb{R}^{n})}\left(\left\|  \mathbf{v}^{\varepsilon}  \right\|^{2} _{L^{2}(\mathbb{R}^{n})}+\left\|\Lambda^{\frac{n}{4}} \mathbf{v}^{\varepsilon}  \right\|^{2} _{L^{2}(\mathbb{R}^{n})}\right)+\left\| \Delta \mathbf{u}^{\varepsilon}\right \|_{L^{2}(\mathbb{R}^{n})}
  \left(\left\|  \nabla \mathbf{v}^{\varepsilon} \right\|^{2} _{L^{2}(\mathbb{R}^{n})}+\left\|\Lambda^{\frac{n}{4}} \nabla \mathbf{v}^{\varepsilon} \right\|^{2} _{L^{2}(\mathbb{R}^{n})}\right)
 \\[0.3cm]
  &\qquad \lesssim \left(\left\| \nabla \mathbf{u}^{\varepsilon} \right\|_{L^{2}(\mathbb{R}^{n})}+\left\| \Delta \mathbf{u}^{\varepsilon}\right \|_{L^{2}(\mathbb{R}^{n})}\right)\left\| \nabla \mathbf{v}^{\varepsilon} \right\| ^{2}_{L^{2}(\mathbb{R}^{n})}
   \\[0.3cm]
  &\qquad\qquad+\left(\left\| \nabla \mathbf{u}^{\varepsilon} \right\|_{L^{2}(\mathbb{R}^{n})}+\left\| \Delta \mathbf{u}^{\varepsilon}\right \|_{L^{2}(\mathbb{R}^{n})}\right)\left\|\Lambda^{\frac{n}{4}} \nabla \mathbf{v}^{\varepsilon} \right\|^{2} _{L^{2}(\mathbb{R}^{n})}.
   \end{array}
\right.
\end{equation}
  By Proposition \ref{p2.4}, the assumptions in (III) of this theorem, once we choose $\left\| \nabla \mathbf{u}^{\varepsilon} \right\|_{L^{2}(\mathbb{R}^{n})}+\left\| \Delta \mathbf{u}^{\varepsilon}\right \|_{L^{2}(\mathbb{R}^{n})}\lesssim \left\|  \mathbf{v}^{\varepsilon} \right\|_{L^{2}(\mathbb{R}^{n})} \lesssim\| \mathbf{v}_{0}^{\varepsilon}   \| _{L^{2}(\mathbb{R}^{n})}\lesssim\| \mathbf{v}_{0}  \| _{L^{2}(\mathbb{R}^{n})}\leq  \varepsilon^{*} \leq \frac{\nu}{4}$, \eqref{3.22}, \eqref{3.26} and \eqref{3.27} yield
     \begin{equation}\label{3.28}
     \frac{d}{dt}\| \nabla  \mathbf{v}^{\varepsilon} \| ^{2}_{L^{2}(\mathbb{R}^{n})}+\frac{\nu}{2}\left\| \Lambda^{\frac{n}{4}}\nabla  \mathbf{v}^{\varepsilon} \right\| ^{2}_{L^{2}(\mathbb{R}^{n})}\lesssim\|  \mathbf{v}^{\varepsilon} \| ^{2}_{L^{2}(\mathbb{R}^{n})}\| \nabla \mathbf{v}^{\varepsilon} \| ^{2}_{L^{2}(\mathbb{R}^{n})}.
     \end{equation}
Using \eqref{3.25} and \eqref{3.28}, Gronwall's lemma yields that for $\displaystyle\frac{n}{4}\leq \beta <1$ with $n=2,3$
\begin{equation}\label{3.29}
    \| \nabla \mathbf{ v}^{\varepsilon} \| ^{2}_{L^{2}(\mathbb{R}^{n})}\leq\| \nabla  \mathbf{v}^{\varepsilon}_{0} \| ^{2}_{L^{2}(\mathbb{R}^{n})}e^{C\|  \mathbf{v}^{\varepsilon}_{0} \| ^{2}_{L^{2}(\mathbb{R}^{n})}}\leq C\varepsilon^{2}.
    \end{equation}
 This gives
  \begin{equation}\label{3.30}
 \left\| \Lambda^{\beta}  \mathbf{u}^{\varepsilon} \right\| ^{2}_{L^{2}(\mathbb{R}^{n})}+\alpha^{2}\left\| \Lambda^{\beta}  \nabla \mathbf{u}^{\varepsilon} \right\| ^{2}_{L^{2}(\mathbb{R}^{n})}  \leq\left\| \nabla \mathbf{u}^{\varepsilon}\right \| ^{2}_{L^{2}(\mathbb{R}^{n})}+\alpha^{2}\left\| \Delta \mathbf{u}^{\varepsilon} \right\| ^{2}_{L^{2}(\mathbb{R}^{n})}  \leq\left\| \nabla \mathbf{v}^{\varepsilon} \right\| ^{2}_{L^{2}(\mathbb{R}^{n})}\leq C\varepsilon^{2}.
 \end{equation}
 It follows from \eqref{2.1} and \eqref{3.30} that
  $$\frac{d}{dt}\left(\| \mathbf{u}^{\varepsilon} \| ^{2}_{L^{2}(\mathbb{R}^{n})}+\alpha^{2}\| \nabla  \mathbf{u}^{\varepsilon} \| ^{2}_{L^{2}(\mathbb{R}^{n})}\right)\geq -C\varepsilon^{2}.$$
  This is the estimate \eqref{3.18}, and thus the proof of (III) is finished.\\
  \indent So far, we finish the proof of Theorem 3.1.\hfill$\Box$\\
 \section{Algebraic Decay }
 \allowdisplaybreaks
Motivated by these works concerning the algebraic decay of the imcompressible Navier-Stokes equations \cite{Kato,Leray}, in this section we shall establish the algebraic decay estimate for the solutions of the Cauchy problem \eqref{1.1}-\eqref{1.2}. From Section 3, we have known that there is no uniform rate of decay for solutions with data exclusively in $D_{\sigma}\left(\Lambda \right)(\mathbb{R}^{n})$ for $\dfrac{n}{4}< \beta < 1$, and in $H_{\sigma}^{1}(\mathbb{R}^{n})$ for $\beta=\dfrac{n}{4}$.
 However, we claim here that there is a uniform rate of decay depending on $D_{\sigma}\left(\Lambda \right)(\mathbb{R}^{n})$ and $L^{1}(\mathbb{R}^{n})$ norms of the initial data for $\dfrac{n}{4}< \beta < 1$, and on $H_{\sigma}^{1}(\mathbb{R}^{n})$ and $L^{1}(\mathbb{R}^{n})$ norms of the initial data for $\beta=\dfrac{n}{4}$. We first in this section establish the decay rate for the filtered velocity $\mathbf{u}$ by applying the Fourier splitting argument introduced in \cite{Kozono-1,Kozono-2} to the natural energy relation \eqref{2.2}. This decay rate is then applied with an inductive argument to achieve deacy rates for the unfiltered velocity $\mathbf{v}$ and all of its derivatives. It should be pointed out that the Fourier splitting method was originally applied to parabolic conservation laws in \cite{Schonbek-3}, and later applied to Navier-Stokes equations in \cite{Schonbek-4}.
\\[0.3cm]
\indent The algebraic decay result is the following.
\begin{thm}\label{t4.1}\rm
 For $\displaystyle\frac{n}{4}\leq \beta <1$, $n=2,3$, let $\mathbf{v}$ be the solution of the Camassa-Holm equations with fractional Laplacian viscosity \eqref{1.1}-\eqref{1.2} constructed in Proposition \ref{p2.4}. Then we have
\\[0.3cm]
(I)\quad If $\mathbf{v}_{0}\in D_{\sigma}\left(\Lambda \right)(\mathbb{R}^{n})\cap L^{1}(\mathbb{R}^{n})$ for $\dfrac{n}{4}< \beta < 1$, and
 $\mathbf{v}_{0}\in H^{1}_{\sigma}(\mathbb{R}^{n})\cap L^{1}(\mathbb{R}^{n})$ for $\beta=\dfrac{n}{4}$ with an additional assumption that there exists an $\varepsilon^{*}=\varepsilon^{*}(\alpha,\nu,n)$ sufficiently small such that $\|\mathbf{v}_{0}\|_{H^{1}_{0}(\mathbb{R}^{n})}\lesssim  \varepsilon^{*} $,~ then
the solution satisfies the "energy" decay rate
$$\displaystyle\int _{\mathbb{R}^{n}}\mathbf{v}\cdot \mathbf{u} dx=\|\mathbf{u}\|_{L^{2} (\mathbb{R}^{n})}^{2}+
\alpha^{2}\|\nabla \mathbf{u}\|_{L^{2} (\mathbb{R}^{n})}^{2}\leq C (1+t)^{-\frac{n}{2\beta}}.$$
\\[0.2cm]
(II)\quad Under the condition of (I), the solution satisfies the decay rate
$$\|\nabla \mathbf{v}\|_{L^{2} (\mathbb{R}^{n})}^{2}
\leq C (1+t)^{-\frac{1}{ \beta}-\frac{n}{2\beta}}.$$
 (III)\quad Under the condition of (I),~ then\\
\\
\indent \qquad (III-1) \quad $|\mathcal{F}(\mathbf{v})|\leq C$, \quad  $|\mathcal{F}(\mathbf{u})|\leq C$,
\\[0.20cm]
\indent \qquad (III-2) \quad $\|\mathbf{v}\|_{L^{2} (\mathbb{R}^{n})}^{2}\leq C (1+t)^{-\frac{n}{2\beta}} $.\\
\\
(IV)\quad Let $\|\nabla^{m}w_{0}\|_{L^{2} (\mathbb{R}^{n})}<\infty$. Given an energy inequality of the form
\begin{equation}\label{4.1}
\frac{1}{2}\frac{d}{dt}\left\|\nabla^{m}w\right\|_{L^{2} (\mathbb{R}^{n})}^{2}+\nu\left\|\Lambda^{\beta}\nabla^{m}w\right\|_{L^{2} (\mathbb{R}^{n})}^{2}\leq C (1+t)^{\gamma},
\end{equation}
 and the bound $|\hat{w}(\xi,t)|\leq C (1+t)^{ \eta}$ which holds for $|\xi|^{2\beta}<\frac{b}{\nu(1+t)}$, we then achieve
 \begin{equation}\label{4.2}
 \left\|\nabla^{m}w\right\|_{L^{2} (\mathbb{R}^{n})}^{2}\leq C\left[(1+t)^{-\frac{m}{\beta} -\frac{n}{2\beta}+2\eta}+(1+t)^{\gamma+1}\right].
 \end{equation}
 \\
(V)\quad Let $\mathbf{v}_{0}\in H^{1}_{\sigma}(\mathbb{R}^{n})\cap L^{1}(\mathbb{R}^{n})$. For $P\geq 1$, if for all $p<P$ and $m=0,1$,
  $$\left\|\partial^{p}_{t}\nabla^{m}\mathbf{v}\right\|_{L^{2} (\mathbb{R}^{n})}^{2}\leq C (1+t)^{-2p-\frac{m}{\beta} -\frac{n}{2\beta}},$$
  then for $|\xi|^{2\beta}\leq \frac{b}{\nu(1+t)}$, there holds
   $$\left|\partial^{P}_{t} \hat{\mathbf{v}}(\xi)\right|\leq C (1+t)^{-P}.$$
   \\
(VI)\quad If $\mathbf{v}_{0}\in D_{\sigma}\left(\Lambda^{K} \right)(\mathbb{R}^{n})\cap L^{1}(\mathbb{R}^{n})$ for $\dfrac{n}{4}< \beta < 1$, and
 $\mathbf{v}_{0}\in H^{K}_{\sigma}(\mathbb{R}^{n})\cap L^{1}(\mathbb{R}^{n})$ for $\beta=\dfrac{n}{4}$ with an additional assumption that there exists an $\varepsilon^{**}=\varepsilon^{**}(\alpha,\nu,n)$ sufficiently small such that $\|\mathbf{v}_{0}\|_{H^{K}_{0}(\mathbb{R}^{n})}\lesssim  \varepsilon^{**} $,~   then\\
   \\
  \indent (VI-1)\quad
  For all $m\leq K$, the solution satisfies the following decay
 $$\left\|\nabla^{m}\mathbf{v}\right\|_{L^{2} (\mathbb{R}^{n})}^{2}\leq C (1+t)^{-\frac{m}{\beta} -\frac{n}{2\beta}}.$$
\indent (VI-2)\quad For all $m+2p\beta\leq K$, the solution satisfies the decay estimate
 $$\left\|\partial^{p}_{t}\nabla^{m}\mathbf{v}\right\|_{L^{2} (\mathbb{R}^{n})}^{2}\leq C (1+t)^{-2p-\frac{m}{\beta} -\frac{n}{2\beta}}.$$
 Here, $m$, $p$ and $P$ are all non-negative integers in (IV), (V) and (VI), the constant $C$ in (I)-(VI) depends only on the initial data, the dimension of space, and the constants in \eqref{1.1}, which may be different on different lines.
\end{thm}
\indent We shall apply the Fourier splitting method and the bootstrap argument to show Theorem \ref{t4.1}. Before going further, we first establish an estimate on
$\left\|  \hat{\mathbf{v}} \right\|  _{L^{\infty}(\mathbb{R}^{n})}$.\\
\begin{lem} \label{l4.2}\rm
 Let $\mathbf{v}$ be the solution of \eqref{1.1}-\eqref{1.2} constructed in Proposition \ref{p2.4} corresponding to $\mathbf{v}_{0}\in \left(L^{2}_{\sigma}\cap L^{1}\right)(\mathbb{R}^{n})$. Then
\begin{equation}\label{4.3}
\left|\mathcal{F}(\mathbf{v})\right|\leq C\left[1+\left(\int^{t}_{0}\|   \mathbf{u}(s)  \|^{2} _{L^{2}(\mathbb{R}^{n})}ds\right)^{\frac{1}{2}}\left(\int^{t}_{0}\| \nabla \mathbf{v}(s)  \|^{2} _{L^{2}(\mathbb{R}^{n})}ds\right)^{\frac{1}{2}}\right].
\end{equation}
Here, $C$ depends only on the initial data, the dimension of space and the constants in (1.1), but not on $\alpha$.\end{lem}
{\bf Proof.} Note that
\begin{equation}\label{4.4}
\sum\limits_{i=1}^{n}\nabla (u_{i}v_{i})=\sum\limits_{i=1}^{n} u_{i}\nabla v_{i} +\sum\limits_{i=1}^{n}v_{i}\nabla u_{i}=\mathbf{u}\cdot \nabla \mathbf{v}^{T}+\mathbf{v}\cdot \nabla \mathbf{u}^{T},
\end{equation}
taking the Fourier transform with respect to $x$ for the first equation in (1.1) yields
$$\hat{\mathbf{v}}_{t}=-\mathcal{F}\left(\mathbf{u}\cdot \nabla \mathbf{v}-\mathbf{u}\cdot \nabla \mathbf{v}^{T}\right)-\mathcal{F}\left[\nabla\left(\sum\limits_{i=1}^{n} u_{i}v_{i}+p\right)\right]-\nu|\xi|^{2\beta}\hat{\mathbf{v}}.$$
A straightforward computation shows that
\begin{equation}\label{4.5}
\mathcal{F}(\mathbf{v})=e^{-\nu t |\xi|^{2\beta} }\mathcal{F}(\mathbf{v}_{0})+\int^{t}_{0}e^{-\nu (t-s)|\xi|^{2\beta}} \psi(\xi,s)ds,
\end{equation}
where
\begin{equation}\label{4.6}
\psi(\xi,t)=-\xi \cdot\mathcal{F}\left(p+\sum\limits_{i=1}^{n} u_{i}v_{i} \right)-\mathcal{F}\left(\mathbf{u}\cdot \nabla \mathbf{v}-\mathbf{u}\cdot \nabla \mathbf{v}^{T} \right).
\end{equation}
We first deal with the term $\psi(\xi,t)$.\\
 \indent Thanks to $\left\| \mathcal{F}(\phi) \right\|  _{L^{\infty}(\mathbb{R}^{n})}\leq \| \phi\|  _{L^{1}(\mathbb{R}^{n})}$ and Young's inequality, one deduces that
\begin{equation}\label{4.7}
 \left|\mathcal{F}(\mathbf{u}\cdot \nabla \mathbf{v})\right| \leq \|\mathbf{u}\cdot \nabla \mathbf{v} \|  _{L^{1}(\mathbb{R}^{n})}
 \leq  C\|\mathbf{u} \|  _{L^{2}(\mathbb{R}^{n})}
\|\nabla \mathbf{v}\|  _{L^{2}(\mathbb{R}^{n})}.
\end{equation}
In the same manner, one achieves
\begin{equation}\label{4.8}
\left|\mathcal{F}\left(\mathbf{u}\cdot \nabla \mathbf{v}^{T}\right)\right|\leq\|\mathbf{u} \|  _{L^{2}(\mathbb{R}^{n})}
\|\nabla \mathbf{v}\|  _{L^{2}(\mathbb{R}^{n})}.
\end{equation}
On the other hand, taking the divergence for the first equation in \eqref{1.1} leads to
\begin{equation}\label{4.9}
-\Delta\left(p+\sum\limits_{i=1}^{n} u_{i}v_{i} \right)=\mbox{div}\left(\mathbf{u}\cdot \nabla \mathbf{v}-\mathbf{u}\cdot \nabla \mathbf{v}^{T} \right).
\end{equation}
Combining \eqref{4.6} with \eqref{4.7}, \eqref{4.8} and \eqref{4.9} yields that
$$\left|\psi(\xi,t)\right|\leq C\|\mathbf{u} \|  _{L^{2}(\mathbb{R}^{n})}\|\nabla \mathbf{v}   \|  _{L^{2}(\mathbb{R}^{n})}.$$
Taking the supremum over $\xi$ for \eqref{4.4} and applying Cauchy-Schwarz inequality, one obtains
$$\left|\mathcal{F}(\mathbf{v})\right|\leq |\mathcal{F}(\mathbf{v}_{0})|+C \left(\int^{t}_{0}\|   \mathbf{u}(s)  \|^{2} _{L^{2}(\mathbb{R}^{n})}ds\right)^{\frac{1}{2}}\left(\int^{t}_{0}\| \nabla \mathbf{v}(s)  \|^{2} _{L^{2}(\mathbb{R}^{n})}ds\right)^{\frac{1}{2}} . $$
In view of  $\left\| \mathcal{F}(\mathbf{v}_{0}) \right\|  _{L^{\infty}(\mathbb{R}^{n})}\leq \| \mathbf{v}_{0} \|  _{L^{1}(\mathbb{R}^{n})}$, the above inequality deduces the desired estimate \eqref{4.3}. \hfill$\Box$
\\[0.3cm]

In the following, we start the proof of Theorem \ref{t4.1}.
 \\[0.3cm]
  {\bf Proof of Theorem \ref{t4.1}}.
 \\[0.3cm]
 {\bf We first show (I)}.
 \\[0.3cm]
\indent Note that the assumption of (I), by Proposition \ref{p2.4}, one obtains
$$\int^{t}_{0}\left\|  \nabla\Lambda^{\beta} \mathbf{v}(s) \right \|^{2} _{L^{2}(\mathbb{R}^{n})}ds\leq C,$$
and
$$\int^{t}_{0}\left\|  \Lambda^{\beta} \mathbf{v}(s)  \right\|^{2} _{L^{2}(\mathbb{R}^{n})}ds\leq C,$$
which imply by interpolation inequality that
$$\int^{t}_{0}\|  \nabla \mathbf{ v}(s)  \|^{2} _{L^{2}(\mathbb{R}^{n})}ds\leq C,$$
where $C$ depends only on $n,\beta,\alpha,\nu$ and
$\left\{
\begin{array}{ll}
\|\mathbf{v}_{0} \|_{D \left(\Lambda \right)(\mathbb{R}^{n})}~~&{for}~\dfrac{n}{4}< \beta < 1,
\\[0.2cm]
\|\mathbf{v}_{0} \|_{H^{1}_{0}(\mathbb{R}^{n})}~~&{for}~\beta=\dfrac{n}{4}.
\end{array}\right. $.
Lemma \ref{l4.2} then gives rise to
\begin{equation}\label{4.10}
\left.
\begin{array}{ll}
\vspace{5pt}
\displaystyle|\hat{\mathbf{v}}|\leq  C\left[1+\left(\int^{t}_{0}\|   \mathbf{u}(s)  \|^{2} _{L^{2}(\mathbb{R}^{n})}ds\right)^{\frac{1}{2}}\left(\int^{t}_{0}\|  \nabla  \mathbf{v}(s)  \|^{2} _{L^{2}(\mathbb{R}^{n})}ds\right)^{\frac{1}{2}}\right]\\
\displaystyle\quad \leq  C\left[1+ \left(\int^{t}_{0}\|   \mathbf{u}(s)  \|^{2} _{L^{2}(\mathbb{R}^{n})}ds \right)^{\frac{1}{2}} \right].
 \end{array}
\right.
\end{equation}
Thanks to the Plancherel's theorem, the energy equality (2.1) is equivalent to
$$\frac{d}{dt}\int _{\mathbb{R}^{n}}\left(1+\alpha^{2}|\xi|^{2}\right)\hat{\mathbf{u}}^{2}d\xi
+2\nu\int _{\mathbb{R}^{n}}|\xi|^{2\beta}\left(1+\alpha^{2}|\xi|^{2}\right)\hat{\mathbf{u}}^{2}d\xi=0.$$
Let $B(\rho)$ be the ball of radius $\rho$ with $\displaystyle\rho^{2\beta}=\frac{\frac{n}{2\beta}+1}{2\nu (1+t)}$. Put \begin{equation}\label{4.11}
E^{2}=\hat{\mathbf{u}}\cdot\hat{\mathbf{v}}=\left(1+\alpha^{2}|\xi|^{2}\right)\hat{\mathbf{u}}^{2}.
\end{equation}  Then
$$\frac{d}{dt}\int _{\mathbb{R}^{n}}E^{2}d\xi
+2\nu\rho^{2\beta}\int _{B^{c}(\rho)}E^{2}d\xi\leq0,$$
or
\begin{equation}\label{4.12}\frac{d}{dt}\int _{\mathbb{R}^{n}}E^{2}d\xi
+2\nu\rho^{2\beta}\int _{\mathbb{R}^{n}}E^{2}d\xi\leq2\nu\rho^{2\beta}\int _{B(\rho)}E^{2}d\xi.\end{equation}
The equation $\mathbf{u}-\alpha^{2}\Delta \mathbf{u}=\mathbf{v}$ implies that $\displaystyle\hat{\mathbf{u}}=
\frac{\hat{\mathbf{v}}}{1+\alpha^{2}|\xi|^{2}}$. This together with \eqref{4.10} and \eqref{4.11} yields that
$$\left\| E^{2}\right \| _{L^{\infty}(\mathbb{R}^{n})}\leq \frac{\left\| \hat{\mathbf{v}}^{2} \right\| _{L^{\infty}(\mathbb{R}^{n})}}{1+\alpha^{2}|\xi|^{2}}\leq C\left[1+\int^{t}_{0}\| u(s)  \|^{2} _{L^{2}(\mathbb{R}^{n})}ds\right].$$
With this, \eqref{4.12} then leads to
$$\frac{d}{dt}\int _{\mathbb{R}^{n}}E^{2}d\xi
+2\nu\rho^{2\beta}\int _{\mathbb{R}^{n}}E^{2}d\xi\leq 2\nu\rho^{2\beta+n}\left[1+\int^{t}_{0}\|  \mathbf{u}(s)  \|^{2} _{L^{2}(\mathbb{R}^{n})}ds\right],$$
which yields a differential inequality by using the integrating factor $\displaystyle f=(1+t)^{\frac{n}{2\beta}+1}$:
$$\frac{d}{dt}\left((1+t)^{\frac{n}{2\beta}+1}\int _{\mathbb{R}^{n}}E^{2}d\xi\right)\leq 2\nu \rho^{2\beta+n}(1+t)^{\frac{n}{2\beta}+1}
 \left[1+\int^{t}_{0}\| \mathbf{u}(s)  \|^{2} _{L^{2}(\mathbb{R}^{n})}ds\right].$$
Integrating this differential inequality in time $t$ from $0$ to $r$ gives rise to
\begin{equation}\label{4.13}
\left.
\begin{array}{ll}
\vspace{5pt}
\displaystyle(1+r)^{\frac{n}{2\beta}+1}\int _{\mathbb{R}^{n}}E^{2}(\xi,r)d\xi\\
\qquad\qquad\displaystyle\leq \int _{\mathbb{R}^{n}}E^{2}(\xi,0)d\xi+\dfrac{\left(\frac{n}{2\beta}+1\right)^{\frac{n}{2\beta}+1}}
{(2\nu)^{\frac{n}{2\beta}}} \int^{r}_{0}\left(1+\int^{t}_{0}\|   \mathbf{u}(s)  \|^{2} _{L^{2}(\mathbb{R}^{n})}ds\right)dt.
\end{array}
\right.
\end{equation}
By the Tonelli theorem, a simple calculation shows that
$$
\left.
\begin{array}{ll}
\vspace{5pt}
&\displaystyle\int^{r}_{0}(r-s) \|   \mathbf{ u}(s)  \|^{2} _{L^{2}(\mathbb{R}^{n})}ds\\
\vspace{5pt}
 &\qquad=r\displaystyle\int^{r}_{0}  \|    \mathbf{u}(s)  \|^{2} _{L^{2}(\mathbb{R}^{n})}ds
-s\displaystyle\int^{s}_{0}  \|    \mathbf{u}(t)  \|^{2} _{L^{2}(\mathbb{R}^{n})}dt|^{r}_{0}
+\displaystyle\int^{r}_{0}  \int^{s}_{0}  \|    \mathbf{u}(t)  \|^{2} _{L^{2}(\mathbb{R}^{n})}dtds\\
\vspace{5pt}
 &\qquad  \geq\displaystyle\int^{r}_{0}  \displaystyle\int^{t}_{0}  \|    \mathbf{u}(s)  \|^{2} _{L^{2}(\mathbb{R}^{n})}dsdt.
 \end{array}
 \right.
$$
Furthermore, it is a simple exercise to obtain the following estimate
\begin{equation*}
\left.
\begin{array}{ll}
\vspace{5pt}
\displaystyle\int^{r}_{0}\left(1+\int^{t}_{0}\|   \mathbf{u}(s)  \|^{2} _{L^{2}(\mathbb{R}^{n})}ds\right)dt
\leq r+\displaystyle\int^{r}_{0}  \displaystyle\int^{t}_{0}  \|    \mathbf{u}(s)  \|^{2} _{L^{2}(\mathbb{R}^{n})}dsdt\\
\qquad\qquad\qquad\qquad\qquad\qquad\displaystyle\leq r+\int^{r}_{0}(r-s) \|    \mathbf{u}(s)  \|^{2} _{L^{2}(\mathbb{R}^{n})}ds.
\end{array}
\right.
\end{equation*}
Due to $\displaystyle\frac{n}{4}\leq \beta<1$, ~ $\hat{\mathbf{u}}^{2}=\dfrac{E^{2}}{1+\alpha^{2}|\xi|^{2}}\leq
E^{2}$ by \eqref{4.11} and $\|  \mathbf{u}(s)  \|^{2} _{L^{2}(\mathbb{R}^{n})}\leq \displaystyle\int _{\mathbb{R}^{n}}E^{2}d\xi$, if follows from \eqref{4.13} that
\begin{equation}\label{4.14}
(1+r)^{\frac{n}{2\beta}+1}\int _{\mathbb{R}^{n}}E^{2}(\xi,r)d\xi\leq
 C(1+r)+C\int^{r}_{0}(r-s)\int _{\mathbb{R}^{n}}E^{2}(\xi,s)d\xi ds.
 \end{equation}
Let $\phi(r)=(1+r)^{\frac{n}{2\beta}+1}\displaystyle\int _{\mathbb{R}^{n}}E^{2}(\xi,r)d\xi$. (4.14) has the following equivalent form:
$$\phi\leq C(1+r)+C\int^{r}_{0}\phi(s)(r-s)(1+s)^{-\frac{n}{2\beta}-1}ds.$$
The Gronwall inequality implies that
\begin{equation}\label{4.15}
(1+r)^{\frac{n}{2\beta}+1}\int _{\mathbb{R}^{n}}E^{2}(\xi,r)d\xi\leq
 C(1+r)exp\left(C\int^{r}_{0}(r-s)(1+s)^{-\frac{n}{2\beta}-1}ds\right).
 \end{equation}
\indent Thanks to the fact that $\displaystyle\frac{n}{2}<\frac{n}{2\beta}\leq 2$ for $\frac{n}{4}\leq\beta< 1$ and $n=2,3$, the integral $\displaystyle\int^{r}_{0}(r-s)(1+s)^{-\frac{n}{2\beta}-1}ds$ is bounded independent of $r$. Applying the Plancherel's theorem finishes the proof of (I).\\
\\
{\bf We next prove (II)}.
\\[0.3cm]
\indent Recalling that $\Delta \mathbf{v}$ is divergence free, thanks to the identity \eqref{4.4} and H\"{o}lder's inequality, multiplying the first equation in \eqref{1.1} by $\Delta \mathbf{v}$,  we obtain
\begin{equation}\label{4.16}
\left.
\begin{array}{ll}
\vspace{5pt}
 \displaystyle \frac{1}{2}\frac{d}{dt}\left\| \nabla  \mathbf{v}  \right\| ^{2}_{L^{2}(\mathbb{R}^{n})}+\nu\left\| \Lambda^{\beta}\nabla \mathbf{ v } \right\| ^{2}_{L^{2}(\mathbb{R}^{n})} &\lesssim \left\langle \mathbf{u} \cdot\nabla \mathbf{v} , \Delta \mathbf{v} \right\rangle+C\left\langle \mathbf{u} \cdot\nabla \mathbf{v}^{T} , \Delta \mathbf{v} \right\rangle\\
  \vspace{5pt}
  & \displaystyle \lesssim\left\| \Lambda^{1-\beta}(\mathbf{u} \cdot\nabla \mathbf{v})  \right\| _{L^{2}(\mathbb{R}^{n})}
 \left\| \Lambda^{\beta}   \nabla \mathbf{v}  \right \| _{L^{2}(\mathbb{R}^{n})}.
 \end{array}
\right.
\end{equation}
 There are two cases to consider for estimating the term $\left\| \Lambda^{1-\beta}(\mathbf{u} \cdot\nabla \mathbf{v})  \right\| _{L^{2}(\mathbb{R}^{n})} $. \\
 \\
  \indent Case {\bf  (I)}\quad $\displaystyle\frac{n}{4}< \beta< 1$ for $n=2,3$;
  \\[0.3cm]
  \indent Case {\bf  (II)}\quad $  \displaystyle\beta=\frac{n}{4}$ for $n=2,3$.\\
  \\
 $\heartsuit$ We first deal with {\bf Case (I)}\quad $\displaystyle\frac{n}{4}<\beta<1$ for $n=2,3$. The following auxiliary computations will be needed for {\bf Case (I)}.
 \\[0.2cm]
 \begin{equation}\label{4.17}
\left\{
\begin{array}{ll}
\displaystyle\frac{n}{\beta}=\dfrac{2n}{n-2(n-2\beta)/2},~~\dfrac{1}{2}=\dfrac{n-2}{2}<\dfrac{n-2\beta}{2}
  <\dfrac{3}{4}~~
  &\mbox{for}~~n=3,
  \\[0.4cm]
    \displaystyle 0=\dfrac{n-2}{2}<\dfrac{n-2\beta}{2}< \dfrac{1}{2}~~
  &\mbox{for}~~n=2,
  \\[0.3cm]
    \displaystyle B=\dfrac{n-2\beta}{2}+1-\beta=\dfrac{n}{2}+1-2\beta,\quad \dfrac{n}{2}-1<B< 1~~  &\mbox{for}~~n=2,3.
 \end{array}
\right.
\end{equation}\\
 A straightforward computation shows that
 \begin{equation}\label{4.18}
\left.
\begin{array}{ll}
\vspace{5pt}
&\left\| \Lambda^{1-\beta} (\mathbf{u}\cdot \nabla \mathbf{v})\right\|_{L^{2} (\mathbb{R}^{n})}\\
 \vspace{5pt}
 &\qquad\leq\left\| \Lambda^{1-\beta} (\mathbf{u}\cdot \nabla \mathbf{v})-\Lambda^{1-\beta}  \mathbf{u}\cdot \nabla \mathbf{v}-  \mathbf{u}\Lambda^{1-\beta} \nabla \mathbf{v} \right\|_{L^{2} (\mathbb{R}^{n})}\\
  \vspace{5pt}
  &\qquad\qquad +\left\|  \Lambda^{1-\beta}  \mathbf{u}\cdot \nabla \mathbf{v } \right\|_{L^{2} (\mathbb{R}^{n})}+\left\|  \mathbf{u}\Lambda^{1-\beta} \nabla \mathbf{v} \right\|_{L^{2} (\mathbb{R}^{n})}.
  \end{array}
\right.
\end{equation}
In view of Lemma \ref{l2.9}, Lemma \ref{l2.10}, Lemma \ref{l2.11} and \eqref{4.17}, for $\displaystyle\frac{1}{2}=\frac{1}{n/\beta}+\frac{1}{2n/(n-2\beta)}$ and $0<1-\beta<\beta<1$, the first term on the right hand side of \eqref{4.18} can be bounded as follows:
\begin{equation}\label{4.19}
\left.
\begin{array}{ll}
\vspace{5pt}
&\left\| \Lambda^{1-\beta} (\mathbf{u}\cdot\nabla \mathbf{v})-\mathbf{u} \Lambda^{1-\beta} \nabla \mathbf{v} - \Lambda^{1-\beta}  \mathbf{u} \nabla \mathbf{v}  \right\|_{L^{2} (\mathbb{R}^{n})}\\
\vspace{5pt}
&\qquad \leq C\left\| \Lambda^{1-\beta}  \mathbf{u}  \right\|_{L^{\frac{n}{\beta}} (\mathbb{R}^{n})}\left\|\nabla \mathbf{v} \right\|_{L^{\frac{2n}{n-2\beta}} (\mathbb{R}^{n})}\\
\vspace{5pt}
&\qquad \leq C\left\| \Lambda^{\frac{n}{2}+1-2\beta}\mathbf{u}   \right\|_{L^{2} (\mathbb{R}^{n})} \left\| \Lambda^{\beta}\nabla \mathbf{v}  \right\|_{L^{2} (\mathbb{R}^{n})}\\
\vspace{5pt}
&\qquad \leq C\left\| \mathbf{u }  \right\|^{\frac{1}{2}}_{L^{2} (\mathbb{R}^{n})}\left\| \nabla \mathbf{u }  \right\|^{\frac{1}{2}}_{L^{2} (\mathbb{R}^{n})} \left\| \Lambda^{\beta} \nabla \mathbf{v}   \right\|_{L^{2} (\mathbb{R}^{n})}.
\end{array}
\right.
\end{equation}
Here, we have used interpolation inequality in the last line. Due to Lemma \ref{l2.11}, a similar estimate to \eqref{4.19} holds for the second term on the right hand side of \eqref{4.18}
 \begin{equation}\label{4.20}
\left.
\begin{array}{ll}
\vspace{5pt}
&\left\|  \Lambda^{1-\beta}  \mathbf{u} \nabla \mathbf{v} \right\|_{L^{2} (\mathbb{R}^{n})}\leq C\left\|  \Lambda^{1-\beta}  \mathbf{u} \right\|_{L^{\frac{n}{\beta}} (\mathbb{R}^{n})} \left\|  \nabla \mathbf{v}\right\|_{L^{\frac{2n}{n-2\beta}} (\mathbb{R}^{n})}\\
\vspace{5pt}
&\qquad\qquad\qquad\quad \leq C\left\| \mathbf{u }  \right\|^{\frac{1}{2}}_{L^{2} (\mathbb{R}^{n})}\left\| \nabla \mathbf{u }  \right\|^{\frac{1}{2}}_{L^{2} (\mathbb{R}^{n})}  \left\| \Lambda^{\beta} \nabla \mathbf{v}  \right\|_{L^{2} (\mathbb{R}^{n})}.
\end{array}
\right.
\end{equation}
In a same manner, recall \eqref{4.17} again, we deduce the estimate for the third term on the right hand side of \eqref{4.18}
\begin{equation}\label{4.21}
\left.
\begin{array}{ll}
\vspace{5pt}
&\left\|  \mathbf{u} \Lambda^{1-\beta}  \nabla \mathbf{v} \right\|_{L^{2} (\mathbb{R}^{n})}\\
\vspace{5pt}
&\qquad \leq C\left\|   \mathbf{u} \right\|_{L^{\frac{n}{2\beta-1}} (\mathbb{R}^{n})} \left\|  \Lambda^{1-\beta}\nabla \mathbf{v}\right\|_{L^{\frac{2n}{n-2(2\beta-1})} (\mathbb{R}^{n})}\\
 \vspace{5pt}
 &\qquad \leq C\left\|    \mathbf{u} \right\|_{L^{\frac{2n}{n-2A}} (\mathbb{R}^{n})} \left\| \Lambda^{\beta}\nabla \mathbf{v}\right\|_{L^{2} (\mathbb{R}^{n})}\\
 \vspace{5pt}
 &\qquad \leq C\left\|  \Lambda^{A}  \mathbf{u} \right\|_{L^{2} (\mathbb{R}^{n})} \left\| \Lambda^{\beta}\nabla \mathbf{v}\right\|_{L^{2} (\mathbb{R}^{n})}\\
 \vspace{5pt}
 &\qquad \leq C\left\| \mathbf{u }  \right\|^{\frac{1}{2}}_{L^{2} (\mathbb{R}^{n})}\left\| \nabla \mathbf{u }  \right\|^{\frac{1}{2}}_{L^{2} (\mathbb{R}^{n})}  \left\|  \Lambda^{\beta} \nabla \mathbf{v}  \right\|_{L^{2} (\mathbb{R}^{n})},
 \end{array}
\right.
\end{equation}
where $A=\frac{n}{2}+1-2\beta$ is given by Lemma \ref{l2.11}.
Note that (I) of this theorem, combining \eqref{4.16} with \eqref{4.17},\eqref{4.18}, \eqref{4.19}, \eqref{4.20} and \eqref{4.21} yields that
\begin{equation}\label{4.22}
\frac{d}{dt} \left\|   \nabla  \mathbf{v} \right\| ^{2}_{L^{2}(\mathbb{R}^{n})}+2\nu\left\|   \Lambda^{\beta}\nabla  \mathbf{v}  \right\| ^{2}_{L^{2}(\mathbb{R}^{n})}\leq C(1+t)^{-\frac{n}{4\beta}}\left\|  \Lambda^{\beta}\nabla  \mathbf{v} \right\| ^{2}_{L^{2}(\mathbb{R}^{n})}.
\end{equation}
 $\heartsuit$ we next consider Case {\bf (II)}~$ \beta= \dfrac{n}{4}$ for $n=2,3$.
 \\[0.3cm]
   \indent In this case, $\left\| \Lambda^{1-\frac{n}{4}} (\mathbf{u}\cdot \nabla \mathbf{v})\right\|_{L^{2} (\mathbb{R}^{n})}$ can be bounded as follows:
\begin{equation}\label{4.23}
\left.
\begin{array}{ll}
\vspace{5pt}
&\left\| \Lambda^{1-\frac{n}{4}} (\mathbf{u}\cdot \nabla \mathbf{v})\right\|_{L^{2} (\mathbb{R}^{n})}\\
\vspace{5pt}
&\qquad\leq\left\|\Lambda^{1-\frac{n}{4}} (\mathbf{u}\cdot \nabla \mathbf{v})-\Lambda^{1-\frac{n}{4}}  \mathbf{u}\cdot \nabla \mathbf{v}-  \mathbf{u}\Lambda^{1-\frac{n}{4}} \nabla \mathbf{v} \right\|_{L^{2} (\mathbb{R}^{n})}\\
 \vspace{5pt}
 &\qquad\quad +\left\| \Lambda^{1-\frac{n}{4}} \mathbf{ u}\cdot \nabla \mathbf{v} \right\|_{L^{2} (\mathbb{R}^{n})}+\left\|  \mathbf{u}\Lambda^{1-\frac{n}{4}} \nabla \mathbf{v} \right\|_{L^{2} (\mathbb{R}^{n})}.
 \end{array}
\right.
\end{equation}
We first bound the first term on the right hand side of estimate \eqref{4.23}. By  Lemma \ref{l2.9}, one attains  that for $\displaystyle 0<\beta_{1}<1-\frac{n}{4}$,
\begin{equation}\label{4.24}
\left.
\begin{array}{ll}
\vspace{5pt}
&\left\| \Lambda^{1-\frac{n}{4}} (\mathbf{u}\cdot\nabla \mathbf{v})-\mathbf{u} \Lambda^{1-\frac{n}{4}} \nabla \mathbf{v} - \Lambda^{1-\frac{n}{4}}  \mathbf{u} \cdot\nabla \mathbf{v }\right\|_{L^{2} (\mathbb{R}^{n})}\\
 \vspace{5pt}
 &\qquad \leq C\left\|\Lambda^{1-\frac{n}{4}-\beta_{1}}  \mathbf{u} \right\|_{L^{\frac{2n}{n-2(n/4+\beta_{1})}} (\mathbb{R}^{n})} \left\| \Lambda^{\beta_{1}}\nabla \mathbf{v}\right\|_{L^{\frac{2n}{n-2(n/4-\beta_{1})}} (\mathbb{R}^{n})}\\
 \vspace{5pt}
 &\qquad \leq C\left\| \nabla \mathbf{u}  \right\|_{L^{2} (\mathbb{R}^{n})} \left\| \Lambda^{\frac{n}{4}}\nabla \mathbf{v}  \right\|_{L^{2} (\mathbb{R}^{n})}
 \\
 \vspace{5pt}
 &\qquad \lesssim\left\|  \mathbf{v}_{0}  \right\|_{L^{2} (\mathbb{R}^{n})} \left\| \Lambda^{\frac{n}{4}}\nabla \mathbf{v}  \right\|_{L^{2} (\mathbb{R}^{n})}
 ,
 \end{array}
\right.
\end{equation}
where we have used the fact that $\displaystyle\beta_{1}\in \left(0,\frac{1}{2}\right), ~~1-\frac{n}{4}-\beta_{1}\in \left(0,\frac{1}{2}\right),~~\frac{1}{2}=\frac{1}{p_{1}}+\frac{1}{p_{2}}$ with $p_{1},p_{2}\in (1,\infty)$, $\displaystyle p_{1}=\frac{2n}{n-2(n/4+\beta_{1})}$, $\displaystyle p_{2}=\frac{2n}{n-2(n/4-\beta_{1})}$. Thanks to Lemma \ref{l4.3}, Agmon's inequality and the interpolation inequality, note that $0<1-\dfrac{n}{4}<\dfrac{n}{4}$, \eqref{2.4} and the assumption of (II) for $\beta=\dfrac{n}{4}$, the second and the third terms on the right hand side of \eqref{4.23} enjoy the similar estimates to \eqref{4.24}
 \begin{equation}\label{4.25}
\left.
\begin{array}{ll}
\vspace{5pt}
\left\| \Lambda^{1-\frac{n}{4}}  \mathbf{u}\cdot \nabla \mathbf{v} \right\|_{L^{2} (\mathbb{R}^{n})}+\left\|   \mathbf{u}\Lambda^{1-\frac{n}{4}} \nabla \mathbf{v} \right\|_{L^{2} (\mathbb{R}^{n})}\\
 \vspace{5pt}
 \qquad\qquad\leq C\left\|  \Lambda^{1-\frac{n}{4}}  \mathbf{u} \right\|_{L^{4} (\mathbb{R}^{n})}\left\| \nabla \mathbf{v} \right\|_{L^{4} (\mathbb{R}^{n})}+C\left\|\mathbf{ u} \right\|_{L^{\infty} (\mathbb{R}^{n})}\left\| \Lambda^{1-\frac{n}{4}} \nabla \mathbf{v} \right\|_{L^{2} (\mathbb{R}^{n})}\\
 \vspace{5pt}
 \qquad\qquad\leq C\left( \left\|\Lambda^{1-\frac{n}{4}}\mathbf{u}\right\| ^{2}_{L^{2} (\mathbb{R}^{n})}+\left\| \Lambda^{1-\frac{n}{4}+\frac{n}{4}}\mathbf{u}\right\|^{2} _{L^{2} (\mathbb{R}^{n})}\right)^{\frac{1}{2}}\\
  \vspace{5pt}
  \qquad\qquad\qquad\qquad\cdot\left( \left\| \nabla \mathbf{v}\right\|^{2} _{L^{2} (\mathbb{R}^{n})}+\left\|\Lambda^{ \frac{n}{4} }\nabla \mathbf{v}\right\|^{2} _{L^{2} (\mathbb{R}^{n})}\right)^{\frac{1}{2}}\\
  \vspace{5pt}
  \qquad\qquad\qquad+
 \left\| \mathbf{ u}  \right\|^{\frac{1}{2}}_{H^{1} (\mathbb{R}^{n})}\left\| \mathbf{ u}  \right\|^{\frac{1}{2}}_{H^{2} (\mathbb{R}^{n})}  \left\|   \nabla \mathbf{v}   \right\|^{\frac{1}{2}} _{L^{2} (\mathbb{R}^{n})}\left\| \Lambda^{\frac{n}{4}} \nabla \mathbf{v}   \right\|^{\frac{1}{2}} _{L^{2} (\mathbb{R}^{n})}\\
  \vspace{5pt}
  \qquad\qquad\lesssim \left\|  \mathbf{v}\right\| _{H^{1}_{0} (\mathbb{R}^{n})}\left\| \Lambda^{ \frac{n}{4} }\nabla \mathbf{v}\right\|_{L^{2} (\mathbb{R}^{n})},
  \end{array}
\right.
\end{equation}
where we have used the fact that $\left\|   \mathbf{v}\right\| _{L^{2} (\mathbb{R}^{n})}\leq C\left(\left\| \mathbf{v}_{0}\right\| _{L^{2} (\mathbb{R}^{n})}\right)$ and   $\left\|  \nabla \mathbf{v}\right\| _{L^{2} (\mathbb{R}^{n})}\leq C\left(\left\| \mathbf{v}_{0}\right\| _{H^{1} _{0} (\mathbb{R}^{n})}\right)$ for $\left\| \mathbf{v}_{0}\right\| _{H^{1} _{0} (\mathbb{R}^{n})}$ small sufficiently from Proposition \ref{p2.4} in the last inequality of \eqref{4.25}. Combining \eqref{4.16} with \eqref{4.18}, \eqref{4.23}, \eqref{4.24} and \eqref{4.25} then yields that
\begin{equation}\label{4.26}
\frac{d}{dt}\left\|  \nabla \mathbf{ v}  \right\| ^{2}_{L^{2}(\mathbb{R}^{n})}+2\nu\left\| \Lambda^{\frac{n}{4}}\nabla  \mathbf{v}  \right\|^{2}_{L^{2}(\mathbb{R}^{n})}\leq C(1+t)^{-\frac{n}{4\beta}}\left\| \Lambda^{\frac{n}{4}}\nabla  \mathbf{v}  \right\| ^{2}_{L^{2}(\mathbb{R}^{n})}.\end{equation}
Therefore, from the above arguments of Case (I) and Case (II), for any $\frac{n}{4}\leq \beta <1$ with $n=2,3$,
choosing $t$ large enough such that $C(1+t)^{-\frac{n}{4\beta}}<\nu$,  one deduces from \eqref{4.22} and \eqref{4.26} that
\begin{equation}\label{4.27}
\frac{d}{dt}\left\| \nabla  \mathbf{v }\right \| ^{2}_{L^{2}(\mathbb{R}^{n})}+ \nu\left\| \Lambda^{\beta}\nabla \mathbf{ v}  \right\| ^{2}_{L^{2}(\mathbb{R}^{n})}\leq 0.\end{equation}
In the following, we continue our proof by applying the Fourier splitting method as used in the proof of (I) of this theorem.\\
  Let $B(\rho)$ be the ball of radius $\rho$, where
$\rho^{2\beta}= \dfrac{\frac{n}{2\beta}+\frac{1}{\beta}+1}{\nu(1+t)}$. Thanks to the Plancherel's theorem, it follows from \eqref{4.27} that
$$\frac{d}{dt}\left\| \xi \hat{\mathbf{v}}  \right\| ^{2}_{L^{2}(\mathbb{R}^{n})}+\nu \rho^{2\beta}\int_{B^{C}(\rho)}\left | \xi \hat{\mathbf{v}} \right  | ^{2}d\xi\leq 0,$$
which gives rise to
\begin{equation}\label{4.28}\frac{d}{dt}\left\| \xi \hat{\mathbf{v}}  \right\| ^{2}_{L^{2}(\mathbb{R}^{n})}+\nu \rho^{2\beta}\left\| \xi \hat{\mathbf{v}}  \right\| ^{2}_{L^{2}(\mathbb{R}^{n})}\leq\nu \rho^{2\beta+2}\int_{B(\rho)} \left| \hat{\mathbf{v}} \right | ^{2}d\xi.\end{equation}
On the other hand, we obtain by Lemma \ref{l4.2} and (I) of this theorem
$$\left| \hat{\mathbf{v}}  \right| ^{2}\leq C\left[1+ \int^{t}_{0}(1+s)^{-\frac{ n}{ 2\beta}}ds \cdot \int^{t}_{0}\| \nabla \mathbf{v}  \| ^{2}_{L^{2}(\mathbb{R}^{n})}ds \right].$$
With this bound and \eqref{4.28}, we arrive at
$$\frac{d}{dt}\left\| \xi \hat{\mathbf{v}} \right \| ^{2}_{L^{2}(\mathbb{R}^{n})}+\nu \rho^{2\beta}\left\| \xi \hat{\mathbf{v}}  \right\| ^{2}_{L^{2}(\mathbb{R}^{n})}\qquad\qquad \qquad\qquad\qquad\qquad\qquad\qquad\qquad$$
$$\leq C\nu \rho^{2\beta+2+n}\left[1+\left(\int^{t}_{0}(1+s)^{-\frac{ n}{ 2\beta}}ds\right)\left(\int^{t}_{0}\|\nabla  \mathbf{v}  \| ^{2}_{L^{2}(\mathbb{R}^{n})}ds\right)\right]. $$
Let $K(t)=(1+t)^{ \frac{ n}{ 2\beta}+\frac{ 1}{ \beta}+1}$. Taking $K(t)$ as an integrating factor, we then have
$$\frac{d}{dt}\left((1+t)^{ \frac{ n}{ 2\beta}+\frac{ 1}{ \beta}+1}\left\| \xi \hat{\mathbf{v}} \right \| ^{2}_{L^{2}(\mathbb{R}^{n})}\right) \qquad\qquad \qquad\qquad\qquad\qquad\qquad\qquad\qquad$$
$$\leq C \left[1+\left(\int^{t}_{0}(1+s)^{-\frac{ n}{ 2\beta}}ds\right)\left(\int^{t}_{0}\| \nabla \mathbf{v}  \| ^{2}_{L^{2}(\mathbb{R}^{n})}ds\right)\right]. $$
Thanks to the Tonelli theorem and the Plancherel's theorem,  we obtain by applying (I) of this theorem  again and integrating in time from $0$ to $r$
 \begin{eqnarray*}
  && (1+r)^{ \frac{ n}{ 2\beta}+\frac{ 1}{ \beta}+1}\| \nabla \mathbf{v} \| ^{2}_{L^{2}(\mathbb{R}^{n})}
  \\[0.3cm]
 &&\qquad \leq C (1+r)\\
 &&\qquad\quad +\int^{r}_{0}\left(\int^{t}_{0}(1+s)^{-\frac{ n}{ 2\beta}}ds\right)
  \cdot\left(\int^{t}_{0}\frac{(1+s)^{ \frac{ n}{ 2\beta}+\frac{ 1}{ \beta}+1}}{(1+s)^{ \frac{ n}{ 2\beta}+\frac{ 1}{ \beta}+1}}\| \nabla \mathbf{v}(s)  \| ^{2}_{L^{2}(\mathbb{R}^{n})}ds\right)dt.
\end{eqnarray*}
The Gronwall inequality then implies that
$$ (1+r)^{ \frac{ n}{ 2\beta}+\frac{ 1}{ \beta}+1}\| \nabla \mathbf{v }\| ^{2}_{L^{2}(\mathbb{R}^{n})} \leq C (1+r)e^{A},$$
where
$$A=\int^{r}_{0}\left(\int^{t}_{0}(1+s)^{-\frac{ n}{ 2\beta}}ds\right)\left(\int^{t}_{0} (1+s)^{ -\frac{ n}{ 2\beta}-\frac{ 1}{ \beta}-1} ds\right)dt.$$
Note that the term $A$ is bounded independent of $r$ for $n=2,3$, we then obtain
$$\left\| \nabla \mathbf{v}(r) \right\| ^{2}_{L^{2}(\mathbb{R}^{n})} \leq C(1+r)^{ -\frac{ n}{ 2\beta}-\frac{ 1}{ \beta}}.$$
This finishes the proof of (II).\\
\\
 {\bf We then prove (III-1)}.
 \\[0.3cm]
\indent Due to Lemma \ref{l2.10}, Lemma \ref{l4.2} with (I) and (II) of this theorem, we have $\displaystyle \left|{\cal{F}}(\mathbf{v})\right|\leq C$. Note that the Helmholtz equation $\mathbf{u}-\alpha^{2}\Delta \mathbf{u}=\mathbf{v}$, simple computation gives $\displaystyle \left|{\cal{F}}(\mathbf{u})\right|\leq\left|{\cal{F}}(\mathbf{v})\right|$ yields the conclusion of (III-1).
 \\[0.3cm]
 {\bf We next show (III-2)}.
 \\[0.3cm]
\indent From (I) of this theorem,  we have shown that
\begin{equation}\label{4.29}
\displaystyle \|\mathbf{u}\|_{L^{2} (\mathbb{R}^{n})}^{2}+
\alpha^{2}\|\nabla \mathbf{u}\|_{L^{2} (\mathbb{R}^{n})}^{2}\leq C (1+t)^{-\frac{n}{2\beta}}.
\end{equation}
Differentiating the Helmholtz equation $\mathbf{u}-\alpha^{2}\Delta \mathbf{u}=\mathbf{v}$ and squaring the resulting equation yields, after some integration by parts,
$$\displaystyle \left\|\nabla \mathbf{u}\right\|_{L^{2} (\mathbb{R}^{n})}^{2}+
2\alpha^{2}\left\|\nabla^{2} \mathbf{u}\right\|_{L^{2} (\mathbb{R}^{n})}^{2}+\alpha^{4}\left\|\nabla^{3} \mathbf{u}\right\|_{L^{2} (\mathbb{R}^{n})}^{2}=\left\|\nabla \mathbf{v}\right\|_{L^{2} (\mathbb{R}^{n})}^{2}.$$
Combining this with (II) of this theorem gives rise to
$$\left \|\nabla^{2} \mathbf{u}\right\|_{L^{2} (\mathbb{R}^{n})}^{2}\leq C (1+t)^{-\frac{n}{2\beta}-\frac{1}{\beta}}\leq C (1+t)^{-\frac{n}{2\beta}}.$$
This together with \eqref{4.29} deduces
$$\left\|\mathbf{v}\right\|_{L^{2} (\mathbb{R}^{n})}^{2}\leq \left\|  \mathbf{u}\right\|_{L^{2} (\mathbb{R}^{n})}^{2}+2\alpha^{2}\left\|\nabla  \mathbf{u}\right\|_{L^{2} (\mathbb{R}^{n})}^{2}+ \alpha^{4}\left\|\Delta \mathbf{u}\right\|_{L^{2} (\mathbb{R}^{n})}^{2}\leq C (1+t)^{ -\frac{n}{2\beta}}.$$
 This ends the proof of (III-2).\\
\\
{\bf In the following, we begin to show (IV)}.
\\[0.3cm]
\indent We will adopt the Fourier splitting argument again. Let $B(\rho)$ be the ball of radius $\rho$. Thanks to the Plancherel's theorem, breaking up the left hand side of the integral \eqref{4.1} deduces that
\begin{equation}\label{4.30}
\frac{1}{2}\frac{d}{dt}\left\| \xi^{m} \hat{w}  \right\| ^{2}_{L^{2}(\mathbb{R}^{n})}+\nu \rho^{2\beta}\left\| \xi^{m} \hat{w}  \right\| ^{2}_{L^{2}(\mathbb{R}^{n})}\leq\nu \rho^{2\beta+2m}\int_{B(\rho)}   |\hat{w}| ^{2}d\xi+C(1+t)^{\gamma}.
\end{equation}
Let $\rho^{2\beta}=\frac{b}{\nu(1+t)}$ for some large $b$. Note that the assumption for the bound on $\hat{w}$,  making direct calculation for the right hand side of \eqref{4.30} gives
\begin{equation}\label{4.31}
\frac{d}{dt}\left[(1+t)^{b}\| \xi^{m} \hat{w}  \| ^{2}_{L^{2}(\mathbb{R}^{n})}\right]
\leq C\left[(1+t)^{-\frac{m}{ \beta}-1+b+2\eta-\frac{n}{2\beta}}+(1+t)^{\gamma+b}\right].
\end{equation}
Integrating both sides of \eqref{4.31} with respect to time $t$, and then applying the plancherel's theorem once again, we arrive at the conclusion of (IV).\\
\\
{\bf We next prove (V)}.
\\[0.3cm]
\indent Note that the chain rule
  $$\frac{d}{dt}\int^{t}_{0}f(t,s)ds=f(t,t)+\int^{t}_{0}\frac{\partial f(t,s)}{\partial t}ds,$$
  one deduces from \eqref{4.5} and \eqref{4.6} that
  \begin{eqnarray*}
  &&\partial^{P}_{t} \mathcal{F}(\mathbf{v})=(-\nu)^{P}|\xi|^{2\beta P}e^{-\nu t|\xi|^{2\beta }}\mathcal{F}(\mathbf{v}_{0})\\
  &&\qquad\qquad\quad+\sum\limits^{P-1}_{p=0}\left(-\nu|\xi|^{2\beta}\right)^{P-1-p}\partial^{p}_{t}\psi(\xi,t)
  \\
   &&\qquad\qquad\quad+\displaystyle\int^{t}_{0}\left(-\nu|\xi|^{2\beta }\right)^{P }e^{-\nu (t-s)|\xi|^{2\beta }}\psi(\xi,s)ds.
      \end{eqnarray*}
With this expression, to achieve (V), the key ingredient is to first bound $\partial^{p}_{t}\psi(\xi,t)$, with $\psi(\xi,t)$ defined by \eqref{4.6}. Applying an argument similar to the proof of Lemma \ref{l4.2}, one obtains
\begin{equation*}
\left.
\begin{array}{ll}
\vspace{5pt}
\partial^{p}_{t} \psi(\xi,t)&=-\partial^{p}_{t}\mathcal{F}\left(\mathbf{u}\cdot \nabla \mathbf{v}\right)
+\partial^{p}_{t}\mathcal{F}\left(\mathbf{u}\cdot \nabla \mathbf{v}^{T}\right)+\partial^{p}_{t}
\left[-\xi \cdot\mathcal{F}\left(p+\sum\limits_{i=1}^{n} u_{i}v_{i} \right)\right]\\
&\triangleq\partial^{p}_{t}A+\partial^{p}_{t}B+\partial^{p}_{t}C.
\end{array}
\right.
\end{equation*}
Thanks to  $\mbox{div}~\mathbf{u}=0$, by the assumptions of (V), $\partial^{p}_{t}A$ can be bounded as follows:
 \begin{eqnarray*}
 \vspace{5pt}
  &&\displaystyle|\partial^{p}_{t}A|=\left|\partial^{p}_{t}\mathcal{F}
  \left(\sum\limits^{n}_{j=1}u_{j}\partial_{j}\mathbf{v}\right)
  \right|=\left|\partial^{p}_{t}\mathcal{F}\left(\sum\limits^{n}_{j=1}\partial_{j}(u_{j}\mathbf{v})\right)
  \right|  =\left|\partial^{p}_{t}\sum\limits^{n}_{j=1}
  \xi_{j}\mathcal{F}(u_{j}\mathbf{v})\right|\\
 \vspace{5pt}
  &&\displaystyle\qquad\quad\leq\sum\limits_{l=0}^{p}C|\xi|\left\|\partial^{l}_{t} \mathbf{v} \right \| _{L^{2}(\mathbb{R}^{n})}\left\|\partial^{p-l}_{t} \mathbf{v}   \right\| _{L^{2}(\mathbb{R}^{n})}\\
  \vspace{5pt}
   &&\displaystyle\qquad\quad\leq C(1+t)^{-\frac{ 1}{ 2\beta}}(1+t)^{-l-\frac{ n}{ 4\beta}}(1+t)^{-(p-l)-\frac{n}{ 4\beta}}\\
   \vspace{5pt}
   &&\displaystyle\qquad\quad\leq C(1+t)^{-p-\frac{ n}{ 2\beta}-\frac{ 1}{ 2\beta}}.
        \end{eqnarray*}
 In the same manner, one attains the bound for $\partial^{p}_{t}B$:
  \begin{eqnarray*}
  \vspace{5pt}
  &&\displaystyle\left|\partial^{p}_{t}B\right|=\left|\partial^{p}_{t}\sum\limits^{n}_{j=1}
  \mathcal{F}\left(u_{j} \nabla v_{j}\right)\right|\\
  \vspace{5pt}
  &&\displaystyle\qquad\quad\leq\sum\limits_{l=0}^{p}C|\xi|\left\|\partial^{l}_{t} \mathbf{v} \right \| _{L^{2}(\mathbb{R}^{n})}\left\|\partial^{p-l}_{t} \nabla \mathbf{v } \right\| _{L^{2}(\mathbb{R}^{n})}\\
   \vspace{5pt}
   &&\displaystyle\qquad\quad\leq C(1+t)^{-l-\frac{ n}{ 4\beta}} (1+t)^{-(p-l)-\frac{1}{ 2\beta}-\frac{n}{ 4\beta}}\\
   \vspace{5pt}
   &&\displaystyle\qquad\quad\leq C(1+t)^{-p-\frac{ n}{ 2\beta}-\frac{ 1}{ 2\beta}},
        \end{eqnarray*}
 Due to \eqref{4.9}, putting together the above estimates for $\partial^{p}_{t}A$ and $\partial^{p}_{t}B$, we get
\begin{eqnarray*}
   \vspace{5pt}
  &&\displaystyle\left|\partial^{p}_{t}C\right|=\left|\partial^{p}_{t}\xi
  \mathcal{F}\left(p+\sum\limits_{i=1}^{n} u_{i}v_{i}\right)\right|\leq \left|\partial^{p}_{t}A\right|+\left|\partial^{p}_{t}B\right|\\
     \vspace{5pt}
     &&\displaystyle\qquad\quad\leq  C(1+t)^{-p-\frac{ n}{ 2\beta}-\frac{ 1}{ 2\beta}}.\qquad\qquad
        \end{eqnarray*}
In view of \eqref{4.6}, the bound $|\hat{\mathbf{v}}|\leq C$ by (III-1) of this theorem, $\displaystyle|\xi|^{2\beta}\leq \frac{ b}{ \nu(1+t)}$ and $\displaystyle 1-\frac{ n}{ 2\beta}-\frac{ 1}{ 2\beta}<0$, we finish the proof of (V).\\
\\
 {\bf We further show (VI-1)}.
 \\[0.3cm]
\indent We shall prove this conclusion by using inductive argument. Due to the regularity of solutions (Proposition \ref{p2.4}), we present the proof only formally. It should be pointed out that the key point of the proof is to establish an inequality in a form satisfing the conclusion in (IV) of this theorem. To achieve this, we shall divide the proof into the following three steps.
\\[0.3cm]
{\bf Step 1.}\quad For $m=0,1$, the inequality holds by (III-2) and (II), respectively. That is,
$$\| \mathbf{v}  \| ^{2}_{L^{2}(\mathbb{R}^{n})}\leq C(1+t)^{-\frac{n}{2\beta}}, ~~\| \nabla  \mathbf{v}  \| ^{2}_{L^{2}(\mathbb{R}^{n})}\leq C(1+t)^{-\frac{1}{\beta}-\frac{n}{2\beta}}.$$
{\bf Step 2.}\quad We now assume (inductive assumption) that the decay
 \begin{equation}\label{4.32}\| \nabla^{m}  \mathbf{v}  \| ^{2}_{L^{2}(\mathbb{R}^{n})}\leq C(1+t)^{-\frac{m}{\beta}-\frac{n}{2\beta}}
 \end{equation}
 holds for all $m<M$. Here, $m$ and $M$ are both non-negative integers.
  \\[0.2cm]
{\bf Step 3.}\quad We will verify that the inequality \eqref{4.32} is true for $m=M$.\\
 \indent Multiplying the first equation in \eqref{1.1} by $\Delta^{M}\mathbf{v}$, and then integrating by parts the resulting equation gives rise to
\begin{equation}\label{4.33}
\left.
\begin{array}{ll}
\vspace{5pt}
\displaystyle
&\frac{d}{dt}\left\| \nabla^{M}  \mathbf{v} \right \| ^{2}_{L^{2}(\mathbb{R}^{n})}+2\nu\left\| \Lambda^{\beta}\nabla^{M}   \mathbf{v}  \right\| ^{2}_{L^{2}(\mathbb{R}^{n})}\\
\vspace{5pt}
&\qquad \leq\left|\left\langle \mathbf{u} \cdot\nabla \mathbf{v} , \Delta^{M} \mathbf{v}\right  \rangle\right|+
 \left|\left\langle \mathbf{v} \cdot\nabla \mathbf{u}^{T} , \Delta^{M} \mathbf{v} \right\rangle\right|\\
\vspace{5pt}&\qquad \triangleq ~~I_{M}+J_{M}.
\end{array}
\right.
\end{equation}
 To bound $I_{M}$ and $J_{M}$, there are two cases to consider.\\
  \\
  \indent  Case {\bf  (I)}\quad $\displaystyle\frac{n}{4}< \beta< 1$ for $n=2,3$;
  ~\\
  \indent Case {\bf  (II)}\quad $ \displaystyle \beta=\frac{n}{4}$ for $n=2,3$.\\
  \\
 We first consider  Case {\bf (I)}\quad $\displaystyle\frac{n}{4}< \beta< 1$ for $n=2,3$.\\
   \indent In this case, it is easy to check that $\displaystyle\frac{n}{2}-\beta<\beta$. Recall that \eqref{4.33} and $\left\langle \mathbf{u} \cdot\nabla \mathbf{v} ,  \mathbf{v} \right\rangle=0$, thanks to Cauchy's inequality, H\"{o}lder's inequality and Gagliardo-Nirenberg-Sobolev inequality, one deduces that\\
\begin{equation}\label{4.34}
\left.
\begin{array}{ll}
\displaystyle
    & \displaystyle I_{M}=  \left|\left\langle \mathbf{u} \cdot\nabla \mathbf{v} , \Delta^{M} \mathbf{v} \right\rangle\right|
    \\[0.3cm]
    &\displaystyle \qquad=    \left|\sum\limits^{M}_{m=1}\left(\begin{array}{l}
      M\\
      m
     \end{array}\right)\left\langle  \nabla^{m}\mathbf{u}\cdot \nabla\nabla^{M-m}\mathbf{v},  \nabla^{M} \mathbf{v} \right\rangle\right|
      \\[0.3cm]
      &\displaystyle\qquad \leq C\sum\limits^{M}_{m=1} \left\|  \nabla^{M+1-m}\mathbf{v} \right\|  _{L^{2}(\mathbb{R}^{n})}\left \| \nabla^{m}\mathbf{u}\cdot  \nabla^{M }\mathbf{v} \right\|  _{L^{2}(\mathbb{R}^{n})}
       \\[0.3cm]
      &\displaystyle\qquad \leq C\sum\limits^{M}_{m=1} \left\|  \nabla^{M+1-m}\mathbf{v} \right\|  _{L^{2}(\mathbb{R}^{n})} \left\| \nabla^{m}\mathbf{u}  \right\|  _{L^{\frac{n}{\beta}}(\mathbb{R}^{n})} \left\| \nabla^{M}\mathbf{v}  \right\|  _{L^{\frac{2n}{n-2\beta}}(\mathbb{R}^{n})}
       \\[0.3cm]
      &\displaystyle\qquad \leq C\sum\limits^{M}_{m=1} \left\|  \nabla^{M+1-m}\mathbf{v}\right \|  _{L^{2}(\mathbb{R}^{n})} \left\| \Lambda^{\beta}\nabla^{m}\mathbf{u}  \right\|  _{L^{2}(\mathbb{R}^{n})} \left\| \Lambda^{\beta}\nabla^{M}\mathbf{v } \right\|  _{L^{2}(\mathbb{R}^{n})}
       \\[0.3cm]
       &\displaystyle\qquad \leq C\sum\limits^{M}_{m=1} \left\|  \nabla^{M+1-m}\mathbf{v} \right\|^{2}  _{L^{2}(\mathbb{R}^{n})} \left\| \Lambda^{\beta}\nabla^{m-1}\mathbf{v}  \right\|^{2}  _{L^{2}(\mathbb{R}^{n})} +\frac{\nu}{2}\left\| \Lambda^{\beta}\nabla^{M}\mathbf{v}  \right\|^{2}    _{L^{2}(\mathbb{R}^{n})},
    \end{array}
\right.
\end{equation}
and
 \begin{equation}\label{4.35}
\left.
\begin{array}{ll}
\displaystyle
    & \displaystyle J_{M}= \left|\left\langle \mathbf{v} \cdot\nabla \mathbf{u}^{T} , \Delta^{M} \mathbf{v} \right\rangle\right|
    \\[0.3cm]
    &\displaystyle\qquad=   \sum\limits^{M}_{m=0}\left(\begin{array}{l}
      M\\
      m
     \end{array}\right) \left|\left\langle  \nabla^{M}\mathbf{v}\cdot \nabla\nabla^{m}\mathbf{u},  \nabla^{M-m} \mathbf{v} \right\rangle\right|
     \\[0.3cm]
     &\displaystyle\qquad \leq C \left|\left\langle  \nabla^{M}\mathbf{v}\cdot \nabla \mathbf{u},  \nabla^{M } \mathbf{v} \right\rangle\right|+C\sum\limits^{M}_{m=1}\left|\left\langle  \nabla^{M}\mathbf{v}\cdot \nabla^{m+1} \mathbf{u},  \nabla^{M -m} \mathbf{v} \right\rangle\right|.
     \end{array}
\right.
\end{equation}
 By Lemma \ref{l2.12}, a straightforward computation shows that \begin{equation}\label{4.36}
\left.
\begin{array}{ll}
\vspace{5pt}
&\displaystyle\left|\left\langle  \nabla^{M}\mathbf{v}\cdot \nabla \mathbf{u},  \nabla^{M } \mathbf{v} \right\rangle\right|\\
\vspace{5pt}
& \displaystyle\qquad\leq \left\|  \nabla^{M}\mathbf{v} \right\| ^{2} _{L^{4}(\mathbb{R}^{n})}\left\|  \nabla \mathbf{u} \right\|  _{L^{2}(\mathbb{R}^{n})}\\
\vspace{5pt}
 &\displaystyle \qquad \leq \left(C(\varepsilon)\left\|  \nabla^{M}\mathbf{v} \right\| ^{2} _{L^{2}(\mathbb{R}^{n})}+\varepsilon\left\| \Lambda^{\beta} \nabla^{M}\mathbf{v} \right\| ^{2} _{L^{2}(\mathbb{R}^{n})}\right) \left\|  \nabla \mathbf{u} \right\|  _{L^{2}(\mathbb{R}^{n})},
     \end{array}
\right.
\end{equation}
and
 \begin{equation}\label{4.37}
\left.
\begin{array}{ll}
\vspace{5pt}
&\displaystyle\sum\limits^{M}_{m=1}\left|\left\langle  \nabla^{M}\mathbf{v}\cdot \nabla^{m+1} \mathbf{u},  \nabla^{M -m} \mathbf{v} \right\rangle\right|\\
\vspace{5pt}
&\displaystyle\qquad\leq \sum\limits^{M}_{m=1}\left\|  \nabla^{M-m}\mathbf{v} \right\|  _{L^{2}(\mathbb{R}^{n})}
\left\|  \nabla^{m+1}\mathbf{u}\cdot\nabla^{M}\mathbf{v} \right\|  _{L^{2}(\mathbb{R}^{n})}\\
\vspace{5pt}
&\displaystyle\qquad\leq \sum\limits^{M}_{m=1}\left\|  \nabla^{M-m}\mathbf{v} \right\|  _{L^{2}(\mathbb{R}^{n})}
\left\|  \nabla^{m+1}\mathbf{u}  \right\|  _{L^{\frac{n}{\beta}}(\mathbb{R}^{n})}\left\|  \nabla^{M}\mathbf{v}  \right\|  _{L^{\frac{2n}{n-2\beta}}(\mathbb{R}^{n})}\\
\vspace{5pt}
&\displaystyle\qquad\leq \sum\limits^{M}_{m=1}\left\|  \nabla^{M-m}\mathbf{v} \right\|  _{L^{2}(\mathbb{R}^{n})}
\left\|  \Lambda^{\beta}\nabla^{m+1}\mathbf{u}  \right\|  _{L^{2}(\mathbb{R}^{n})}\left\| \Lambda^{\beta} \nabla^{M}\mathbf{v } \right\|  _{L^{2}(\mathbb{R}^{n})}\\
\vspace{5pt}
&\displaystyle\qquad\leq \sum\limits^{M}_{m=1}\left\|  \nabla^{M-m}\mathbf{v} \right\|  _{L^{2}(\mathbb{R}^{n})}
\left\|  \Lambda^{\beta}\nabla^{m-1}\mathbf{v} \right \|  _{L^{2}(\mathbb{R}^{n})}\left\| \Lambda^{\beta} \nabla^{M}\mathbf{v} \right \|  _{L^{2}(\mathbb{R}^{n})}\\
\vspace{5pt}
&\displaystyle\qquad\leq C\sum\limits^{M}_{m=1}\left\|  \nabla^{M-m}\mathbf{v} \right\|^{2}  _{L^{2}(\mathbb{R}^{n})}
\left\|  \Lambda^{\beta}\nabla^{m-1}\mathbf{v}  \right\|^{2}  _{L^{2}(\mathbb{R}^{n})}+\frac{\nu}{4}\left\| \Lambda^{\beta} \nabla^{M}\mathbf{v}  \right\| ^{2} _{L^{2}(\mathbb{R}^{n})}.
   \end{array}
\right.
\end{equation}
 We have used the relation $\mathbf{u}-\alpha^{2}\Delta \mathbf{u}=\mathbf{v}$ in the estimates \eqref{4.34} and \eqref{4.37}. Choosing $\varepsilon\leq \dfrac{\nu}{4\|\nabla \mathbf{u}_{0}  \|^{2}  _{L^{2}(\mathbb{R}^{n})}}$, if follows from \eqref{4.33}, \eqref{4.34}, \eqref{4.35}, \eqref{4.36} and \eqref{4.37} that
 \begin{equation}\label{4.38}
\left.
\begin{array}{ll}
&\displaystyle\frac{d}{dt}\left\| \nabla^{M}  \mathbf{v}  \right\| ^{2}_{L^{2}(\mathbb{R}^{n})}+ \nu\left\| \Lambda^{\beta}\nabla^{M}   \mathbf{v}  \right\| ^{2}_{L^{2}(\mathbb{R}^{n})}
 \\[0.3cm]
 &\displaystyle\qquad \leq C\left\| \nabla^{M}  \mathbf{v}  \right\| ^{2}_{L^{2}(\mathbb{R}^{n})}\left\| \Lambda^{\beta}\mathbf{v}     \right\| ^{2}_{L^{2}(\mathbb{R}^{n})}+C(\varepsilon)\left\| \nabla^{M}  \mathbf{v } \right\| ^{2}_{L^{2}(\mathbb{R}^{n})}\left\| \nabla \mathbf{u}     \right\| ^{2}_{L^{2}(\mathbb{R}^{n})}
  \\[0.3cm]
 &\displaystyle\qquad \quad+\left\| \nabla^{M-1}  \mathbf{v}  \right\| ^{2}_{L^{2}(\mathbb{R}^{n})}\left\| \Lambda^{\beta}\mathbf{v}  \right \| ^{2}_{L^{2}(\mathbb{R}^{n})}+C \sum\limits^{M}_{m=2}\left\| \nabla^{M-m}  \mathbf{v} \right \| ^{2}_{L^{2}(\mathbb{R}^{n})}\left\| \Lambda^{\beta}\nabla^{m-1}\mathbf{v} \right\| ^{2}_{L^{2}(\mathbb{R}^{n})}.
  \end{array}
\right.
\end{equation}
We next consider Case {\bf  (II)}\quad $ \displaystyle \beta=\frac{n}{4}$ for $n=2,3$.
\\[0.3cm]
\indent In this case, thanks to Lemma \ref{l2.12}, note that \eqref{4.33} and $\left\langle \mathbf{u}\cdot \nabla \mathbf{v},   \mathbf{v} \right\rangle=0$, we deduce the following two estimates:
\begin{equation}\label{4.39}
\left.
\begin{array}{ll}
\displaystyle
    &\displaystyle I_{M}=  \left|\left\langle \mathbf{u} \cdot\nabla \mathbf{v} , \Delta^{M} \mathbf{v} \right\rangle\right|
    \\[0.3cm]
    &\displaystyle\qquad=    \left|\sum\limits^{M}_{m=1}\left(\begin{array}{l}
      M\\
      m
     \end{array}\right)\left\langle  \nabla^{m}\mathbf{u}\cdot \nabla\nabla^{M-m}\mathbf{v},  \nabla^{M} \mathbf{v} \right \rangle\right|
     \\[0.3cm]
      &\displaystyle\qquad \leq C\sum\limits^{M}_{m=1} \left\|  \nabla^{M+1-m}\mathbf{v} \right\|  _{L^{2}(\mathbb{R}^{n})} \left\| \nabla^{m}\mathbf{u}\cdot  \nabla^{M }\mathbf{v } \right\|  _{L^{2}(\mathbb{R}^{n})}
      \\[0.3cm]
      &\displaystyle\qquad \leq C\sum\limits^{M}_{m=1}\left \|  \nabla^{M+1-m}\mathbf{v} \right\|  _{L^{2}(\mathbb{R}^{n})}\left \| \nabla^{m}\mathbf{u} \right \|  _{L^{4}(\mathbb{R}^{n})} \left\| \nabla^{M}\mathbf{v}  \right\|  _{L^{4}(\mathbb{R}^{n})}
      \\[0.3cm]
      &\displaystyle\qquad \leq C\sum\limits^{M}_{m=1}\left \|  \nabla^{M+1-m}\mathbf{v }\right\|  _{L^{2}(\mathbb{R}^{n})}\left(\left\| \nabla^{m}\mathbf{u}  \right\|  ^{2}_{L^{2}(\mathbb{R}^{n})}+ \left\| \Lambda^{\frac{n}{4}}\nabla^{m-1}\mathbf{v}  \right\| ^{2} _{L^{2}(\mathbb{R}^{n})}\right)^{\frac{1}{2}}
      \\[0.3cm]
        &\displaystyle\qquad\qquad\qquad\qquad\qquad
      \cdot\left( \left\| \nabla^{M}\mathbf{v} \right \|  ^{2}_{L^{2}(\mathbb{R}^{n})}+ \left\| \Lambda^{\frac{n}{4}}\nabla^{M}\mathbf{v} \right \|  ^{2}_{L^{2}(\mathbb{R}^{n})}\right)^{\frac{1}{2}},
    \end{array}
\right.
\end{equation}
 and
 \begin{equation}\label{4.40}
\left.
\begin{array}{ll}
\displaystyle
    &\displaystyle J_{M}= \left|\left\langle \mathbf{v} \cdot\nabla \mathbf{u}^{T} , \Delta^{M} \mathbf{v} \right\rangle\right|
    \\[0.3cm]
    &\displaystyle\qquad=   \sum\limits^{M}_{m=0}\left(\begin{array}{l}
      M\\
      m
     \end{array}\right) \left|\left\langle  \nabla^{M}\mathbf{v}\cdot \nabla\nabla^{m}\mathbf{u},  \nabla^{M-m} \mathbf{v } \right\rangle\right|
     \\[0.3cm]
     &\displaystyle\qquad \leq C \sum\limits^{M}_{m=0}\left\|  \nabla^{M-m}\mathbf{v} \right\|  _{L^{2}(\mathbb{R}^{n})} \left\| \nabla^{m+1}\mathbf{u}\cdot \nabla^{M}\mathbf{v}  \right\|  _{L^{2}(\mathbb{R}^{n})}
     \\[0.3cm]
       &\displaystyle\qquad \leq C \sum\limits^{M}_{m=0}\left\|  \nabla^{M-m}\mathbf{v} \right\|  _{L^{2}(\mathbb{R}^{n})}\left \| \nabla^{m+1}\mathbf{u} \right\|  _{L^{4}(\mathbb{R}^{n})} \left\| \nabla^{M}\mathbf{v}\right\|  _{L^{4}(\mathbb{R}^{n})}
       \\[0.3cm]
      &\displaystyle\qquad \leq C \sum\limits^{M}_{m=0}\left\|  \nabla^{M-m}\mathbf{v}\right \|  _{L^{2}(\mathbb{R}^{n})}
      \left(\left\| \nabla^{m+1}\mathbf{u}\right\|  ^{2}_{L^{2}(\mathbb{R}^{n})}+ \left\| \Lambda^{\frac{n}{4}}\nabla^{m+1}\mathbf{u}  \right\| ^{2} _{L^{2}(\mathbb{R}^{n})}\right)^{\frac{1}{2}}
      \\[0.3cm]
        &\displaystyle\qquad\qquad\qquad\qquad\qquad
      \cdot
      \left(\left\| \nabla^{M}\mathbf{v}\right\|  ^{2}_{L^{2}(\mathbb{R}^{n})}+ \left\| \Lambda^{\frac{n}{4}} \nabla^{M}\mathbf{v}  \right\| ^{2} _{L^{2}(\mathbb{R}^{n})}\right)^{\frac{1}{2}}
      \\[0.3cm]
        &\displaystyle\qquad \leq C \sum\limits^{M}_{m=0}\left\|  \nabla^{M-m}\mathbf{v}\right \|  _{L^{2}(\mathbb{R}^{n})}
      \left(\left\| \nabla^{m}\mathbf{v}\right\|  ^{2}_{L^{2}(\mathbb{R}^{n})}+ \left\| \Lambda^{\frac{n}{4}}\nabla^{m}\mathbf{v}  \right\| ^{2} _{L^{2}(\mathbb{R}^{n})}\right)^{\frac{1}{2}}
       \\[0.3cm]
        &\displaystyle\qquad\qquad\qquad\qquad
      \qquad
      \cdot\left(\left\| \nabla^{M}\mathbf{v}\right\|  ^{2}_{L^{2}(\mathbb{R}^{n})}+ \left\| \Lambda^{\frac{n}{4}} \nabla^{M}\mathbf{v}  \right\| ^{2} _{L^{2}(\mathbb{R}^{n})}\right)^{\frac{1}{2}}.
     \end{array}
\right.
\end{equation}
 Combining \eqref{4.39} with \eqref{4.40} yields that
 \begin{equation}\label{4.41}
\left.
\begin{array}{ll}
&\displaystyle I_{M}+J_{M}\leq C \sum\limits^{M}_{m=0} \left\|  \nabla^{M-m}\mathbf{v} \right\|  _{L^{2}(\mathbb{R}^{n})}\left(\left\| \nabla^{m}\mathbf{v} \right\|  ^{2}_{L^{2}(\mathbb{R}^{n})}+ \left\| \Lambda^{\frac{n}{4}} \nabla^{m}\mathbf{v}   \right\| ^{2} _{L^{2}(\mathbb{R}^{n})}\right)^{\frac{1}{2}}
 \\[0.3cm]
        &\displaystyle\qquad\qquad\qquad\qquad
      \qquad
      \cdot\left(\left\| \nabla^{M}\mathbf{v}\right\|  ^{2}_{L^{2}(\mathbb{R}^{n})}+ \left\| \Lambda^{\frac{n}{4}} \nabla^{M}\mathbf{v}  \right\| ^{2} _{L^{2}(\mathbb{R}^{n})}\right)^{\frac{1}{2}}
\\[0.4cm]
 &\displaystyle\qquad\qquad \leq C \left\|\nabla^{M}\mathbf{v}\right \|  _{L^{2}(\mathbb{R}^{n})}
  \left(\left\|  \mathbf{v}\right\|  ^{2}_{L^{2}(\mathbb{R}^{n})}+ \left\| \Lambda^{\frac{n}{4}} \mathbf{v}  \right\| ^{2} _{L^{2}(\mathbb{R}^{n})}\right)^{\frac{1}{2}}
  \\[0.3cm]
        &\displaystyle\qquad\qquad\qquad\qquad
      \qquad
      \cdot\left(\left\| \nabla^{M} \mathbf{v}\right\|  ^{2}_{L^{2}(\mathbb{R}^{n})}+ \left\| \Lambda^{\frac{n}{4}}\nabla^{M} \mathbf{v}  \right\| ^{2} _{L^{2}(\mathbb{R}^{n})}\right)^{\frac{1}{2}}
 \\[0.3cm]
 &\displaystyle\qquad\qquad\quad+ C \sum\limits^{M}_{m=1} \left\|\nabla^{M-m}\mathbf{v} \right\|  _{L^{2}(\mathbb{R}^{n})}
 \left(\left\|\nabla^{m}\mathbf{v}\right\|  ^{2}_{L^{2}(\mathbb{R}^{n})}+ \left\| \Lambda^{\frac{n}{4}} \nabla^{m}\mathbf{v} \right\| ^{2} _{L^{2}(\mathbb{R}^{n})}\right)^{\frac{1}{2}}
  \\[0.3cm]
        &\displaystyle\qquad\qquad\qquad\qquad
      \qquad
      \cdot\left(\left\|\nabla^{M}\mathbf{v}\right\|  ^{2}_{L^{2}(\mathbb{R}^{n})}+ \left\| \Lambda^{\frac{n}{4}} \nabla^{M}\mathbf{v} \right\| ^{2} _{L^{2}(\mathbb{R}^{n})}\right)^{\frac{1}{2}}
 \\[0.3cm]
 &\displaystyle\qquad\qquad\quad+ C \left  \|\mathbf{v} \right\|  _{L^{2}(\mathbb{R}^{n})}  \left(\left\|\nabla^{M}\mathbf{v}\right\|  ^{2}_{L^{2}(\mathbb{R}^{n})}+ \left\| \Lambda^{\frac{n}{4}} \nabla^{M}\mathbf{v} \right\| ^{2} _{L^{2}(\mathbb{R}^{n})}\right)^{\frac{1}{2}}
  \\[0.3cm]
 &\displaystyle\qquad\qquad \lesssim \left\|\nabla^{M}\mathbf{v} \right\|^{2}_{L^{2}(\mathbb{R}^{n})} \left(\left\| \mathbf{v}\right\|  ^{2}_{L^{2}(\mathbb{R}^{n})}+ \left\| \Lambda^{\frac{n}{4}}  \mathbf{v} \right\| ^{2} _{L^{2}(\mathbb{R}^{n})}\right)
 \\[0.3cm]
 &\displaystyle\qquad\qquad\quad+ \frac{\nu}{2} \left(\left\|\nabla^{M}\mathbf{v}\right\|  ^{2}_{L^{2}(\mathbb{R}^{n})}+ \left\| \Lambda^{\frac{n}{4}} \nabla^{M}\mathbf{v} \right\| ^{2} _{L^{2}(\mathbb{R}^{n})}\right)
   \\[0.3cm]
 &\displaystyle\qquad\qquad\quad+ C \sum\limits^{M-1}_{m=1} \left\|\nabla^{M-m}\mathbf{v} \right\| ^{2} _{L^{2}(\mathbb{R}^{n})} \left(\left\|\nabla^{m-1}\mathbf{v}\right\|  ^{2}_{L^{2}(\mathbb{R}^{n})}+ \left\| \Lambda^{\frac{n}{4}} \nabla^{m-1}\mathbf{v} \right\| ^{2} _{L^{2}(\mathbb{R}^{n})}\right)
 \\[0.3cm]
 &\displaystyle\qquad\qquad\quad+ C \left \|\mathbf{v}\right \|  _{L^{2}(\mathbb{R}^{n})}  \left(\left\|\nabla^{M}\mathbf{v}\right\|  ^{2}_{L^{2}(\mathbb{R}^{n})}+ \left\| \Lambda^{\frac{n}{4}} \nabla^{M}\mathbf{v} \right\| ^{2} _{L^{2}(\mathbb{R}^{n})}\right).
 \end{array}
\right.
\end{equation}
This together with \eqref{4.33} and the smallness assumption of the initial data for $\beta=\dfrac{n}{4}$ in (VI) with $M\leq K$ ensures that
\begin{equation}\label{4.42}
\left.
\begin{array}{ll}
&\displaystyle\frac{d}{dt}\left \| \nabla^{M}  \mathbf{v} \right \| ^{2}_{L^{2}(\mathbb{R}^{n})}+ \nu\left \| \Lambda^{\frac{n}{4}}\nabla^{M}   \mathbf{v} \right \| ^{2}_{L^{2}(\mathbb{R}^{n})}
 \\[0.3cm]
&\displaystyle\qquad \lesssim\left \| \nabla^{M}  \mathbf{v} \right \| ^{2}_{L^{2}(\mathbb{R}^{n})}\left \| \Lambda^{\frac{n}{4}}\mathbf{v}     \right\| ^{2}_{L^{2}(\mathbb{R}^{n})}
+\sum\limits^{M-1}_{m=1} \left \| \nabla^{M-m}  \mathbf{v}\right  \| ^{2}_{L^{2}(\mathbb{R}^{n})}\left \| \Lambda^\frac{n}{4}\nabla^{m-1}\mathbf{v}  \right  \| ^{2}_{L^{2}(\mathbb{R}^{n})}.
\end{array}
\right.
\end{equation}
Note that $\left\| \nabla^{M}  \mathbf{v } \right\| ^{2}_{L^{2}(\mathbb{R}^{n})}\leq C$ by Proposition \ref{p2.4}, (I), (II) of this theorem and the inductive assumption \eqref{4.32}, applying interpolation inequality and a bootstrap argument, it follows from \eqref{4.38} and \eqref{4.42} that for $\displaystyle \frac{n}{4}\leq \beta<1$ with $n=2,3$,
\begin{equation}\label{4.43}
\left.
\begin{array}{ll}
&\displaystyle\frac{d}{dt}\left\| \nabla^{M}  \mathbf{v}  \right\| ^{2}_{L^{2}(\mathbb{R}^{n})}+ \nu\left\| \Lambda^{\beta}\nabla^{M}  \mathbf{ v}  \right\| ^{2}_{L^{2}(\mathbb{R}^{n})}
 \\[0.3cm]
&\displaystyle\qquad \leq C(1+t)^{-1-\frac{n}{2\beta} }\left\| \nabla^{M}  \mathbf{v}  \right\| ^{2}_{L^{2}(\mathbb{R}^{n})}+C(1+t)^{-\frac{M}{\beta} -\frac{n}{2\beta}-1 }(1+t)^{ -\frac{n}{2\beta}+\frac{1}{\beta} }
 \\[0.3cm]
&\displaystyle\qquad\leq C(1+t)^{-\frac{M}{\beta}-\frac{n}{2\beta}-1 }.
\end{array}
\right.
\end{equation}
Here, we applied the fact that for $n=2,3$, $ -\frac{n}{2\beta}+\frac{1}{\beta} \leq 0$ and $(1+t)^{ -\frac{n}{2\beta}+\frac{1}{\beta} }\leq C$. Since $\left|\mathcal{F}(\mathbf{v})\right|\leq C$ by (III-1) of this theorem, applying (IV) of this theorem to estimate \eqref{4.43} gives rise to
$$\left\| \nabla^{M}  \mathbf{v}  \right\| ^{2}_{L^{2}(\mathbb{R}^{n})}\leq C(1+t)^{-\frac{M}{\beta}-\frac{n}{2\beta}}.  $$
With the inductive assumption \eqref{4.32}, thanks to (I) and (II) of this theorem, applying interpolation inequality, it follows form \eqref{4.38} and \eqref{4.42} that for
$\displaystyle \frac{n}{4}\leq \beta<1$ with $n=2,3$,
\begin{equation*}
\left.
\begin{array}{ll}
&\displaystyle\frac{d}{dt}\left\| \nabla^{M} \mathbf{ v}  \right\| ^{2}_{L^{2}(\mathbb{R}^{n})}+ \nu\left\| \Lambda^{\beta}\nabla^{M}   \mathbf{v}  \right\| ^{2}_{L^{2}(\mathbb{R}^{n})}
 \\[0.3cm]
&\displaystyle\qquad \leq C(1+t)^{-1-\frac{n}{2\beta} }\left\| \nabla^{M} \mathbf{ v}  \right\| ^{2}_{L^{2}(\mathbb{R}^{n})}+C(1+t)^{-\frac{M-m}{\beta} -\frac{n}{2\beta} }
(1+t)^{-\frac{m-1+\beta}{\beta} -\frac{n}{2\beta} }
 \\[0.3cm]
&\displaystyle\qquad\leq C(1+t)^{-1-\frac{n}{2\beta} }\left\| \nabla^{M} \mathbf{ v}  \right\| ^{2}_{L^{2}(\mathbb{R}^{n})}+C(1+t)^{-\frac{M}{\beta} -\frac{n}{\beta}- \frac{-1+\beta}{\beta}},
\end{array}
\right.
\end{equation*}
which together with the fact that $\displaystyle \frac{1}{\beta}-\frac{n}{2\beta}\leq 0$ implies that
\begin{equation}\label{4.44}
 \displaystyle\frac{d}{dt}\left\| \nabla^{M}  \mathbf{v}  \right\| ^{2}_{L^{2}(\mathbb{R}^{n})}+ \nu\left\| \Lambda^{\beta}\nabla^{M}   \mathbf{v } \right\| ^{2}_{L^{2}(\mathbb{R}^{n})}
 \leq C(1+t)^{-1-\frac{n}{2\beta} }.
\end{equation}
Note that $\displaystyle \left\| \nabla^{M}  \mathbf{v}  \right\| ^{2}_{L^{2}(\mathbb{R}^{n})}\leq C$ by Proposition \ref{p2.4}, with a bootstrap argument, and $\left|\mathcal{F}(v)\right|\leq C$ by (III-1) of this theorem, applying (IV) of this theorem to the estimate \eqref{4.44} deduces
\begin{equation}\label{4.45}
\left.
\begin{array}{ll}
\displaystyle\left\| \nabla^{M}  \mathbf{v}  \right\| ^{2}_{L^{2}(\mathbb{R}^{n})}&\leq C\left[(1+t)^{ -\frac{M}{\beta}-\frac{n}{2\beta} }+(1+t)^{ -\frac{n}{2\beta} }
(1+t)^{-\frac{m-1+\beta}{\beta} -\frac{n}{2\beta} }\right]
 \\[0.3cm]
\displaystyle&\leq C(1+t)^{ -\frac{n}{2\beta} }.
\end{array}
\right.
\end{equation}
In the same manner, using a bootstrap argument again and placing \eqref{4.45} into \eqref{4.43} yields that
\begin{equation*}
\left.
\begin{array}{ll}
&\displaystyle\frac{d}{dt}\left\| \nabla^{M}  \mathbf{v } \right\| ^{2}_{L^{2}(\mathbb{R}^{n})}+ \nu\left\| \Lambda^{\beta}\nabla^{M}   \mathbf{v}  \right\| ^{2}_{L^{2}(\mathbb{R}^{n})}
 \\[0.3cm]
&\displaystyle\qquad \leq C(1+t)^{-1-\frac{n}{\beta} } +C(1+t)^{-\frac{M }{\beta} -\frac{n}{ \beta}-1 }
 \\[0.3cm]
&\displaystyle\qquad\leq C(1+t)^{-1-\frac{n}{\beta} }.
\end{array}
\right.
\end{equation*}
Making the same argument as that used in \eqref{4.44} and \eqref{4.45}, we obtain
\begin{equation*}\displaystyle\left\| \nabla^{M}  \mathbf{v } \right\| ^{2}_{L^{2}(\mathbb{R}^{n})} \leq C (1+t)^{ -\frac{n}{2\beta}-\frac{n}{2\beta} }.
  \end{equation*}
Continuing with a bootstrap argument again, and using (IV) of this theorem, we deduce that
\begin{equation*}\displaystyle\left\| \nabla^{M}  \mathbf{v}  \right\| ^{2}_{L^{2}(\mathbb{R}^{n})} \leq C (1+t)^{ -\frac{M}{\beta}-\frac{n}{2\beta} }.
  \end{equation*}
Combinig Step 1 with Step 2 and Step 3 finishes the proof of (VI-1).\\
\\
{\bf We finally show the last conclusion (VI-2)}.
\\[0.3cm]
 We will adopt an inductive argument as above. The inductive assumption is as follows. \\
\indent For $p\leq\frac{K}{2\beta}$, the decay rate
\begin{equation}\label{4.46}
\left\|\partial^{p}_{t} \nabla^{m}  \mathbf{v} \right \| ^{2}_{L^{2}(\mathbb{R}^{n})}\leq C(1+t)^{-2p-\frac{n}{2\beta}-\frac{m}{\beta}} \end{equation}
holds for all $p<P$ and $m$ such that $2p\beta+m\leq K$. Here, $p$, $P$ and $m$ are all non-negative integers.
\\[0.3cm]
In the following, based on the inductive assumption \eqref{4.46}, we divided the proof into four steps. In Step 1, we show that for $\displaystyle|\xi|^{2\beta}\leq\frac{b}{\nu(1+t)}$, $\displaystyle\left\|\partial^{p}_{t} \hat{\mathbf{v}}(\xi)\right\|\leq C(1+t)^{-P}$. In the second step, we verify that the decay rate \eqref{4.46} holds for $p=P$ and $m=0$ by an inductive argument on $p$. We will check the decay rate \eqref{4.46} holds for any $m>0$ by another inductive argument on $m$ in the third step. In the fourth step, we conclude the expected result by a bootstrap argument.\\
\indent We begin to show (VI-2) step by step in detail.\\
{\bf Step 1}\quad We show for $\displaystyle|\xi|^{2\beta}\leq\frac{b}{\nu(1+t)}$, $\displaystyle\left\|\partial^{p}_{t} \hat{\mathbf{v}}(\xi)\right\|\leq C(1+t)^{-P}$.\\
\indent By \eqref{4.46} we get for all $p<P$ and $m=0,1$,
\begin{equation}\label{4.47}
\left\|\partial^{p}_{t} \nabla^{m}  \mathbf{v} \right \| ^{2}_{L^{2}(\mathbb{R}^{n})}\leq C(1+t)^{-2p-\frac{m}{\beta}-\frac{n}{2\beta}}. \end{equation}
By the aid of (V) of this theorem, \eqref{4.47} implies that for $\displaystyle|\xi|^{2\beta}\leq\frac{b}{\nu(1+t)}$,
\begin{equation}\label{4.48}
\left\| \partial^{P}_{t}   \hat{\mathbf{v}}(\xi) \right \| ^{2}_{L^{2}(\mathbb{R}^{n})}\leq C(1+t)^{-P}. \end{equation}
{\bf Step 2}\quad We now show that the decay rate \eqref{4.46} holds for $p=P$ and $m=0$ by an inductive argument on $p$.
\\[0.3cm]
\indent Note that $\mathbf{v}\cdot\nabla \mathbf{u}^{T}=\nabla(\mathbf{u}\mathbf{v})-\mathbf{u}\cdot \nabla \mathbf{v}^{T}$ by \eqref{4.4} and $\mbox{div} ~\mathbf{v}=0$, choosing $P$ and $M$ such that $M+2P\leq K$, then applying $\partial^{P}_{t}$ to the first equation in \eqref{1.1}, multiplying the resulting equation by  $\partial^{P}_{t}\Delta^{M}\mathbf{v}$ and integrating in space variable $x$ yields, after some integration by parts,
\begin{equation}\label{4.49}
\left.
\begin{array}{ll}
 &\displaystyle\frac{d}{dt}\left\| \partial^{P}_{t} \nabla^{M}\mathbf{v} \right \| ^{2}_{L^{2}(\mathbb{R}^{n})}+2\nu\left\| \partial^{P}_{t}\nabla^{M}\Lambda^{\beta} \mathbf{v} \right \| ^{2}_{L^{2}(\mathbb{R}^{n})}
\\[0.3cm]
&\displaystyle\qquad \leq\left|\left\langle \partial^{P}_{t}(\mathbf{u} \cdot\nabla \mathbf{v}) , \partial^{P}_{t}\Delta^{M} \mathbf{v} \right\rangle\right|+\left|\left\langle \partial^{P}_{t}\left(\mathbf{v} \cdot\nabla \mathbf{u}^{T}\right) , \partial^{P}_{t}\Delta^{M} \mathbf{v }\right\rangle\right|
\\[0.3cm]
&\displaystyle\qquad =\left|\left\langle \partial^{P}_{t}(\mathbf{u} \cdot\nabla \mathbf{v}) , \partial^{P}_{t}\Delta^{M} \mathbf{v} \right\rangle\right|+\left|\left\langle \partial^{P}_{t}\left(\mathbf{u} \cdot\nabla \mathbf{v}^{T}\right) , \partial^{P}_{t}\Delta^{M} \mathbf{v} \right\rangle\right|
\\[0.3cm]
&\displaystyle\qquad \triangleq I_{M,P}+J_{M,P}.
\end{array}
\right.
\end{equation}
In the following, we deal with the two terms on the right hand side of \eqref{4.49} by considering two cases:
\\[0.3cm]
  \indent Case {\bf  (1)}\quad $\displaystyle\frac{n}{4}< \beta< 1$ for $n=2,3$;
  \\[0.15cm]
  \indent Case {\bf  (2)}\quad $  \displaystyle\beta=\frac{n}{4}$ for $n=2,3$.\\
  \\
$\heartsuit$ We first consider  Case {\bf  (1)}\quad $\displaystyle\frac{n}{4}<\beta<1$ for $n=2,3$.
\\[0.2cm]
\indent In this case, a straightforward computation shows that
  \begin{equation}\label{4.50}
\left\{
\begin{array}{ll}
\displaystyle\dfrac{n}{\beta}=\dfrac{2n}{n-2\cdot\frac{n-2\beta}{2}},~~\dfrac{1}{2}=\dfrac{n-2}{2}<\dfrac{n-2\beta}{2}
  <\dfrac{3}{4}~~
  &\mbox{for}~~n=3,
  \\[0.4cm]
    \displaystyle 0=\dfrac{n-2}{2}<\dfrac{n-2\beta}{2}< \dfrac{1}{2}~~
  &\mbox{for}~~n=2,
  \\[0.3cm]
    \displaystyle B=\dfrac{n-2\beta}{2}+1-\beta=\dfrac{n}{2}+1-2\beta,&
   \\[0.3cm]
   \displaystyle\dfrac{n}{2}-1<B< 1,~~  &\mbox{for}~~n=2,3.
   \end{array}
\right.
\end{equation}
Note that \eqref{4.49}, one attains
 \begin{equation}\label{4.51}
\left.
\begin{array}{ll}
\displaystyle
    & \displaystyle I_{M,P}= \left|\left\langle \partial^{P}_{t}(\mathbf{u} \cdot\nabla \mathbf{v}) , \partial^{P}_{t}\Delta^{M} \mathbf{v} \right\rangle\right|
    \\[0.3cm]
     &\displaystyle\qquad=   \sum\limits^{P}_{p=0} \sum\limits^{M-1}_{m=0}\left(\begin{array}{l}
      P\\
      p
     \end{array}\right)\left(\begin{array}{l}
      M-1\\
      \quad m
     \end{array}\right)  \left|\left\langle  \Lambda^{1-\beta}\left(\partial^{p}_{t}\nabla^{m}\mathbf{u}\cdot \partial^{P-p}_{t}\nabla^{M-m}\mathbf{v}\right), \partial^{P}_{t} \nabla^{M }\Lambda^{\beta} \mathbf{v} \right\rangle\right|
     \\[0.3cm]
     &\displaystyle\qquad \leq  \sum\limits^{P}_{p=0} \sum\limits^{M-1}_{m=0}\left(\begin{array}{l}
      P\\
      p
     \end{array}\right)\left(\begin{array}{l}
      M-1\\
       \quad m
     \end{array}\right)\left\| \Lambda^{1-\beta}\left(\partial^{p}_{t}\nabla^{m}\mathbf{u}\cdot \partial^{P-p}_{t}\nabla^{M-m}\mathbf{v}\right)\right \|  _{L^{2}(\mathbb{R}^{n})}
     \\[0.3cm]
    &\displaystyle\qquad \qquad \qquad\qquad\qquad \qquad \qquad\qquad  \cdot \left\| \partial^{P}_{t} \nabla^{M }\Lambda^{\beta} \mathbf{v} \right\|  _{L^{2}(\mathbb{R}^{n})}.
          \end{array}
\right.
\end{equation}
Thanks to higher order fractional Leibniz's rule \cite{FGO},  $\displaystyle\left\| \Lambda^{1-\beta}\left(\partial^{p}_{t}\nabla^{m}\mathbf{u}\cdot \partial^{P-p}_{t}\nabla^{M-m}\mathbf{v}\right) \right\|  _{L^{2}(\mathbb{R}^{n})}$ can be bounded as follows:
\begin{equation}\label{4.52}
\left.
\begin{array}{ll}
&\displaystyle\left\| \Lambda^{1-\beta}\left(\partial^{p}_{t}\nabla^{m}\mathbf{u}\cdot \partial^{P-p}_{t}\nabla^{M-m}\mathbf{v}\right) \right\|  _{L^{2}(\mathbb{R}^{n})}
\\[0.3cm]
&\displaystyle\qquad \leq\left\| \Lambda^{1-\beta}\left(\partial^{p}_{t}\nabla^{m}\mathbf{u}\cdot \partial^{P-p}_{t}\nabla^{M-m}\mathbf{v}\right)\right.
\\[0.3cm]
&\displaystyle\qquad\quad-\partial^{p}_{t}\nabla^{m}\mathbf{u}\cdot \partial^{P-p}_{t}\nabla^{M-m}\Lambda^{1-\beta}\mathbf{v}
\\[0.3cm]
   &\displaystyle\qquad\quad\left.-\partial^{p}_{t}\Lambda^{1-\beta}\nabla^{m}\mathbf{u}\cdot \partial^{P-p}_{t}\nabla^{M-m}\mathbf{v}\right\|  _{L^{2}(\mathbb{R}^{n})}
   \\[0.3cm]
  &\displaystyle\qquad\quad +\left\|  \partial^{p}_{t}\nabla^{m}\mathbf{u}\cdot \partial^{P-p}_{t}\nabla^{M-m}\Lambda^{1-\beta}\mathbf{v} \right\|  _{L^{2}(\mathbb{R}^{n})}
  \\[0.3cm]
  &\displaystyle\qquad\quad +\left\|  \partial^{p}_{t}\Lambda^{1-\beta}\nabla^{m}\mathbf{u}\cdot \partial^{P-p}_{t}\nabla^{M-m}\mathbf{v}\right\|  _{L^{2}(\mathbb{R}^{n})}.
   \end{array}
\right.
\end{equation}
 Due to Lemma \ref{l2.10}, (I) of Lemma \ref{l2.11} and \eqref{4.50}, for  $\displaystyle\frac{1}{2}=\frac{1}{n/\beta}+\frac{1}{2n/(n-2\beta)}$ and $0<1-\beta<\beta<1$, the first term on the right hand side of \eqref{4.52} can be bounded by
 \begin{equation}\label{4.53}
\left.
\begin{array}{ll}
&\displaystyle \left\| \Lambda^{1-\beta}\left(\partial^{p}_{t}\nabla^{m}\mathbf{u}\cdot \partial^{P-p}_{t}\nabla^{M-m}\mathbf{v}\right)\right.
\\[0.3cm]
   &\displaystyle\qquad\quad-\partial^{p}_{t}\nabla^{m}\mathbf{u}\cdot \partial^{P-p}_{t}\nabla^{M-m}\Lambda^{1-\beta}\mathbf{v}
   \\[0.3cm]
   &\displaystyle\qquad\quad\left.-\partial^{p}_{t}\Lambda^{1-\beta}\nabla^{m}\mathbf{u}\cdot \partial^{P-p}_{t}\nabla^{M-m}\mathbf{v}\right\|  _{L^{2}(\mathbb{R}^{n})}
   \\[0.3cm]
 &\displaystyle\qquad\leq C \left\|  \partial^{p}_{t}\Lambda^{1-\beta} \nabla^{m}\mathbf{u} \right \|  _{L^{\frac{n}{\beta}}(\mathbb{R}^{n})}
  \left\|  \partial^{P-p}_{t}  \nabla^{M-m}  \mathbf{v} \right \|  _{L^{\frac{2n}{n-2\beta}}(\mathbb{R}^{n})}
  \\[0.3cm]
 &\displaystyle\qquad\leq C \left\|  \partial^{p}_{t}\Lambda^{\frac{n}{2}+1-2\beta} \nabla^{m}\mathbf{u}  \right\|  _{L^{2}(\mathbb{R}^{n})}
  \left\|  \partial^{P-p}_{t}  \nabla^{M-m} \Lambda^{ \beta}\mathbf{ v }\right \|  _{L^{2}(\mathbb{R}^{n})}
  \\[0.3cm]
  &\displaystyle\qquad\leq C \left\|  \partial^{p}_{t} \nabla^{m+1}\mathbf{u}  \right\|  _{L^{2}(\mathbb{R}^{n})}
  \left\|  \partial^{P-p}_{t}  \nabla^{M-m} \Lambda^{ \beta} \mathbf{v}  \right\|  _{L^{2}(\mathbb{R}^{n})}.
  \end{array}
\right.
\end{equation}
On the other hand, note that Lemma \ref{l2.8}, Lemma \ref{l2.9} and Lemma \ref{l2.11}, the second and the third terms on the right hand side of \eqref{4.52} can be bounded by
\begin{equation}\label{4.54}
\left.
\begin{array}{ll}
&\displaystyle \left\|  \partial^{p}_{t}\nabla^{m}\mathbf{u}\cdot \partial^{P-p}_{t}\nabla^{M-m}\Lambda^{1-\beta}\mathbf{v} \right\|  _{L^{2}(\mathbb{R}^{n})}
\\[0.3cm]
 &\displaystyle\qquad\leq C\left \|  \partial^{p}_{t}  \nabla^{m}\mathbf{u} \right \|  _{L^{\frac{n}{ 2\beta-1 }}(\mathbb{R}^{n})}
  \left\|  \partial^{P-p}_{t}  \nabla^{M-m}  \Lambda^{1-\beta}\mathbf{v}  \right\|  _{L^{\frac{2n}{n-2(2\beta-1)}}(\mathbb{R}^{n})}
  \\[0.3cm]
  &\displaystyle\qquad\leq C \left\|  \partial^{p}_{t}  \nabla^{m}\mathbf{u}  \right\|  _{L^{\frac{2n}{n-2\left(\frac{n}{2}+1-2\beta\right)}}(\mathbb{R}^{n})}
  \left\|  \partial^{P-p}_{t}  \nabla^{M-m}  \Lambda^{\beta}\mathbf{v}\right  \|  _{L^{2}(\mathbb{R}^{n})}
  \\[0.4cm]
  &\displaystyle\qquad\leq C \left\|  \partial^{p}_{t}  \nabla^{m}\Lambda^{\frac{n}{2}+1-2\beta}\mathbf{u } \right\|  _{L^{2}(\mathbb{R}^{n})}
 \left \|  \partial^{P-p}_{t}  \nabla^{M-m}  \Lambda^{\beta}\mathbf{v } \right\|  _{L^{2}(\mathbb{R}^{n})}
  \\[0.3cm]
  &\displaystyle\qquad\leq C \left\|  \partial^{p}_{t}  \nabla^{m+1} \mathbf{u }\right \|  _{L^{2}(\mathbb{R}^{n})}
  \left\|  \partial^{P-p}_{t}  \nabla^{M-m}  \Lambda^{\beta}\mathbf{v} \right \|  _{L^{2}(\mathbb{R}^{n})},
  \end{array}
\right.
\end{equation}
and
\begin{equation}\label{4.55}
\left.
\begin{array}{ll}
& \displaystyle\left\|  \partial^{p}_{t}\Lambda^{1-\beta}\nabla^{m}\mathbf{u}\cdot \partial^{P-p}_{t}\nabla^{M-m}\mathbf{v}\right \|  _{L^{2}(\mathbb{R}^{n})}
\\[0.3cm]
 &\displaystyle\qquad\leq C \left\|  \partial^{p}_{t} \Lambda^{1-\beta} \nabla^{m}\mathbf{u}  \right \|  _{L^{\frac{n}{\beta}}(\mathbb{R}^{n})}
  \left\|  \partial^{P-p}_{t}  \nabla^{M-m} \mathbf{ v} \right  \|  _{L^{\frac{2n}{n-2\beta}}(\mathbb{R}^{n})}
  \\[0.3cm]
  &\displaystyle\qquad\leq C \left\|  \partial^{p}_{t}  \nabla^{m+1}\mathbf{u}  \right \|  _{L^{2}(\mathbb{R}^{n})}
 \left \|  \partial^{P-p}_{t}  \nabla^{M-m} \Lambda^{ \beta} \mathbf{v} \right  \|  _{L^{2}(\mathbb{R}^{n})}.
  \end{array}
\right.
\end{equation}
Hence, combining \eqref{4.51} with \eqref{4.52}, \eqref{4.53}, \eqref{4.54} and \eqref{4.55} gives rise to
\begin{equation}\label{4.56}
\left.
\begin{array}{ll}
&\displaystyle I_{M,P}= \left|\left\langle \partial^{P}_{t}\left(\mathbf{u} \cdot\nabla \mathbf{v}\right) , \partial^{P}_{t}\Delta^{M} \mathbf{v} \right\rangle\right|
\\[0.3cm]
&\displaystyle\qquad\leq C\sum\limits^{P}_{p=0} \sum\limits^{M-1}_{m=0} \left\|  \partial^{p}_{t} \nabla^{m+1}\mathbf{u}  \right\|  _{L^{2}(\mathbb{R}^{n})} \left\|  \partial^{P-p}_{t}  \nabla^{M-m} \Lambda^{\beta} \mathbf{v}  \right\|  _{L^{2}(\mathbb{R}^{n})}
\\[0.3cm]
   &\displaystyle\qquad\qquad\qquad\qquad\qquad\qquad\cdot
   \left\|  \partial^{P}_{t}  \nabla^{M}\Lambda^{\beta}\mathbf{v}  \right\|  _{L^{2}(\mathbb{R}^{n})}
   \\[0.3cm]
   &\displaystyle\qquad\leq C\sum\limits^{P}_{p=0}  \left \|  \partial^{p}_{t}\mathbf{v}\right \|^{2}  _{L^{2}(\mathbb{R}^{n})} \left\|  \partial^{P-p}_{t}  \nabla^{M } \Lambda^{\beta}\mathbf{ v } \right\|^{2}  _{L^{2}(\mathbb{R}^{n})}
   \\[0.3cm]
   &\displaystyle\qquad\quad +C\sum\limits^{P}_{p=0}   \left\|  \partial^{p}_{t}\nabla \mathbf{v} \right \|^{2}  _{L^{2}(\mathbb{R}^{n})} \left\|  \partial^{P-p}_{t}  \nabla^{M -1} \Lambda^{\beta} \mathbf{v} \right \|^{2}  _{L^{2}(\mathbb{R}^{n})}
   \\[0.3cm]
    &\displaystyle\qquad\quad +C\sum\limits^{P}_{p=0}\sum\limits^{M-1}_{m=2}   \left\|  \partial^{p}_{t}\nabla^{m-1} \mathbf{v} \right\|^{2}  _{L^{2}(\mathbb{R}^{n})} \left\|  \partial^{P-p}_{t}  \nabla^{M -m} \Lambda^{\beta} \mathbf{v } \right\|^{2}  _{L^{2}(\mathbb{R}^{n})}
    \\[0.5cm]
    &\displaystyle\qquad\quad +\frac{\nu}{4} \left\|  \partial^{P }_{t}  \nabla^{M } \Lambda^{\beta} \mathbf{v } \right\|^{2}  _{L^{2}(\mathbb{R}^{n})}.
   \end{array}
\right.
\end{equation}
Similar estimates to that used in the estimates for $I_{M,P}$ are valid for $J_{M,P}$ in \eqref{4.49}:
\begin{equation}\label{4.57}
\left.
\begin{array}{ll}
& \displaystyle J_{M,P}= \left|\left\langle \partial^{P}_{t}\left(\mathbf{u} \cdot\nabla \mathbf{v}^{T}\right) , \partial^{P}_{t}\Delta^{M} \mathbf{v} \right\rangle\right|
\\[0.3cm]
  &\displaystyle\qquad\leq C\sum\limits^{P}_{p=0}  \left\|  \partial^{p}_{t}\mathbf{v}  \right\|^{2}  _{L^{2}(\mathbb{R}^{n})}\left \|  \partial^{P-p}_{t}  \nabla^{M } \Lambda^{\beta} \mathbf{v}  \right\|^{2}  _{L^{2}(\mathbb{R}^{n})}
  \\[0.3cm]
  &\displaystyle\qquad\quad+ C\sum\limits^{P}_{p=0}  \left\|  \partial^{p}_{t}\nabla \mathbf{v}  \right\|^{2}  _{L^{2}(\mathbb{R}^{n})}\left \|  \partial^{P-p}_{t}  \nabla^{M-1 } \Lambda^{\beta} \mathbf{v}  \right\|^{2}  _{L^{2}(\mathbb{R}^{n})}
  \\[0.3cm]
    &\displaystyle\qquad\quad +C\sum\limits^{P}_{p=0}\sum\limits^{M-1}_{m=2}   \left\|  \partial^{p}_{t}\nabla^{m-1} \mathbf{v} \right\|^{2}  _{L^{2}(\mathbb{R}^{n})}\left \|  \partial^{P-p}_{t}  \nabla^{M -m} \Lambda^{\beta} \mathbf{v}  \right\|^{2}  _{L^{2}(\mathbb{R}^{n})}
    \\[0.5cm]
    &\displaystyle\qquad\quad +\frac{\nu}{4} \left\|  \partial^{P }_{t}  \nabla^{M } \Lambda^{\beta} \mathbf{v }\right \|^{2}  _{L^{2}(\mathbb{R}^{n})}.
    \end{array}
\right.
\end{equation}
Substituting the above estimates \eqref{4.56} and \eqref{4.57} into \eqref{4.49} leads to the following estimate under the case $\displaystyle\frac{n}{4}< \beta< 1$ with $n=2,3$:
\begin{equation}\label{4.58}
\left.
\begin{array}{ll}
& \displaystyle\frac{d}{dt} \left\|  \partial^{P}_{t}\nabla^{M} \mathbf{v}  \right\|^{2}  _{L^{2}(\mathbb{R}^{n})}+
  \nu \left\|  \partial^{P }_{t}  \nabla^{M } \Lambda^{\beta} \mathbf{v}  \right\|^{2}  _{L^{2}(\mathbb{R}^{n})}
  \\[0.3cm]
   &\displaystyle\qquad\leq C\sum\limits^{P}_{p=0}  \left\|  \partial^{p}_{t}\mathbf{v} \right \|^{2}  _{L^{2}(\mathbb{R}^{n})}\left \|  \partial^{P-p}_{t}  \nabla^{M } \Lambda^{\beta} \mathbf{v} \right \|^{2}  _{L^{2}(\mathbb{R}^{n})}
   \\[0.3cm]
  &\displaystyle\qquad\quad+ C\sum\limits^{P}_{p=0} \left \|  \partial^{p}_{t}\nabla \mathbf{v}  \right\|^{2}  _{L^{2}(\mathbb{R}^{n})} \left\|  \partial^{P-p}_{t}  \nabla^{M-1 } \Lambda^{\beta} \mathbf{v} \right \|^{2}  _{L^{2}(\mathbb{R}^{n})}
  \\[0.3cm]
    &\displaystyle\qquad\quad +C\sum\limits^{P}_{p=0}\sum\limits^{M-1}_{m=2} \left  \|  \partial^{p}_{t}\nabla^{m-1}\mathbf{ v} \right\|^{2}  _{L^{2}(\mathbb{R}^{n})} \left\|  \partial^{P-p}_{t}  \nabla^{M -m} \Lambda^{\beta} \mathbf{v}  \right\|^{2}  _{L^{2}(\mathbb{R}^{n})}.
  \end{array}
\right.
\end{equation}
$\heartsuit$ We now tackle \eqref{4.49} under the Case  {\bf (2)}\quad $ \displaystyle \beta=\frac{n}{4}$ for $n=2,3$.
\\[0.3cm]
\indent In this case, it is easy to check that $\displaystyle 1-\frac{n}{4}\leq \frac{n}{4}$. Recall \eqref{4.49}, we shall estimate $I_{M,P}$ and $J_{M,P}$, respectively. We first handle $I_{M,P}$.\\
 \indent By a similar proof to that for case (1), one deduces the following:
\begin{equation}\label{4.59}
\left.
\begin{array}{ll}
\displaystyle
    & \displaystyle I_{M,P}= \left|\left\langle \partial^{P}_{t}(\mathbf{u} \cdot\nabla \mathbf{v}) , \partial^{P}_{t}\Delta^{M} \mathbf{v} \right\rangle\right|
    \\[0.3cm]
    &\displaystyle\qquad=   \sum\limits^{P}_{p=0} \sum\limits^{M-1}_{m=0}\left(\begin{array}{l}
      P\\
      p
     \end{array}\right)\left(\begin{array}{l}
      M-1\\
      \quad m
     \end{array}\right) \left \langle  \Lambda^{1-\frac{n}{4}}\left(\partial^{p}_{t}\nabla^{m}\mathbf{u}\cdot \partial^{P-p}_{t}\nabla^{M-m}\mathbf{v}\right), \partial^{P}_{t} \nabla^{M }\Lambda^{\frac{n}{4}}\mathbf{ v} \right\rangle
     \\[0.3cm]
     &\displaystyle\qquad \leq  \sum\limits^{P}_{p=0} \sum\limits^{M-1}_{m=0}\left(\begin{array}{l}
      P\\
      p
     \end{array}\right)\left(\begin{array}{l}
      M-1\\
       \quad m
     \end{array}\right)\left\| \Lambda^{1-\frac{n}{4}}\left(\partial^{p}_{t}\nabla^{m}\mathbf{u}\cdot \partial^{P-p}_{t}\nabla^{M-m}\mathbf{v}\right)\right \|  _{L^{2}(\mathbb{R}^{n})}
     \\[0.3cm]
    &\displaystyle\qquad \qquad \qquad\qquad\qquad \qquad \qquad\qquad  \cdot\left \| \partial^{P}_{t} \nabla^{M }\Lambda^{\frac{n}{4}} \mathbf{v} \right\|  _{L^{2}(\mathbb{R}^{n})}.
          \end{array}
\right.
\end{equation}
However, direct calculation gives
\begin{equation}\label{4.60}
\left.
\begin{array}{ll}
&\displaystyle\left\| \Lambda^{1-\frac{n}{4}}\left(\partial^{p}_{t}\nabla^{m}\mathbf{u}\cdot \partial^{P-p}_{t}\nabla^{M-m}\mathbf{v}\right) \right\|  _{L^{2}(\mathbb{R}^{n})}
\\[0.3cm]
 &\displaystyle\qquad \leq\left\| \Lambda^{1-\frac{n}{4}}\left(\partial^{p}_{t}\nabla^{m}\mathbf{u}\cdot \partial^{P-p}_{t}\nabla^{M-m}\mathbf{v}\right)\right.
 \\[0.3cm]
   &\displaystyle\qquad\quad-\partial^{p}_{t}\nabla^{m}\mathbf{u}\cdot \partial^{P-p}_{t}\nabla^{M-m}\Lambda^{1-\frac{n}{4}}\mathbf{v}
   \\[0.3cm]
   &\displaystyle\qquad\quad-\left.\partial^{p}_{t}\Lambda^{1-\frac{n}{4}}\nabla^{m}\mathbf{u}\cdot \partial^{P-p}_{t}\nabla^{M-m}\mathbf{v}\right\|  _{L^{2}(\mathbb{R}^{n})}
   \\[0.3cm]
  &\displaystyle\qquad\quad +\left\|  \partial^{p}_{t}\nabla^{m}\mathbf{u}\cdot \partial^{P-p}_{t}\nabla^{M-m}\Lambda^{1-\frac{n}{4}}\mathbf{v} \right\|  _{L^{2}(\mathbb{R}^{n})}
  \\[0.3cm]
  &\displaystyle\qquad\quad +\left\|  \partial^{p}_{t}\Lambda^{1-\frac{n}{4}}\nabla^{m}\mathbf{u}\cdot \partial^{P-p}_{t}\nabla^{M-m}\mathbf{v}\right\|  _{L^{2}(\mathbb{R}^{n})}.
 \end{array}
\right.
\end{equation}
Due to Lemma \ref{l2.9} and Lemma \ref{l2.11}, one deduces that for $\displaystyle 0<\beta_{1}<1-\frac{n}{4}$, the first term on the right hand side of \eqref{4.60} can be estimated by
\begin{equation}\label{4.61}
\left.
\begin{array}{ll}
& \displaystyle\left\| \Lambda^{1-\frac{n}{4}}\left(\partial^{p}_{t}\nabla^{m}\mathbf{u}\cdot \partial^{P-p}_{t}\nabla^{M-m}\mathbf{v}\right)\right.
\\[0.3cm]
 &\displaystyle\qquad\quad-\partial^{p}_{t}\nabla^{m}\mathbf{u}\cdot \partial^{P-p}_{t}\nabla^{M-m}\Lambda^{1-\frac{n}{4}}\mathbf{v}
 \\[0.3cm]
  &\displaystyle\qquad\quad\left.-\partial^{p}_{t}\Lambda^{1-\frac{n}{4}}\nabla^{m}\mathbf{u}\cdot \partial^{P-p}_{t}\nabla^{M-m}\mathbf{v}\right\|  _{L^{2}(\mathbb{R}^{n})}
  \\[0.3cm]
  &\displaystyle\qquad\leq C \left\|  \partial^{p}_{t}\Lambda^{1- \frac{n}{4}-\beta_{1}} \nabla^{m}\mathbf{u}  \right\|  _{L^{\frac{2n}{n-2(n/4+\beta_{1})}}(\mathbb{R}^{n})}
  \left\|  \partial^{P-p}_{t}  \nabla^{M-m} \Lambda^{ \beta_{1}} \mathbf{v} \right \|  _{L^{\frac{2n}{n-2(n/4-\beta_{1}) }}(\mathbb{R}^{n})}
  \\[0.3cm]
  &\displaystyle\qquad\leq C \left\|  \partial^{p}_{t}  \nabla^{m-1}\mathbf{v}  \right\|  _{L^{2}(\mathbb{R}^{n})}
  \left\|  \partial^{P-p}_{t}  \nabla^{M-m} \Lambda^{\frac{n}{4}} \mathbf{v}  \right\|  _{L^{2}(\mathbb{R}^{n})},
  \end{array}
\right.
\end{equation}
where $\displaystyle\beta_{1}\in \left(0,\frac{1}{2}\right)$, $\displaystyle 1-\frac{n}{4}-\beta_{1}\in \left(0,\frac{1}{2}\right)$,
  $\displaystyle\frac{1}{2}=\frac{1}{p_{1}}+\frac{1}{p_{2}}$ with $p_{1},p_{2}\in (1,\infty)$, $\displaystyle p_{1}=\frac{2n}{n-2(n/4+\beta_{1})}$,   $\displaystyle p_{2}=\frac{2n}{n-2(n/4-\beta_{1})}$.\\
   \indent Let us turn to estimate the second and the third terms on the right hand side of \eqref{4.60}. Thanks to Agmon's inequality and Lemma \ref{l2.12} and the assumption of (VI) for $\beta=\dfrac{n}{4}$, we have
 \begin{equation}\label{4.62}
\left.
\begin{array}{ll}
  &\displaystyle \left\|  \partial^{p}_{t}\nabla^{m}\mathbf{u}\cdot \partial^{P-p}_{t}\nabla^{M-m}\Lambda^{1-\frac{n}{4}}\mathbf{v}\right\|  _{L^{2}(\mathbb{R}^{n})}
  \\[0.3cm]
  &\displaystyle\qquad\quad+
   \left\| \partial^{p}_{t}\Lambda^{1-\frac{n}{4}}\nabla^{m}\mathbf{u}\cdot \partial^{P-p}_{t}\nabla^{M-m}\mathbf{v}\right\|  _{L^{2}(\mathbb{R}^{n})}
   \\[0.3cm]
  &\displaystyle\qquad\leq C \left\|  \partial^{p}_{t}  \nabla^{m}\mathbf{u}  \right\|  _{L^{\infty}(\mathbb{R}^{n})}
  \left\|  \partial^{P-p}_{t}  \nabla^{M-m} \Lambda^{ 1-\frac{n}{4}} \mathbf{v} \right \|  _{L^{2}(\mathbb{R}^{n})}
  \\[0.3cm]
  &\displaystyle\qquad\quad +\left\|  \partial^{p}_{t} \Lambda^{ 1-\frac{n}{4}} \nabla^{m}\mathbf{u} \right \|  _{L^{4}(\mathbb{R}^{n})}
  \left\|  \partial^{P-p}_{t}  \nabla^{M-m}  \mathbf{v}  \right\|  _{L^{4}(\mathbb{R}^{n})}
  \\[0.3cm]
  &\displaystyle\qquad\leq C \left\|  \partial^{p}_{t}  \nabla^{m}\mathbf{u}  \right\| ^{\frac{1}{2}} _{H^{1}(\mathbb{R}^{n})}\left\|  \partial^{p}_{t}  \nabla^{m}\mathbf{u} \right \| ^{\frac{1}{2}} _{H^{2}(\mathbb{R}^{n})}
  \left\|  \partial^{P-p}_{t}  \nabla^{M-m} \Lambda^{\frac{n}{4}} \mathbf{v} \right \|  _{L^{2}(\mathbb{R}^{n})}
  \\[0.3cm]
  &\displaystyle\qquad\quad +C \left( \left\|  \partial^{p}_{t}\Lambda^{ 1-\frac{n}{4}}  \nabla^{m}\mathbf{u} \right \|  _{L^{2}(\mathbb{R}^{n})}+\left\|  \partial^{p}_{t}\Lambda^{ 1-\frac{n}{4}+\frac{n}{4}}  \nabla^{m}\mathbf{u}  \right\|  _{L^{2}(\mathbb{R}^{n})}\right)
  \\[0.3cm]
  &\displaystyle\qquad\qquad \times  \left( \left\|  \partial^{P-p}_{t}  \nabla^{M-m}\mathbf{v} \right \|  _{L^{2}(\mathbb{R}^{n})}+\left\|  \partial^{P-p}_{t}\Lambda^{ \frac{n}{4}}  \nabla^{M-m}\mathbf{v}\right \|  _{L^{2}(\mathbb{R}^{n})}\right)
  \\[0.3cm]
  &\displaystyle\qquad\leq C\left\|  \partial^{p}_{t}  \nabla^{m}\mathbf{v}  \right\|  _{L^{2}(\mathbb{R}^{n})}\left\|  \partial^{P-p}_{t}  \nabla^{M-m}\Lambda^{ \frac{n}{4}} \mathbf{v}  \right\|  _{L^{2}(\mathbb{R}^{n})}.
  \end{array}
\right.
\end{equation}
 This together with \eqref{4.59}, \eqref{4.60} and \eqref{4.61} gives rise to
 \begin{equation}\label{4.63}
\left.
\begin{array}{ll}
  &\displaystyle I_{M,P}= \left|\left\langle \partial^{P}_{t}(\mathbf{u} \cdot\nabla \mathbf{v}) , \partial^{P}_{t}\Delta^{M} \mathbf{v} \right\rangle\right|
  \\[0.3cm]
  &\displaystyle\quad\quad\leq C\sum\limits^{P}_{p=0} \sum\limits^{M-1}_{m=0}\left\|  \partial^{p}_{t}\nabla^{m} \mathbf{v}\right\|  _{L^{2}(\mathbb{R}^{n})}
   \left\| \partial^{P-p}_{t}\Lambda^{\frac{n}{4}}\nabla^{M-m} \mathbf{v}\right\|  _{L^{2}(\mathbb{R}^{n})}
   \left\| \partial^{P-p }_{t}\Lambda^{\frac{n}{4}}\nabla^{M } \mathbf{v}\right\|  _{L^{2}(\mathbb{R}^{n})}
   \\[0.3cm]
  &\displaystyle\quad\quad\leq C \sum\limits^{P}_{p=0} \left\|  \partial^{p}_{t}  \mathbf{v}\right\| ^{2} _{L^{2}(\mathbb{R}^{n})}
   \left\| \partial^{P-p}_{t}\Lambda^{\frac{n}{4}}\nabla^{M } \mathbf{v}\right\|^{2}  _{L^{2}(\mathbb{R}^{n})}
   \\[0.3cm]
  &\displaystyle\quad\quad\quad+ C \sum\limits^{P}_{p=0} \left\|  \partial^{p}_{t}  \nabla \mathbf{v}\right\| ^{2} _{L^{2}(\mathbb{R}^{n})}
  \left \| \partial^{P-p}_{t}\Lambda^{\frac{n}{4}}\nabla^{M-1 } \mathbf{v}\right\|^{2}  _{L^{2}(\mathbb{R}^{n})}
   \\[0.3cm]
  &\displaystyle\quad\quad\quad+ C \sum\limits^{P}_{p=0} \sum\limits^{M-1}_{m=2} \left\|  \partial^{p}_{t}  \nabla^{m} \mathbf{v}\right\| ^{2} _{L^{2}(\mathbb{R}^{n})}
   \left\| \partial^{P-p}_{t}\Lambda^{\frac{n}{4}}\nabla^{M-m } \mathbf{v}\right\|^{2}  _{L^{2}(\mathbb{R}^{n})}
   \\[0.5cm]
  &\displaystyle\quad\quad\quad+ \frac{\nu}{4}\left\| \partial^{P }_{t}\Lambda^{\frac{n}{4}}\nabla^{M } \mathbf{v}\right\|^{2}  _{L^{2}(\mathbb{R}^{n})}.
  \end{array}
\right.
\end{equation}
 Here we have used Cauchy-Schwarz inequality in the third inequality.\\
 \indent In the same manner, one may deduce the following estimate for $J_{M,P}$ in \eqref{4.49}:
 \begin{equation}\label{4.64}
\left.
\begin{array}{ll}
  &\displaystyle
 J_{M,P}= \left|\left\langle \partial^{P}_{t}\left(\mathbf{u} \cdot\nabla \mathbf{v}^{T}\right) , \partial^{P}_{t}\Delta^{M} \mathbf{v} \right\rangle\right|
 \\[0.3cm]
  &\displaystyle\quad\quad\leq C \sum\limits^{P}_{p=0} \left\|  \partial^{p}_{t}  \mathbf{v}\right\| ^{2} _{L^{2}(\mathbb{R}^{n})}
   \left\| \partial^{P-p}_{t}\Lambda^{\frac{n}{4}}\nabla^{M } \mathbf{v}\right\|^{2}  _{L^{2}(\mathbb{R}^{n})}
   \\[0.3cm]
  &\displaystyle\quad\quad\quad+ C \sum\limits^{P}_{p=0} \left\|  \partial^{p}_{t}  \nabla \mathbf{v}\right\| ^{2} _{L^{2}(\mathbb{R}^{n})}
   \left\| \partial^{P-p}_{t}\Lambda^{\frac{n}{4}}\nabla^{M-1 }\mathbf{ v}\right\|^{2}  _{L^{2}(\mathbb{R}^{n})}
   \\[0.3cm]
  &\displaystyle\quad\quad\quad+ C \sum\limits^{P}_{p=0} \sum\limits^{M-1}_{m=2} \left\|  \partial^{p}_{t}  \nabla^{m} \mathbf{v}\right\| ^{2} _{L^{2}(\mathbb{R}^{n})}
   \left\| \partial^{P-p}_{t}\Lambda^{\frac{n}{4}}\nabla^{M-m } \mathbf{v}\right\|^{2}  _{L^{2}(\mathbb{R}^{n})}
   \\[0.5cm]
  &\displaystyle\quad\quad\quad+ \frac{\nu}{4}\left\| \partial^{P }_{t}\Lambda^{\frac{n}{4}}\nabla^{M } \mathbf{v}\right\|^{2}  _{L^{2}(\mathbb{R}^{n})}.
  \end{array}
\right.
\end{equation}
 Therefore, under the case of $\displaystyle\beta=\frac{n}{4}$ with $n=2,3$, substituting \eqref{4.63} and \eqref{4.64} into \eqref{4.49} yields that
 \begin{equation}\label{4.65}
\left.
\begin{array}{ll}
& \displaystyle\frac{d}{dt} \left\| \partial^{P}_{t}\nabla^{M } \mathbf{v}  \right\| ^{2}_{L^{2}(\mathbb{R}^{n})}+\nu \left\| \partial^{P}_{t}\Lambda^{\frac{n}{4}} \nabla^{M } \mathbf{v} \right \| ^{2}_{L^{2}(\mathbb{R}^{n})}
\\[0.3cm]
  &\displaystyle\qquad\leq C \sum\limits^{P}_{p=0} \left \|  \partial^{p}_{t}  \mathbf{v}\right\| ^{2} _{L^{2}(\mathbb{R}^{n})}
    \left\| \partial^{P-p}_{t}\Lambda^{\frac{n}{4}}\nabla^{M } \mathbf{v}\right\|^{2}  _{L^{2}(\mathbb{R}^{n})}
   \\[0.3cm]
  &\displaystyle\quad\quad\quad+ C \sum\limits^{P}_{p=0}  \left\|  \partial^{p}_{t}  \nabla \mathbf{v}\right\| ^{2} _{L^{2}(\mathbb{R}^{n})}
    \left\| \partial^{P-p}_{t}\Lambda^{\frac{n}{4}}\nabla^{M-1 } \mathbf{v}\right\|^{2}  _{L^{2}(\mathbb{R}^{n})}
   \\[0.3cm]
  &\displaystyle\quad\quad\quad+ C \sum\limits^{P}_{p=0} \sum\limits^{M-1}_{m=2}  \left\|  \partial^{p}_{t}  \nabla^{m}\mathbf{ v}\right\| ^{2} _{L^{2}(\mathbb{R}^{n})}
    \left\| \partial^{P-p}_{t}\Lambda^{\frac{n}{4}}\nabla^{M-m } \mathbf{v}\right\|^{2}  _{L^{2}(\mathbb{R}^{n})}.
  \end{array}
\right.
\end{equation}
 With \eqref{4.49}, \eqref{4.58} and \eqref{4.65}, for $\displaystyle\frac{n}{4}\leq\beta< 1$ with $n=2,3$, we always have
\begin{equation}\label{4.66}
\left.
\begin{array}{ll}
  &\displaystyle \frac{d}{dt}\left\| \partial^{P}_{t}\nabla^{M } \mathbf{v}  \right\| ^{2}_{L^{2}(\mathbb{R}^{n})}+\nu\left\| \partial^{P}_{t}\Lambda^{\beta} \nabla^{M } \mathbf{v}  \right\| ^{2}_{L^{2}(\mathbb{R}^{n})}
  \\[0.3cm]
  &\displaystyle\qquad\leq C \sum\limits^{P}_{p=0} \left\|  \partial^{p}_{t}  \mathbf{v}\right\| ^{2} _{L^{2}(\mathbb{R}^{n})}
  \left \| \partial^{P-p}_{t}\Lambda^{\beta}\nabla^{M } \mathbf{v}\right\|^{2}  _{L^{2}(\mathbb{R}^{n})}
   \\[0.3cm]
  &\displaystyle\quad\quad\quad+ C \sum\limits^{P}_{p=0} \left\|  \partial^{p}_{t}  \nabla \mathbf{v}\right\| ^{2} _{L^{2}(\mathbb{R}^{n})}
   \left\| \partial^{P-p}_{t}\Lambda^{\beta}\nabla^{M-1 } \mathbf{v}\right\|^{2}  _{L^{2}(\mathbb{R}^{n})}
   \\[0.3cm]
  &\displaystyle\quad\quad\quad+ C \sum\limits^{P}_{p=0} \sum\limits^{M}_{m=2} \left\|  \partial^{p}_{t}  \nabla^{m} \mathbf{v}\right\| ^{2} _{L^{2}(\mathbb{R}^{n})}
   \left\| \partial^{P-p}_{t}\Lambda^{\beta}\nabla^{M-m } \mathbf{v}\right\|^{2}  _{L^{2}(\mathbb{R}^{n})}.
  \end{array}
\right.
\end{equation}
Note that the inductive assumption \eqref{4.46}, one deduces that for $p=P$ and $m=0$
\begin{equation*}
\left.
\begin{array}{ll}
  &\displaystyle\frac{d}{dt}\left\| \partial^{P}_{t} \mathbf{v} \right \| ^{2}_{L^{2}(\mathbb{R}^{n})}+\nu\left\| \partial^{P}_{t}\Lambda^{\beta}   \mathbf{v} \right \| ^{2}_{L^{2}(\mathbb{R}^{n})}
  \\[0.3cm]
  &\displaystyle\qquad\leq C (1+t)^{-\frac{n}{2\beta}}\left\|  \partial ^{p}_{t}\Lambda^{\beta}  \mathbf{v}\right\| ^{2} _{L^{2}(\mathbb{R}^{n})}+ C (1+t)^{-2p-\frac{n}{2\beta}}
      (1+t)^{-2(P-p) -\frac{\beta}{\beta}-\frac{n}{2\beta}}
      \\[0.3cm]
  &\displaystyle\qquad\leq C (1+t)^{-\frac{n}{2\beta}}\left\| \partial ^{p}_{t}  \Lambda^{\beta}  \mathbf{v}\right\| ^{2} _{L^{2}(\mathbb{R}^{n})}+ C (1+t)^{-2P-1-\frac{n}{\beta}}.
\end{array}
\right.
\end{equation*}
Takeing $t$ large enough such that $\displaystyle C(1+t)^{- \frac{n}{ 2 \beta} }\leq \frac{\nu}{2}$, the above inequality then implies that
\begin{equation}\label{4.67}\frac{d}{dt}\left\| \partial^{P}_{t} \mathbf{v}  \right\| ^{2}_{L^{2}(\mathbb{R}^{n})}+\frac{\nu}{2}\left\| \partial^{P}_{t}\Lambda^{\beta} \mathbf{v} \right \| ^{2}_{L^{2}(\mathbb{R}^{n})}\leq C(1+t)^{-2P-1 -\frac{ n}{ \beta}}.\end{equation}
This together with (IV) of this theorem ensures that \eqref{4.46} and \eqref{4.48} hold for $p=P$ and $m=0$.\\
\indent So far we have shown that for $m=0$, and $\forall P\leq \frac{K}{2}$, there holds
\begin{equation*} \frac{d}{dt}\left\| \partial^{P}_{t} \mathbf{v} \right \| ^{2}_{L^{2}(\mathbb{R}^{n})}+\frac{\nu}{2}\left\| \partial^{P}_{t}\Lambda^{\beta} \mathbf{v} \right  \| ^{2}_{L^{2}(\mathbb{R}^{n})}\leq C(1+t)^{-2P-1 -\frac{ n}{ \beta}}.\end{equation*}
This deduces by Gronwall's inequality that
\begin{equation}\label{4.68}\frac{d}{dt}\left\| \partial^{P}_{t} \mathbf{v}  \right\| ^{2}_{L^{2}(\mathbb{R}^{n})} \leq C(1+t)^{-2P-\frac{ n}{ \beta}}.\end{equation}
{\bf Step 3}\quad We show that the decay rate \eqref{4.46} holds for any $m\leq M+1$ for $p<P$, and $m<M$ for $p=P$. That is,
\begin{equation}\label{4.69}\left\| \partial^{p}_{t}\nabla^{m} \mathbf{v}  \right \| ^{2}_{L^{2}(\mathbb{R}^{n})}\leq C(1+t)^{-2p-\frac{m}{\beta}-\frac{n}{2\beta}}.\end{equation}
The base case is \eqref{4.68} where \eqref{4.69} holds for $p=P$ and $m=0$. In the following, based on the inductive assumption \eqref{4.69}, we will show that the decay rate \eqref{4.69} holds for $m=M$ and $p=P$.\\
\indent Recall (I) and (II) of this theorem, applying the inductive assumption \eqref{4.69} to \eqref{4.66}, one deduces that
\begin{equation}\label{4.70}
\left.
\begin{array}{ll}
  &\displaystyle\frac{d}{dt}\left \| \partial^{P}_{t}\nabla^{M} \mathbf{v}  \right \| ^{2}_{L^{2}(\mathbb{R}^{n})}+\nu\left \| \partial^{P}_{t}\Lambda^{\beta}  \nabla^{M}\mathbf{ v } \right\| ^{2}_{L^{2}(\mathbb{R}^{n})}
   \\[0.3cm]
  &\displaystyle\qquad\leq C (1+t)^{-\frac{n}{2\beta}}\left \|  \partial^{P}_{t}  \Lambda^{\beta}\nabla^{M}  \mathbf{v}\right\| ^{2} _{L^{2}(\mathbb{R}^{n})}
   \\[0.3cm]
  &\displaystyle\qquad\quad+ C (1+t)^{ -\frac{1}{ \beta}-\frac{n}{2\beta}}\left \|  \partial^{P}_{t}  \Lambda^{\beta}\nabla^{M-1}  \mathbf{v}\right\| ^{2} _{L^{2}(\mathbb{R}^{n})}
   \\[0.3cm]
  &\displaystyle\qquad\quad+ C (1+t)^{-2P-\frac{M}{ \beta}-\frac{n}{ \beta}-1}.
       \end{array}
\right.
\end{equation}
Taking $t$ large enough such that $\displaystyle C(1+t)^{ -\frac{n}{ 2\beta}}\leq \frac{\nu}{2}$, thanks to \eqref{4.48}, using (IV) once again deduces that
\begin{equation}\label{4.71}\left\| \partial^{P}_{t}\nabla^{M} \mathbf{v} \right \| ^{2}_{L^{2}(\mathbb{R}^{n})} \leq C(1+t)^{-2P  -\frac{ M}{ \beta}-\frac{ n}{\beta}}\leq C(1+t)^{-2P  -\frac{ M}{ \beta}-\frac{ n}{2\beta}}.\end{equation}
This implies that the inductive assumption \eqref{4.69} holds for $m=M$ and $p=P$. By another bootstrap argument, we obtain for all $m+2p\beta\leq K$, the following optimal decay holds:
\begin{equation*}\label{4.71}\left\| \partial^{p}_{t}\nabla^{m} \mathbf{v } \right\| ^{2}_{L^{2}(\mathbb{R}^{n})} \leq C(1+t)^{-2p  -\frac{m}{ \beta}-\frac{ n}{2\beta}}.\end{equation*}
This completes the proof of (VI-2).\\
\indent So far, we finish the proof of Theorem 4.1. \hfill$\Box$

\section{Convergence to the NSE with nonlocal viscosity}
\allowdisplaybreaks
We observe that for $\alpha=0$, the system \eqref{1.1} formally reduces to the incompressible Navier-Stokes equations with fractional Laplacian viscosity
\begin{equation}\label{5.1}
\left\{
\begin{array}{ll}
 \displaystyle \mathbf{v}_{t}+\mathbf{v}\cdot \nabla \mathbf{v} +\nabla p=-\nu (-\Delta)^{\beta}\mathbf{v}
 \\[0.3cm]
 \displaystyle\mbox{div}~\mathbf{v}=0.
\end{array}
\right.
\end{equation}
By means of the fractional heat kernel estimates \cite{Dong-Li-1} and Leray projection, we figure out the relation between the nonlocal system \eqref{1.1} and \eqref{5.1}. In particular, we investigate the convergence of the solution of \eqref{1.1} as the filter parameter $\alpha\rightarrow 0$ to a solution of \eqref{5.1}, and relate the limit to \eqref{5.1}. To achieve this, we need first exploring  how a solution $\mathbf{u}$ of the Helmholtz equation
\begin{equation}\label{5.2}\mathbf{u}-\alpha^{2}\Delta \mathbf{u}=\mathbf{v}\end{equation}
 approaches $\mathbf{v}$ as $\alpha$ tends to zero. In \cite{Foias-1,Foias-2}, the authors clarified that how the solutions of the Camassa-Holm equations \eqref{1.3} approach solutions of the corresponding imcompressible Navier-Stokes equations \eqref{1.4} weakly when the filter parameter $\alpha$ tends to zero. In \cite{Clayton-08}, the authors established how solutions to the viscous Camassa-Holm equations \eqref{1.3} approach solutions to \eqref{1.4} strongly as $\alpha\rightarrow 0$ when the solution to \eqref{1.4} is known to be regular enough. Here, we expect to establish a similar result for the nonlocal Camassa-Holm equations \eqref{1.1} to that for \eqref{1.3}  mentioned as above . Precisely, we hope to make sure how solutions to \eqref{1.1} approach solutions to \eqref{5.1} strongly as $\alpha\rightarrow 0$ when the solutions to \eqref{5.1} are to be sufficiently regular. To attain this goal, we must establish some a priori estimates on the solutions of \eqref{1.1} which are independent of $\alpha$, but on regions of time where a solution to the nonlocal Navier-Stokes equations \eqref{5.1} is known to be regular by the functional analytic argument.
 \\[0.3cm]
 \indent The object of this section is to prove the following convergence theorem for \eqref{1.1}:
 \begin{thm}\label{t5.1}\rm
  For $\frac{n}{4}\leq \beta   <1$, $n=2,3$, let $\{\alpha_{i}\}$ be a sequence of filter coefficients tending to zero, and let $\mathbf{v}_{\alpha_{i}}$ be the solutions of \eqref{1.1}  constructed in Proposition \ref{p2.4} corresponding to the initial data
  $\mathbf{w}_{0}\in D_{\sigma}\left(\Lambda^{\beta} \right)(\mathbb{R}^{n})$ for $\dfrac{n}{4}< \beta < 1$, and $\mathbf{w}_{0}\in H_{\sigma}^{\beta}(\mathbb{R}^{n})$ for $\beta=\dfrac{n}{4}$. Let $\mathbf{w}$ be the solution of \eqref{5.1}
with the same initial data $\mathbf{w}_{0}$. In any time interval $[0,T]$, where a solution to \eqref{5.1} is known to be sufficiently regular, if there exists a bound
$$  \sup \limits_{\alpha_{i}} \sup\limits_{t\in [0,T] }\left(\left\|\mathbf{v}_{\alpha_{i}}\right\|_{L^l(\mathbb{R}^{n})} +\left\|\Lambda^{\beta}\mathbf{v}_{\alpha_{i}}\right\|_{L^l(\mathbb{R}^{n})}\right)\leq C,$$
 which is independent of $\alpha$, then $\mathbf{v}_{\alpha }~\mbox{approaches}~\mathbf{w}~\mbox{strongly~ in }~L^{\infty}\left([0,T],L^q(\mathbb{R}^{n})\right)~\mbox{as}~\alpha\rightarrow 0$,
 where $\displaystyle q=\frac{2s}{s-2}$, $\displaystyle  s=\frac{ln}{n-l\beta}$ and $\displaystyle  l> \frac{n}{ 3\beta-1}$.
 \end{thm}
 \indent Before proving this theorem, we first make some preliminary remarks and preparations.
 \begin{rem}\label{r5.2} \rm
 By a similar proof to that for the Camassa-Holm equations without any fractional viscosity term \eqref{1.3}, we deduce that a solution $\mathbf{u}$ of \eqref{5.2} approaches $\mathbf{v}$ weakly as the filter parameter $\alpha$ tends to zero. That is, fix $\mathbf{v}\in L^p(\mathbb{R}^{n})$, let $\{\alpha_{i}\}$ be a sequence of filter coefficients tending to zero, for each $\alpha_{i}$ there is a weak solution $\mathbf{u}_{\alpha_{i}}\in W^{1,p}(\mathbb{R}^{n})$ of \eqref{5.2} such that
$$\mathbf{u}_{\alpha_{i}}\rightharpoonup \mathbf{v}~\mbox{weakly~in}~L^p(\mathbb{R}^{n})~~\mbox{as}~\alpha_{i}\rightarrow 0.\eqno\Box$$
\end{rem}
  \indent Due to Remark \ref{r5.2}, we claim a stronger result if $\mathbf{v}$ is sufficiently differentiable.
  \begin{prop} \label{p5.3}\rm
   For $\displaystyle\frac{n}{4}\leq \beta   <1$, $n=2,3$, let $\mathbf{v}\in W^{\beta,p}(\mathbb{R}^{n})$ and $\mathbf{u}$ be the solution of \eqref{5.2}. Then for $ \alpha \in (0,1)$, there holds
 $$\left\|\mathbf{u}-\mathbf{v}\right\|_{L^q(\mathbb{R}^{n})}\leq C(n,p,q)\alpha^{\frac{\beta}{2}-\gamma}\left\|\Lambda^{\beta}\mathbf{v}\right\|_{L^p(\mathbb{R}^{n})}~\mbox{for}~
 \gamma=\frac{n}{2}\left(\frac{1}{p}-\frac{1}{q}\right)<\frac{\beta}{2}.$$
 In particular, if $\{\alpha_{i}\}$ is a sequence tending to zero, and $\mathbf{u}_{\alpha_{i}}$ are solutions of \eqref{5.2}, then for $\displaystyle\frac{1}{p}-\frac{1}{q}<\frac{\beta}{n}$,
 $$\mathbf{u}_{\alpha_{i}}\rightarrow \mathbf{v}~\mbox{strongly~in}~L^q(\mathbb{R}^{n})~\mbox{as}~\alpha_{i}\rightarrow 0.$$
  Here, $W^{\beta,p}(\mathbb{R}^{n})$ is defined by Definition \ref{d1.2}.
  \end{prop}
 {\bf Proof.}\quad If $\mathbf{u}$ and $\mathbf{v}$ satisfy \eqref{5.2}, then a preliminary calculation gives rise to
 \begin{equation}\label{5.3}
 \|\mathbf{u}-\mathbf{v}\|_{L^q(\mathbb{R}^{n})}\leq \alpha^{2}\|\Delta \mathbf{u}\|_{L^q(\mathbb{R}^{n})}.\end{equation}
 Since \eqref{5.2} is linear, the derivatives of the functions obey the relation
 \begin{equation}\label{5.4}
 \Lambda^{\beta}\mathbf{u}-\alpha^{2}\Lambda^{\beta}\Delta \mathbf{u}=\Lambda^{\beta}\mathbf{v}.\end{equation}
Applying Lemma \ref{l2.6} to \eqref{5.4} with $\displaystyle\gamma=\gamma_{3}=\frac{n}{2}
\left(\frac{1}{p}-\frac{1}{q}\right)<\frac{\beta}{2}$, we arrive at the following:
$$\|\Delta \mathbf{u}\|_{L^q(\mathbb{R}^{n})}\leq\left\|\Lambda^{\beta}\Delta \mathbf{u}\right\|_{L^p(\mathbb{R}^{n})}
\leq\dfrac{C(n,p,q)}{\alpha^{2-\frac{\beta}{2}+\gamma}}
 \left\|\Lambda^{\beta} \mathbf{v}\right\|_{L^p(\mathbb{R}^{n})}.$$
 This together with \eqref{5.3} yields that
 \begin{equation}\label{5.5}\| \mathbf{u}-\mathbf{v}\|_{L^q(\mathbb{R}^{n})}\leq \alpha^{ \frac{\beta}{2}-\gamma}
 \left\|\Lambda^{\beta} \mathbf{v}\right\|_{L^p(\mathbb{R}^{n})}.\end{equation}
 Note that  $\displaystyle\gamma<\frac{\beta}{2}$, replacing $\alpha$ with $\alpha_{i}$ in \eqref{5.5} and then letting $\alpha_{i}\rightarrow 0$, we immediately deduce the second statement that $u_{\alpha_{i} }\rightarrow v$ strongly in $L^q(\mathbb{R}^{n})$ for $\dfrac{1}{p}-\dfrac{1}{q} <\dfrac{\beta}{n}$.\hfill$\Box$
 \\[0.3cm]
 \indent We now mention some estimates concerning the fundamental solution of the linear nonlocal operator $\displaystyle\partial_{t}+(-\Delta)^{\frac{\gamma_{0}}{2}}$in \cite{Dong-Li-1}, which are key to the proof of Theorem \ref{t5.1}.
 \begin{lem}[ \cite{Dong-Li-1}]\label{l5.4}\rm
  Let $\gamma_{0}\in (1,2]$. Define $G_{\gamma_{0}}(t,x)$ by its Fourier transform $ \widehat{ G}_{\gamma_{0}}  (t,\xi)=e^{-t|\xi|^{\gamma_{0}}}$ for $t>0$. Then $G_{\gamma_{0}}(t,x)$ is the fundamental solution of the linear operator: $\partial_{t}+(-\Delta)^{\frac{\gamma_{0}}{2}}$. In addition, it enjoys the scaling property:
 $$G_{\gamma_{0}}(t,x)
 =t^{-\frac{n}{\gamma_{0}}}G_{\gamma_{0}}\left(1,t^{-\frac{1}{\gamma_{0}}}x\right).\eqno \Box$$
  \end{lem}
\begin{lem}[ \cite{Dong-Li-1}]\label{l5.5}\rm
   For $\gamma_{0}\in (0,\infty)$ and $p\in [1,\infty]$, let $k\geq 0$ be an integer and $\varepsilon\in (0,1]$. Then for some constant $C=C(n,\gamma_{0},\varepsilon)$, there holds that
  $$\left\|D^{k}_{x}\Lambda^{\alpha}G_{\gamma_{0}} (t,\cdot)\right\| _{L^{p}_{x}(\mathbb{R}^{n})}\leq C^{k+1}k^{\frac{k}{\gamma_{0}}}t^{-\frac{k+\alpha}{\gamma_{0}}-\frac{n}{\gamma_{0}}
  \left(1-\frac{1}{p}\right)}$$
 for any $\alpha$ satisfying
 \begin{equation*}
\left\{
\begin{array}{ll}
\varepsilon-1\leq\alpha\leq 1~~&\mbox{if}~~k\geq 1
\\[0.2cm]
\varepsilon \leq\alpha\leq 1~~\mbox{or}~~\alpha=0~~&\mbox{if}~~k=0.
 \end{array}
\right.
\end{equation*}
Here the constant $C$ can be taken to be independent of $p$.\hfill$\Box$
 \end{lem}
\indent In addition, the following auxiliary lemma will be needed for the proof of Theorem \ref{t5.1}.
\begin{lem} \label{l5.6}\rm
 For $\displaystyle\frac{n}{4}\leq\beta<1$, $n=2,3$, let $\displaystyle\frac{1}{p}+\frac{1}{2}=\frac{1}{q}+1$, $\displaystyle q=\frac{2s}{s-2}$, $\displaystyle s=\frac{ln}{n-l\beta}$ and $\displaystyle l>\frac{n}{3\beta-1}$. It follows that
 $$
1\leq p<\left\{
\begin{array}{ll}
  \frac{3}{2}~~\mbox{for}~~n=3,
 \\[0.2cm]
 2~~\mbox{for}~~ n=2.
 \end{array}
\right.
$$
 \end{lem}
{\bf Proof.} Direct calculation gives $\displaystyle p=\frac{ln}{ln-n+l\beta}$. Thanks to $\displaystyle l>\frac{n}{3\beta-1}$,  $\displaystyle\frac{n}{4}\leq\beta<1$ and $n=2,3$, it follows that $\displaystyle 1\leq p<\frac{n}{n-1}$. This completes the proof of this lemma.\hfill$\Box$
\\[0.2cm]
 \indent With the previous preparations, we begin to show Theorem \ref{t5.1}.
 \\[0.2cm]
{\bf Proof of Theorem \ref{t5.1}.}
 \\[0.3cm]
\indent We will work in a time interval with known regularity of the solutions to the Camassa-Holm equations with fractional Laplacian viscosity \eqref{1.1} and the imcompressible Navier-Stokes equations with fractional Laplacian viscosity \eqref{5.1}. Hence these are the unique solutions. Note that \eqref{1.1} and \eqref{5.1}, if $\mathbb{P}$ is the Leray projector onto the divergence free subspace of $L^2(\mathbb{R}^{n})$ and $\phi(t)$ is the fractional power heat kernel $\phi(t)=e^{-(-\Delta)^{\beta}t}$, then
 \begin{equation}\label{5.6}\mathbf{w}(t)=\phi(t)\ast \mathbf{w}_{0}-\int^{t}_{0}\phi(t-s)\ast \mathbb{P}\left[\mathbf{w}\cdot \nabla \mathbf{ w}\right](s)ds,\end{equation}
\begin{equation}\label{5.7}\mathbf{v}(t)=\phi(t)\ast\mathbf{ w}_{0}-\int^{t}_{0}\phi(t-s)\ast \mathbb{P}\left[\mathbf{u}\cdot \nabla \mathbf{v}-\mathbf{u}\cdot \nabla \mathbf{v}^{T}\right](s)ds.\end{equation}
 Thanks to \eqref{4.4}, a straightforward computation shows that
 \begin{equation}\label{5.8}
\left.
\begin{array}{ll}
\mathbf{w}(t)-\mathbf{v}(t)&=-\displaystyle\int^{t}_{0}\phi(t-s)\ast  \mathbb{P}\left[(\mathbf{w}-\mathbf{u})\cdot \nabla \mathbf{w}+\mathbf{u}\cdot \nabla(\mathbf{w}-\mathbf{v})\right] (s)ds
\\[0.3cm]
&\quad-\displaystyle\int^{t}_{0}\phi(t-s)\ast  \mathbb{P} \left[\sum\limits^{n}_{j=1}u_{j}\nabla(v_{j}-w_{j})
+\sum\limits^{n}_{j=1}(u_{j}-w_{j}) \nabla w_{j} \right ]  (s)ds.
 \end{array}
\right.
\end{equation}
 Note that the definition of the Leray projector, and the fact that the projector commutes with derivative for smooth functions in the entire space, using Young's inequality and Gagliardo-Nirenberg-Sobolev inequality, one deduces the following estimate for the first term of the integrand in \eqref{5.8}:
 \\[-0.3cm]
 \begin{equation}\label{5.9}
\left.
\begin{array}{ll}
& \left \|\phi(t-s)\ast  \mathbb{P} \left[(\mathbf{w}-\mathbf{u})\cdot \nabla \mathbf{w} \right](s)\right\| _{L^{q}(\mathbb{R}^{n})}
 \\[0.3cm]
    &\qquad\leq  \left \|\nabla\phi(t-s)\ast  \mathbb{P} \left[(\mathbf{w}-\mathbf{u})\cdot   \mathbf{w} \right](s)\right\| _{L^{q}(\mathbb{R}^{n})}
     \\[0.3cm]
    &\qquad\leq  \left \|\nabla\phi(t-s) \right\| _{L^{p}(\mathbb{R}^{n})} \left\| (\mathbf{w}-\mathbf{u})\cdot   \mathbf{w} \right\| _{L^{2}(\mathbb{R}^{n})}
     \\[0.3cm]
     &\qquad\leq   \left\|\nabla\phi(t-s) \right\| _{L^{p}(\mathbb{R}^{n})}\|  \mathbf{w}   \| _{L^{\frac{ln}{n-l\beta}}(\mathbb{R}^{n})} \left\|  \mathbf{w}-\mathbf{u}  \right\| _{L^{q}(\mathbb{R}^{n})}
      \\[0.3cm]
           &\qquad\leq   \left\|\nabla\phi(t-s) \right\| _{L^{p}(\mathbb{R}^{n})} \left\|  \Lambda^{\beta}\mathbf{w}   \right\| _{L^{l}(\mathbb{R}^{n})} \left\| \mathbf{ w}-\mathbf{u}  \right\| _{L^{q}(\mathbb{R}^{n})}.
      \end{array}
\right.
\end{equation}
Here and hereafter, $\displaystyle\frac{1}{q}+1=\frac{1}{p}+\frac{1}{2}$, $\displaystyle\frac{1}{2} =\frac{1}{q}+\frac{n-l\beta}{ln}$ and $\displaystyle 1-\frac{1}{p}=\frac{1}{2}-\frac{1}{q}=\frac{n-l\beta}{ln}<\frac{2\beta-1}{n}$ for $ \displaystyle l>\frac{n}{3\beta-1}$. \\
\indent Due to Proposition \ref{p5.3} with $\displaystyle\gamma=\frac{n}{2}\left(\frac{1}{2}-\frac{1}{q}\right)<\frac{\beta}{2}$, \eqref{5.9} can be bounded as follows:
\\[-0.3cm]
\begin{equation}\label{5.10}
\left.
\begin{array}{ll}
& \displaystyle \left\|\phi(t-s)\ast  \mathbb{P}\left[(\mathbf{w}-\mathbf{u})\cdot \nabla \mathbf{w} \right](s)\right\| _{L^{q}(\mathbb{R}^{n})}
\\[0.3cm]
          &\displaystyle\qquad\leq  \left\|\nabla\phi(t-s) \right\| _{L^{p}(\mathbb{R}^{n})}\left\|  \Lambda^{\beta}\mathbf{w}  \right \| _{L^{l}(\mathbb{R}^{n})}
          \\[0.3cm]
              &\displaystyle\qquad\quad\cdot
          \left(\| \mathbf{ w}-\mathbf{v}  \| _{L^{q}(\mathbb{R}^{n})}+C\alpha^{\frac{\beta}{2}-\gamma}\left\|\Lambda^{\beta}\mathbf{v} \right\| _{L^{2}(\mathbb{R}^{n})}\right).  \end{array}
\right.
\end{equation}
Making a similar derivation to \eqref{5.10} for the second term of \eqref{5.8}, one achieves
 \\[-0.3cm]
\begin{equation}\label{5.11}
\left.
\begin{array}{ll}
& \displaystyle \left\|\phi(t-s)\ast  \mathbb{P}\left[\mathbf{u}\cdot \nabla (\mathbf{w}-\mathbf{v})\right ](s)\right\| _{L^{q}(\mathbb{R}^{n})}
\\[0.3cm]
          &\displaystyle\qquad=\left\|\nabla\phi(t-s)\ast  \mathbb{P}\left[\mathbf{u}\cdot (\mathbf{w}-\mathbf{v}) \right](s)\right\| _{L^{q}(\mathbb{R}^{n})}
          \\[0.3cm]
          &\displaystyle\qquad\leq
            \left\|\nabla\phi(t-s) \right\| _{L^{p}(\mathbb{R}^{n})}\left\| \mathbf{u} \cdot (\mathbf{w}-\mathbf{v}) \right\| _{L^{2}(\mathbb{R}^{n})}
            \\[0.3cm]
          &\displaystyle\qquad\leq
           \left \|\nabla\phi(t-s) \right\| _{L^{p}(\mathbb{R}^{n})}\| \mathbf{u }  \| _{L^{\frac{ln}{n-l\beta}}(\mathbb{R}^{n})}\|  \mathbf{w}-\mathbf{v} \| _{L^{q}(\mathbb{R}^{n})}
            \\[0.3cm]
                &\displaystyle\qquad\leq\left\|\nabla\phi(t-s) \right\| _{L^{p}(\mathbb{R}^{n})}\left\|  \Lambda^{\beta}\mathbf{u}\right \| _{L^{l}(\mathbb{R}^{n})}\|   \mathbf{w}-\mathbf{v}  \| _{L^{q}(\mathbb{R}^{n})}.
\end{array}
\right.
\end{equation}
 In the same manner, one can deduce the following estimate for the third term of \eqref{5.8}:
 \begin{equation}\label{5.12}
\left.
\begin{array}{ll}
& \displaystyle \left\|\phi(t-s)\ast  \mathbb{P}\left[\sum\limits^{n}_{j=1}u_{j}\nabla(v_{j}-w_{j})\right]\right \| _{L^{q}(\mathbb{R}^{n})}
  \\[0.4cm]
 &\displaystyle\qquad\leq \left\|\nabla\phi(t-s)\ast  \mathbb{P}\left[\sum\limits^{n}_{j=1} u_{j} (v_{j}-w_{j})\right] \right\| _{L^{q}(\mathbb{R}^{n})}
   \\[0.4cm]
&\displaystyle\qquad\quad+\left\| \phi(t-s)\ast  \mathbb{P}\left[\sum\limits^{n}_{j=1}\nabla u_{j} (v_{j}-w_{j})\right] \right\| _{L^{q}(\mathbb{R}^{n})}
  \\[0.3cm]
 &\displaystyle\qquad\leq \left\|\nabla\phi(t-s)\right\| _{L^{p}(\mathbb{R}^{n})} \sum\limits^{n}_{j=1}  \left\|u_{j}  \right\| _{L^{\frac{ln}{n-l\beta}}(\mathbb{R}^{n})} \left\| v_{j}-w_{j}   \right\| _{L^{q}(\mathbb{R}^{n})}
   \\[0.3cm]
 &\displaystyle\qquad\quad+\left\| \phi(t-s)\right\| _{L^{p}(\mathbb{R}^{n})}\sum\limits^{n}_{j=1}  \left\| \nabla u_{j}\right\| _{L^{\frac{ln}{n-l\beta}}(\mathbb{R}^{n})} \left\| w_{j}-v_{j} \right \| _{L^{q}(\mathbb{R}^{n})}
   \\[0.3cm]
 &\displaystyle\qquad\leq \left\|\nabla\phi(t-s)\right\| _{L^{p}(\mathbb{R}^{n})} \sum\limits^{n}_{j=1}  \left\|\Lambda^{\beta}u_{j}\right\| _{L^{l}(\mathbb{R}^{n})}\left \| v_{j}-w_{j}\right \| _{L^{q}(\mathbb{R}^{n})}
   \\[0.3cm]
 &\displaystyle\qquad\quad+\left\| \phi(t-s)\right\| _{L^{p}(\mathbb{R}^{n})}\sum\limits^{n}_{j=1} \left \| \Lambda^{\beta}\nabla u_{j}\right\|_{L^{l}(\mathbb{R}^{n})}\left \| w_{j}-v_{j} \right \| _{L^{q}(\mathbb{R}^{n})}.
 \end{array}
\right.
\end{equation}
 The fourth term of \eqref{5.8} can also be bounded as
 \begin{equation}\label{5.13}
\left.
\begin{array}{ll}
& \displaystyle \left\|\phi(t-s)\ast  \mathbb{P}\left[\sum\limits^{n}_{j=1}(v_{j}-w_{j})\nabla w_{j}\right] \right\| _{L^{q}(\mathbb{R}^{n})}
  \\[0.4cm]
&\displaystyle\qquad\leq \left\|\nabla\phi(t-s)\ast  \mathbb{P}\left[\sum\limits^{n}_{j=1} (u_{j}-w_{j})w_{j} \right]\right \| _{L^{q}(\mathbb{R}^{n})}
  \\[0.4cm]
  &\displaystyle\qquad\quad+\left\| \phi(t-s)\ast  \mathbb{P}\left[\sum\limits^{n}_{j=1}\nabla  (u_{j}-w_{j})w_{j} \right] \right\| _{L^{q}(\mathbb{R}^{n})}
    \\[0.4cm]
  &\displaystyle\qquad\leq \left\|\nabla\phi(t-s)\ast  \mathbb{P}\left[\sum\limits^{n}_{j=1} (u_{j}-w_{j})w_{j} \right] \right\| _{L^{q}(\mathbb{R}^{n})}
    \\[0.4cm]
  &\displaystyle\qquad\leq\left \|\nabla\phi(t-s)\right \| _{L^{p}(\mathbb{R}^{n})}\left \| \sum\limits^{n}_{j=1} (u_{j}-w_{j})w_{j} \right\| _{L^{2}(\mathbb{R}^{n})}
    \\[0.4cm]
  &\displaystyle\qquad\leq \left\|\nabla\phi(t-s)\right \| _{L^{p}(\mathbb{R}^{n})}\left\|\mathbf{u}-\mathbf{w} \right\| _{L^{q}(\mathbb{R}^{n})} \left\| \mathbf{w}\right\| _{L^{\frac{ln}{n-l\beta}}(\mathbb{R}^{n})}
    \\[0.3cm]
  &\displaystyle\qquad\leq \left\|\nabla\phi(t-s)\right \| _{L^{p}(\mathbb{R}^{n})}\left\|\mathbf{u}-\mathbf{w} \right\| _{L^{q}(\mathbb{R}^{n})} \left\| \Lambda^{\beta}\mathbf{w}\right\| _{L^{l}(\mathbb{R}^{n})}
    \\[0.3cm]
  &\displaystyle\qquad\leq \left\|\nabla\phi(t-s) \right\| _{L^{p}(\mathbb{R}^{n})}\left\| \Lambda^{\beta}\mathbf{w}\right\| _{L^{l}(\mathbb{R}^{n})}
    \\[0.3cm]
  &\displaystyle\qquad\qquad\cdot\left(\left\|\mathbf{w}-\mathbf{v} \right\| _{L^{q}(\mathbb{R}^{n})} +C\alpha^{\frac{\beta}{2}-\gamma}\left\|\Lambda^{\beta}\mathbf{v} \right\| _{L^{2}(\mathbb{R}^{n})}\right).
.   \end{array}
\right.
\end{equation}
Putting together estimates \eqref{5.8}, \eqref{5.9}, \eqref{5.10}, \eqref{5.11}, \eqref{5.12} and \eqref{5.13}, we conclude
 \begin{equation}\label{5.14}
\left.
\begin{array}{ll}
& \displaystyle \left\| \mathbf{v}-\mathbf{w}\right \| _{L^{q}(\mathbb{R}^{n})}\leq C\alpha^{\frac{\beta}{2}-\gamma}\int^{t}_{0}\left\|\Lambda^{\beta}\mathbf{v} \right\| _{L^{2}(\mathbb{R}^{n})}ds
 \\[0.3cm]
 &\displaystyle\qquad\qquad\qquad +C\sup\limits_{t\in [0,T]}\left(\left\|\Lambda^{\beta}\mathbf{w} \right \| _{L^{l}(\mathbb{R}^{n})}+ \left\|\Lambda^{\beta}\mathbf{u} \right\| _{L^{l}(\mathbb{R}^{n})}+\left\|\Lambda^{\beta}\nabla \mathbf{u} \right\| _{L^{l}(\mathbb{R}^{n})}\right)
  \\[0.4cm]
 &\displaystyle\qquad\qquad\qquad\qquad \cdot\int^{t}_{0}\left(\left\| \phi(t-s)\right  \| _{L^{p}(\mathbb{R}^{n})}+\left\| \nabla\phi(t-s) \right \| _{L^{p}(\mathbb{R}^{n})}\right)\left\| \mathbf{v}-\mathbf{w}  \right\| _{L^{q}(\mathbb{R}^{n})}(s)ds.
    \end{array}
\right.
\end{equation}
Finally, thanks to Lemma \ref{l5.4}, Lemma \ref{l5.5} and Lemma \ref{l5.6}, we deduce that for $\phi(t)=e^{-(-\Delta)^{\beta}t}$
\begin{equation}\label{5.15}\left\|  \phi(t-s)  \right\| _{L^{p}(\mathbb{R}^{n})}\leq\frac{1}{(t-s)^{\delta_{1}}},
~~\delta_{1}=\left(1-\frac{1}{p}\right)\frac{n}{2\beta},\end{equation}
\begin{equation}\label{5.16}\left\|  \nabla\phi(t-s) \right \| _{L^{p}(\mathbb{R}^{n})}\leq\frac{1}{(t-s)^{\delta_{2}}},
~~\delta_{2}=\frac{1}{2\beta}+\left(1-\frac{1}{p}\right)\frac{n}{2\beta}.\end{equation}
 As a consequence, for $\displaystyle\gamma=\frac{n}{2}\left(\frac{1}{p}-\frac{1}{q}\right)<\frac{\beta}{2}$, we infer from \eqref{5.14}, \eqref{5.15} and \eqref{5.16} that
\begin{equation}\label{5.17}
\| \mathbf{v}-\mathbf{w} \| _{L^{q}(\mathbb{R}^{n})}\leq C\alpha^{\frac{\beta}{2}-\gamma}+B\int^{t}_{0}\frac{1}{(t-s)^{\delta}}\| \mathbf{v}-\mathbf{w}  \| _{L^{q}(\mathbb{R}^{n})}(s)ds,\end{equation}
 where
 \begin{equation}\label{5.18}
\left\{
\begin{array}{ll}
 \displaystyle A=C\int^{t}_{0}\left\|\Lambda^{\beta}\mathbf{v} \right\| _{L^{2}(\mathbb{R}^{n})}ds,
 \\[0.3cm]
 \displaystyle B=C\sup\limits_{t\in [0,T]}\left(\left\|\Lambda^{\beta}\mathbf{w} \right\| _{L^{l}(\mathbb{R}^{n})}+\left\|\Lambda^{\beta}\mathbf{u} \right\| _{L^{l}(\mathbb{R}^{n})}+\left\|\Lambda^{\beta}\nabla \mathbf{u}\right \| _{L^{l}(\mathbb{R}^{n})}\right),
 \\[0.3cm]
 \displaystyle \delta=\max\left\{\delta_{1},\delta_{2}\right\}=\frac{1}{2\beta}+\left(1-\frac{1}{p}\right)\frac{n}{2\beta}<\frac{1}{2\beta}+ \frac{2\beta-1}{n}\frac{n}{2\beta}<1.
\end{array}
\right.
\end{equation}
 Here, we have used some known facts:
 $$l>\frac{n}{3\beta-1},~\frac{1}{q}+1=\frac{1}{p}+\frac{1}{2},~~\frac{1}{2} =\frac{1}{q}+\frac{n-l\beta}{ln},~~
 1-\frac{1}{p}=\frac{1}{2}-\frac{1}{q}=\frac{n-l\beta}{ln}<\frac{2\beta-1}{n}.$$
  Note that \eqref{5.17} and \eqref{5.18}, the Gronwall inequality then implies that
\begin{equation}\label{5.19}\displaystyle\| \mathbf{v}-\mathbf{w} \| _{L^{q}(\mathbb{R}^{n})}\leq A\alpha^{\frac{\beta}{2}-\gamma}\exp\left( \int^{t}_{0}\frac{B}{(t-s)^{\delta}} ds \right),\end{equation}
  where $\displaystyle \int^{t}_{0} (t-s)^{-\delta}ds=\left.\frac{(t-s)^{-\delta+1}}{-\delta+1}\right|^{t}_{0}=
  -\frac{t^{-\delta+1}}{-\delta+1}$, which is finite for $\delta<1$ and $t\in [0,T]$. Hence, for $t\in [0,T]$, letting $\alpha\rightarrow 0$ deduces that $\mathbf{v}\rightarrow \mathbf{w}$ strongly in $L^{q}(\mathbb{R}^{n})$. \\
  \indent This finishes the proof of Theorem \ref{t5.1}.\hfill$\Box$
\section*{Acknowledgments }
The authors would like to thank Professor Dong Li and Dr. Yuan Cai for their helpful comments. Zaihui Gan is partially supported by the National Natural Science Foundation
of China ( No. 11571254).

\end{document}